\documentclass[11pt]{article}

\usepackage{amsfonts,amssymb,amsmath,mathrsfs}
\usepackage[hidelinks,colorlinks=true,linkcolor=blue,citecolor=blue]{hyperref}

\usepackage{xcolor,latexsym,amsfonts}

\numberwithin{equation}{section}

\newtheorem{theorem}{Theorem}[section]
\newtheorem{lemma}[theorem]{Lemma}

\newtheorem{remark}{Remark}

\newtheorem{corollary}[theorem]{Corollary}
\newtheorem{proposition}[theorem]{Proposition}
\newtheorem{example}[theorem]{Example}

\def\<{\langle}\def\>{\rangle}

\def\beqnn{\begin{eqnarray*}}\def\eeqnn{\end{eqnarray*}}
\def\<{\langle}\def\>{\rangle}

\def\beqlb{\begin{eqnarray}}\def\eeqlb{\end{eqnarray}}

\begin{document}
	
	\bigskip\bigskip
	\noindent{\Large\bf Quasi-stationary distributions for subcritical branching Markov chains\footnote{
			\noindent  This work was supported in part by NSFC (No.~11971062) and  the National
			Key Research and Development Program of China (No.~2020YFA0712900).} }
	
	\noindent{
		Wenming Hong\footnote{ School of Mathematical Sciences
			\& Laboratory of Mathematics and Complex Systems, Beijing Normal
			University, Beijing 100875, P.R. China. Email: wmhong@bnu.edu.cn} ~
		Dan Yao\footnote{ School of Mathematical Sciences
			\& Laboratory of Mathematics and Complex Systems, Beijing Normal
			University, Beijing 100875, P.R. China. Email: dyao@mail.bnu.edu.cn}~

	}
	
	\begin{center}
		\begin{minipage}{12cm}
			\begin{center}\textbf{Abstract}\end{center}
			\footnotesize
			Consider a subcritical branching Markov chain. Let $Z_n$ denote the counting measure of particles of generation $n$. Under some conditions, we give a probabilistic proof for the existence of the Yaglom limit of $(Z_n)_{n\in\mathbb{N}}$  by the moment method, based on the spinal decomposition and the many-to-few formula. As a result, we give explicit integral representations of all quasi-stationary distributions of $(Z_n)_{n\in\mathbb{N}}$, whose proofs are direct and probabilistic, and don't rely on Martin boundary theory.
			\bigskip
			
			\mbox{}\textbf{Keywords:} branching Markov chain, Yaglom limit, quasi-stationary distribution, conditional reduced Galton-Watson process. \\
			\mbox{}\textbf{Mathematics Subject Classification}:  Primary 60J80; 60F05;
			secondary 60F10.
			
		\end{minipage}
	\end{center}
	\section{Backgrounds}
	\label{sec1}
	Let $\mathbb{N}:=\{0,1,2,\cdots\}$ and $\mathbb{N}_+:=\mathbb{N}\backslash\{0\}$. Let $p=\{p_j\}_{j\in\mathbb{N}}$ be a probability measure on $\mathbb{N}$ and $m:=\sum jp_j$ be the mean of $p$. Suppose that $N=(N_n:n\in\mathbb{N})$ is a subcritical Galton-Watson process (GW-process) with offspring distribution $p$ (i.e., $m<1$). We denote by $\mathbf{P}_k$ the law of the GW-process $N$ starting from $N_0=k$. It is well known from \cite{AN72,G99,HS67,J67,Y47} that 
	\begin{align*}
		\lim_{n\rightarrow\infty}\mathbf{P}_k(N_n=j\mid N_n>0)=\nu_{\text{min}}(j),\quad j\in\mathbb{N}_+
	\end{align*}
	exists and $\nu_{\text{min}}=\{\nu_{\text{min}}(j)\}_{j\in\mathbb{N}_+}$ is a probability measure on $\mathbb{N}_+$, which is called the {\it Yaglom limit} of $N$.
	
	A probability measure $\nu=\{\nu(j)\}_{j\in\mathbb{N}_+}$ on $\mathbb{N}_+$ is called a {\it quasi-stationary distribution} (QSD) of $N$ if $\nu$ satisfies
	\begin{align}\label{1.1}
		\mathbf{P}_{\nu}(N_n=j\mid N_n>0)=\nu(j),\quad j\in\mathbb{N}_+,\quad n\in\mathbb{N},
	\end{align}
	where $\mathbf{P}_{\nu}(\cdot):=\sum_{k=1}^{\infty}\nu(k)\mathbf{P}_k(\cdot)$. If $f(s)$ is the generating function of the offspring distribution and $g(s)$ is the generating function of the probability measure $\nu$, then (\ref{1.1}) is equivalent to 
	\begin{align}\label{1.2}
			g\big(f(s)\big)-g\big(f(0)\big)=m^{\alpha} g(s),\quad s\in[0,1],
	\end{align}
	where $\alpha\in(0,\infty)$. Indeed, we have (\ref{1.2}) $\Rightarrow$ (\ref{1.1}) by easy calculations. (\ref{1.1}) $\Rightarrow$ (\ref{1.2}) follows from the fact that if $\nu$ is a QSD of $N$, there exists $\alpha\in(0,\infty)$ such that $\mathbf{P}_{\nu}(N_n>0)=m^{\alpha n}$ for all $n\in\mathbb{N}$; see, e.g.,  \cite[Proposition 2]{MV12}. In this case, the probability measure $\nu$ is called a QSD of $N$ with eigenvalue $m^{\alpha}$. 
	
	There is a large literature on invariant measures for subcritical GW-processes. An explicit characterization of $1$-invariant measures was given by Kesten and Spitzer \cite{KS67}, motivated by the Martin boundary theory. Hoppe \cite{H80} used functional equations to show that there exists a bijection between QSDs with eigenvalue $m^{\alpha}$ (where $\alpha\in(0,1]$) and $1$-invariant measure. Recently, Maillard \cite{M18} provided a new probabilistic approach to $m^{\alpha}$-invariant measures (and QSDs) for subcritical GW-processes, yielding complete characterizations of these measures involving explicit formulas.
	
	
	
	In this paper, we are interested in branching  Markov chains, which can be regarded as the Markov chain indexed by a discrete tree.
	As we know, the branching  Markov chain can be regarded as a multitype branching process with an infinite set of types. Hering \cite{H77} used an analytical method to obtain the Yaglom's theorem for subcritical continuous-time Markov branching processes, whose proof was relying on the mean operator and its spectral properties. Based on the Yaglom's theorem, Asmussen and Hering \cite{AH83} gave the explicit integral representation of the $1$-invariant measure for subcritical continuous-time Markov branching processes by using Martin boundary theory, and they also proved that there exists a one-to-one correspondence between QSDs and $1$-invariant measures. In particular, Asmussen and Hering \cite{AH83} claimed that the above results about the Yaglom's theorem and invariant measures remain valid for branching Markov chains.
	
	
	Our object in this article is to use a probabilistic method to deal with the above two topics. Firstly, based on the conditional reduced GW-process,
	the spinal decomposition and the many-to-few formula, we will prove the Yaglom's theorem for subcritical branching Markov chains by the moment method. Secondly, we give explicit integral representations of all quasi-stationary distributions of subcritical branching Markov chains (similar to the proof of subcritical GW-processes in Maillard \cite{M18}). Recently, the moment method have been used in the proof of  the conditional central limit theorem for the critical branching random walk (see \cite{HL23}).

	The related results of the Yaglom's theorem and QSDs for other branching models have been extensively studied in the literature, see, for example, \cite{EP90,H22} for the Yaglom's theorem and \cite{V23} for the cluster-invariant point processes of critical branching spatial processes on $\mathbb{R}^{d}$. In addition, Lambert \cite{L07} proved that the Yaglom limit exists and characterized all QSDs for subcritical continuous-state branching processes by their Laplace transforms. For subcritical superprocesses, Liu et al. \cite{LR21} gave the existence of the Yaglom limit and identified all QSDs under some conditions for cumulant semigroups. 
	
	
	
	{\it Notations} \  Let $\mathcal{B}(\mathbb{R})$ (resp. $\mathcal{B}((0,\infty))$) denote the $\sigma$-algebra on the field of real numbers $\mathbb{R}$ (resp. $(0,\infty)$) generated by the class of open sets. 
	We write $\|\cdot\|$ for the supremum/uniform norm. A {\it point measure} is an integer-valued measure on $\mathbb{R}$ and is finite on all compact sets of $\mathbb{R}$. Let $\mathfrak{N}(\mathbb{R})$ and $\mathfrak{M}(\mathbb{R})$ denote respectively the space of finite point measures and finite Borel measures on $\mathbb{R}$ equipped with the topology of weak convergence. Let $\mathbf{0}$ denote the null measure on $\mathbb{R}$. We write $\mathfrak{N}^{\circ}(\mathbb{R}):=\mathfrak{N}(\mathbb{R}) \backslash\{\mathbf{0}\}$ and $\mathfrak{M}^{\circ}(\mathbb{R}):=\mathfrak{M}(\mathbb{R}) \backslash\{\mathbf{0}\}$. Let $\mathcal{B}(\mathfrak{N}(\mathbb{R}))$, $\mathcal{B}(\mathfrak{N}^{\circ}(\mathbb{R}))$, $\mathcal{B}(\mathfrak{M}(\mathbb{R}))$ and $\mathcal{B}(\mathfrak{M}^{\circ}(\mathbb{R}))$ be the Borel $\sigma$-algebra generated on $\mathfrak{N}(\mathbb{R})$, $\mathfrak{N}^{\circ}(\mathbb{R})$, $\mathfrak{M}(\mathbb{R})$ and $\mathfrak{M}^{\circ}(\mathbb{R})$ by this weak topology, respectively. We use $\mu(h)$ to denote the integral of a function $h$ with respect to a measure $\mu$. 
	
	
	\section{Introduction and main results}
	\label{sec2}
	\subsection{Definition of the model}
	\label{sec2.1}
	We use the classical Ulam-Harris-Neveu notation for discrete trees. Each individual
	in the population is an element of 
	$$
	\mathcal{U}=\{\varnothing\}\cup\bigcup_{n\geq 1}(\mathbb{N}_+)^n,
	$$
	which is written as a finite string $u_1\cdots u_n$ of positive integers. As before, for $u=u_1\cdots u_n\in\mathcal{U}$, we denote by $|u|=n$ the generation of the individual $u$ (with $u_1\cdots u_n=\varnothing$ for $n=0$). If $u\neq \varnothing$, we write  $\overleftarrow{u}=u_1\cdots u_{n-1}$ for its parent. For elements $u=u_1\cdots u_n$ and $v=v_1\cdots v_l$ of $\mathcal{U}$, let $uv=u_1\cdots u_nv_1\cdots v_l$ be the concatenation of $u$ and $v$ with the convention $u\varnothing=\varnothing u=u$. We write $v\prec u$ if $v$ is an ancestor of $u$, i.e., $\exists\ \bar{v}\in\mathcal{U}$ such that $u=v\bar{v}$.
	
	Let $p=\{p_j\}_{j\in\mathbb{N}}$ be a probability measure on $\mathbb{N}$ called the offspring distribution. Let us consider a family of independent, identically distributed integer-valued random variables $\{N(u);u\in\mathcal{U}\}$ with the same law $p$. We say a random subset $\mathbb{T}\subset\mathcal{U}$ is a Galton-Watson tree (GW-tree) with offspring distribution $p$ rooted at $\varnothing$ if it is defined by
	$$
	\mathbb{T}:=\{u=u_1\cdots u_n\in\mathcal{U}:n\geq0,1\leq u_j\leq N_{u_1\cdots u_{j-1}} \text{ for every }1\leq j\leq n\}.
	$$
	We write $\mathbb{G}_n$ for the set of individuals in the $n$-th generation. 
	
	Let $P$ be the transition probability on $\mathbb{R}\times\mathcal{B}(\mathbb{R})$. Let $V(u)\in\mathbb{R}$ be the position of the individual $u$. The {\it branching Markov chain} $(V(u):u\in\mathbb{G}_n,n\geq0)$ starting from one single individual at $x\in\mathbb{R}$ is defined recursively as follows under the probability $\mathbb{P}_{\delta_{x}}$: 
	\begin{itemize}
		\item[(1)] At time $0$, there is one individual $\varnothing$ located at $x$, which forms the $0$-th generation, i.e., $\mathbb{G}_0=\{\varnothing\}$ and $V(\varnothing)=x$.
		\item[(2)] For each $n\geq0$ and $u\in\mathbb{G}_n$, the individual $u$ (after one unit time) splits independently of the other individuals of the $n$-th generation into a random number $N(u)$ of children with law $p$, which implies $\mathbb{G}_{n+1}=\{uj:u\in\mathbb{G}_n,1\leq j\leq N(u)\}$. Then these children {are positioned independently 
			according to the one-step transition probability} $P(V(u),\cdot)$, i.e., 
		\begin{align*}
			\mathbb{P}_{\delta_{x}}(V(u1)\in\mathrm{d}x_1,\cdots,V(uk)\in\mathrm{d}x_k\mid N(u)=k,V(u)=y)=\prod_{i=1}^{k}P(y,\mathrm{d}x_i),\  x_1,\cdots,x_k\in\mathbb{R}.
		\end{align*}
	\end{itemize}
	Throughout this paper,  We denote by $\mathbb{P}:=\mathbb{P}_{{\delta_0}}$ the probability measure for this system started from a single particle at $0$, and the corresponding expectation is $\mathbb{E}$.
	
	Furthermore, we can define a branching Markov chain starting from a point measure on $\mathbb{R}$. For any $\mu\in\mathfrak{N}^{\circ}(\mathbb{R})$, then it can be written as
	$\mu=\sum_{i\in I}\delta_{x_i}$, where $I\subseteq\mathbb{N}_+$ and $x_i\in\mathbb{R}$ ($x_i$ may take the same value for the different choices of $i$). Then the branching Markov chain starting from a point measure $\mu$ is defined as follows under the probability $\mathbb{P}_{\mu}$: 
	\begin{itemize}
		\item[(1)] At time $0$, there are $\mu(\mathbb{R})$ particles with positions $x_1,\cdots,x_{\mu(\mathbb{R})}$ respectively. These particles and their positions form the $0$-th generation.
		\item[(2)] For every particle {of generation $0$} located at $x_i\in\mathbb{R}$, it evolves as a branching Markov chain with offspring distribution $p$ and one-step transition probability $P$ starting from one particle located at $x_i$ according to the previous construction.
		\item[(3)] These branching Markov chains starting from one particle located at $x_i$ are independent of each other for the different choices of $i$. 
	\end{itemize}
	
	The branching Markov chain will be studied in the framework of measure-valued Markov processes. To this end, let $Z_n(\cdot)$ be the counting measure of particles in generation $n$. Then the process $Z=(Z_n:n\in\mathbb{N})$ is a $\mathfrak{N}(\mathbb{R})$-valued Markov process. For each $B\in\mathcal{B}(\mathbb{R})$, $Z_n(B)$ denotes the number of individuals in the $n$-th generation located in $B$. By convention, we write $N_n:=Z_n(\mathbb{R})$. Clearly, under $\mathbb{P}_{\mu}$, $(N_n:n\in\mathbb{N})$ is a GW-process with offspring distribution $p$ starting from $N_0=\mu(\mathbb{R})$, 
	which is called the {\it associated GW-process}. {Let} $m=\sum_j j p_j$ and $f$ denote respectively the mean and the generating function of $p$. We write $f_n$ for the generating function of $N_n$, then $f_1=f$. For any $h\in\mathcal{B}(\mathbb{R})$ we write $Z_n(h):=\int_{\mathbb{R}} h(x) Z_n(\mathrm{d}x)$.
	
	In this paper, we always assume $0<m<1$, which guarantees that the associated GW-process is subcritical and non-trivial (i.e., $p_0<1$). 
	\subsection{ A warm-up about finite point processes}
	\label{sub2.2}
	We recall some properties about finite point processes. One can refer to \cite{AH83,DV03,DV08} for more details. We first give some function classes. There is no ambiguity to use $\mathcal{B}(\mathbb{R})$ to denote the set of real-valued Borel measurable functions on $\mathbb{R}$.
	\begin{center}
		\begin{tabular}{rl}
			$C(\mathbb{R})$ :& the set of continuous real functions on $\mathbb{R}$.\\
			$\mathcal{B}_b^+(\mathbb{R})$ :& the set of bounded positive Borel measurable functions on $\mathbb{R}$.\\
			$\mathcal{V}(\mathbb{R})$ :& $\{h\in\mathcal{B}(\mathbb{R}):0\leq h(x)\leq 1$ \text{for every} $x\in\mathbb{R}\}$.\\
			$\mathcal{V}_0(\mathbb{R})$ :& $\{h\in\mathcal{B}(\mathbb{R}):\inf_{x\in\mathbb{R}}h(x)>0,\ h(x)\leq 1$ \text{for every} $x\in\mathbb{R}\}$.\\
			$\mathcal{C}_0(\mathbb{R})$ :& $\{h\in C(\mathbb{R}):\inf_{x\in\mathbb{R}}h(x)>0,\ h(x)\leq1$ \text{for every} $x\in\mathbb{R}\}$.
		\end{tabular}
	\end{center}
	
	A {\it finite point process} on $\mathbb{R}$ is a random variable taking values in $\mathfrak{N}(\mathbb{R})$. We will now introduce the probability generating functional of a point process, which is an important tool in handling point processes. Suppose that $\{X;\mathbb{P}\}$ is a finite point process on $\mathbb{R}$ throughout this section. 
	The {\it probability generating functional} (p.g.fl.) $G_{X}$ of $X$ is defined by
	\begin{align*}
		G_{X}[h]:=\mathbb{E}\big[\mathrm{e}^{X(\log h)}\big],\quad h\in\mathcal{V}(\mathbb{R}),
	\end{align*}
	where $X(\log h):=\int_{\mathbb{R}} \log h(x) X(\mathrm{d}x)$.
	\begin{remark}
		Daley and Vere-Jones gave the definition of the p.g.fl. for more general point processes (see \cite{DV08}, p.59). For finite point processes, we remove the requirement ``$1-h$ vanishes outside some bounded set'', whose reason is that the total number of points is finite almost surely in this case. 
	\end{remark}
	
	For a probability measure $\mathbf{Q}$ on $\mathfrak{N}(\mathbb{R})$,  {its {\it p.g.fl.} $G_{\mathbf{Q}}$} and {\it Laplace functional} $L_{\mathbf{Q}}$ are defined as follows
	\begin{align*}
		&G_{\mathbf{Q}}[h]:=\int_{\mathfrak{N}(\mathbb{R})}\mathrm{e}^{\mu(\log h)} \mathbf{Q}(\mathrm{d}\mu),\quad h\in\mathcal{V}(\mathbb{R}).\\
		&L_{\mathbf{Q}}[h] := \int_{\mathfrak{N}(\mathbb{R})} \mathrm{e}^{-\mu(h)} \mathbf{Q}(\mathrm{d}\mu),\quad h\in\mathcal{B}_b^+(\mathbb{R}).
	\end{align*}
	Since $-\log h\in\mathcal{B}_b^+(\mathbb{R})$ if and only if $h\in\mathcal{V}_0(\mathbb{R})$ (see \cite{DV08}, p.59), we get the following lemma, which
	characterizes the relationship between the p.g.fl. and the Laplace functional.
	\begin{lemma}\label{lem2.1}
		For any $\mathbf{Q}\in\mathfrak{N}(\mathbb{R})$, suppose that $G_{\mathbf{Q}}$ and $L_{\mathbf{Q}}$ denote respectively the the p.g.fl. and the Laplace functional of $\mathbf{Q}$. Then we have
		\begin{align*}
			G_{\mathbf{Q}}[h]=\int_{\mathfrak{N}(\mathbb{R})} \mathrm{e}^{\mu(\log h)} \mathbf{Q}(\mathrm{d}\mu) =L_{\mathbf{Q}}[-\log h],\quad h\in\mathcal{V}_0(\mathbb{R}).
		\end{align*}
	\end{lemma}
	
	{Based on this result, we only need to restrict the domain of the p.g.fl. to $\mathcal{V}_0(\mathbb{R})$ in the remainder of this paper. }
	By a combination of Lemma \ref{lem2.1} and \cite[Theorem 1,17]{L11}, we obtain the following result. 
	\begin{lemma}\label{lem2.2}
		Suppose that $\mathbf{Q}_1,\mathbf{Q}_2$ are two probability measures on $\mathfrak{N}(\mathbb{R})$. Then $\mathbf{Q}_1=\mathbf{Q}_2$ if and only if $G_{\mathbf{Q}_1}[h]=G_{\mathbf{Q}_2}[h]$ for all $h\in\mathcal{C}_0(\mathbb{R})$.
	\end{lemma}
	
	The following lemma characterizes the weak convergence of finite point processes in terms of their probability generating functionals.
	\begin{lemma}\label{lem2.4}{\rm(\cite[Proposition 11.1.VIII.]{DV08})}
		Suppose that $\mathbf{Q}_1,\cdots,\mathbf{Q}_n$ and $\mathbf{Q}$ are probability measures on $\mathfrak{N}(\mathbb{R})$. Then $\mathbf{Q}_n\rightarrow\mathbf{Q}$ weakly if and only if 
		$G_{\mathbf{Q}_n}[h]\rightarrow G_{\mathbf{Q}}[h]$ for all $h\in\mathcal{V}_0(\mathbb{R})$.
	\end{lemma}
	
	{In addition, we give some characterizations of the weak convergence of finite point processes in terms of their Laplace functions. }
	
	\begin{lemma}\label{lem2.5}{\rm (\cite[Theorem 1.18]{L11})}
		Suppose that $\mathbf{Q}_1,\cdots,\mathbf{Q}_n$ and $\mathbf{Q}$ are probability measures on $\mathfrak{N}(\mathbb{R})$. Then $\mathbf{Q}_n\rightarrow\mathbf{Q}$ weakly if and only if 
		$L_{\mathbf{Q}_n}[h]\rightarrow L_{\mathbf{Q}}[h]$ for all $h\in\mathcal{B}_b^+(\mathbb{R})$.
	\end{lemma}
	
	{One can refer to \cite[Lemma A.4.]{LR21} and \cite[Theorem 1.20]{L11} for} the continuity theorem for Laplace functionals of finite random measures. From \cite[Proposition 3.14]{R87}, it is easy to obtain that when the weak limit of a sequence of point measures is finite, it is still a point measure. Hence by a similar proof of \cite[Theorem 1.20]{L11} we obtain:
	
	\begin{lemma}\label{lem2.3}{\rm (The continuity theorem)}
		Let $(\mathbf{Q}_n)_{n\in\mathbb{N}}$ be a sequence of probability measures on $\mathfrak{N}(\mathbb{R})$. Suppose that
		\begin{itemize}
			\item[{\rm (1)}] for each $h\in\mathcal{B}_b^+(\mathbb{R})$, the limit $L[h]=\lim_{n\rightarrow\infty} L_{\mathbf{Q}_n}[h]$ exists;
			\item[{\rm (2)}] for each $h_n\rightarrow h$ in $\mathcal{B}_b^+(\mathbb{R})$ and $\sup_n \|h_n\|<\infty$, we have $L[h_n]\rightarrow L[h]$ as $n\rightarrow\infty$.  
		\end{itemize}
		Then there exists a unique probability measure $\mathbf{Q}$ on $\mathfrak{N}(\mathbb{R})$ such that $\lim_{n\rightarrow\infty}\mathbf{Q}_n=\mathbf{Q}$ by weak convergence and $L_{\mathbf{Q}}=L$ on $\mathcal{B}_b^+(\mathbb{R})$. 
	\end{lemma}

	Let us go back to the measure-valued process $Z=(Z_n:n\in\mathbb{N})$. Let $G_n$ be the p.g.fl. of $Z_n$ under $\mathbb{P}_{\delta_{x}}$, that is, 
	\begin{align}\label{4}
		G_n[h](x)=\mathbb{E}_{\delta_{x}}\big[\mathrm{e}^{Z_n(\log h)}\big],\quad x\in\mathbb{R},\quad h\in\mathcal{V}_0(\mathbb{R}).
	\end{align}
	In particular, we set $G[h]:= G_1[h]$ for simplicity. 
	\begin{proposition}\label{prop1}{\rm (\cite[Example 5.5(c)]{DV03})}
		Suppose that $G_n$ is the p.g.fl. of $Z_n$ under $\mathbb{P}_{\delta_{x}}$. Let $G_0$ be the identity operator and let $G:=G_1$. Then $(G_n)_{n\in\mathbb{N}}$ has the semigroup property, which means that 
		\begin{align}\label{2}
			G_{n+1}[h](x)=G_n\big[G[h](\cdot)\big](x),\quad x\in\mathbb{R},\quad h\in\mathcal{V}_0(\mathbb{R}).
		\end{align}
		In this case, we say $(G_n)_{n\in\mathbb{N}}$ is the generating semigroup.
	\end{proposition}
	\begin{remark}
		Let $G_n[G[h](\cdot)](x)$ be defined by (\ref{4}) with the function $h$ replaced by the function $y\mapsto G[h](y)$, where $G[h](y)=G_1[h](y)=\mathbb{E}_{\delta_y}[\mathrm{e}^{Z_1(\log h)}]$. In fact, the definition of $G_n[G[h](\cdot)](x)$ is valid. Since the assumption $m<1$ implies that $\inf_{y\in\mathbb{R}} G_n[h](y)\geq \mathbb{P}(N_n=0) \geq p_0>0$, the p.g.fl. $G_n$  is a mapping from $\mathcal{V}_0(\mathbb{R})$ to $\mathcal{V}_0(\mathbb{R})$. 
	\end{remark}
	\begin{remark}
		If we take $h\equiv s$, then (\ref{2}) can be written as
		\begin{align}\label{b2}
			f_n(s)=f\big(f_{n-1}(s)\big),\quad f_0(s)=s,\quad f_1(s)=f(s),\quad s\in[0,1] .
		\end{align}
		which describes the iteration relation for the generating functions $f_n$ of the numbers of particles in generation $n$ of the GW-process.
	\end{remark}
	\subsection{Main results}
	\label{sub2.3}
	Recall that $P$ is the one-step transition probability, then its $n$-step transition probability $P_n$ is given by the iterate 
	\begin{align*}
		P_1:=P;\quad P_n(x,\cdot):=\int_{\mathbb{R}} P_{n-1} (y,\cdot) P(x,\mathrm{d}y),\quad n\geq2.
	\end{align*}
	We write $P_nh(x)$ for the integral of a function $h$ with respect to the measure $P_n(x,\cdot)$. For the convenience of statement of the results, we formulate the following conditions:
	\begin{itemize}
		\item[(C1):] The offspring distribution $p$ satisfies $0<m<1$ and has finite moments of all orders.
		\item[(C2):] $\varlimsup_{n\rightarrow\infty}(p_n)^{\frac{1}{n}}\in[0,1).$
		\item[(C3):] For any $x\in\mathbb{R}$ and any $h\in\mathcal{B}_b^+(\mathbb{R})$, there exists a probability measure $\pi$ on $\mathbb{R}$ such that $P_nh(x)\rightarrow \pi(h)$ as $n\rightarrow\infty$. 
		\item[(C4):] For any $x\in\mathbb{R}$ and any $h\in\mathcal{B}_b^+(\mathbb{R})$, there exists a probability measure $\pi$ on $\mathbb{R}$ such that $P_nh(x)\rightarrow \pi(h)$ uniformly in $x$ as $n\rightarrow\infty$.
	\end{itemize}

	The conditions {\rm(C1)} and {\rm(C2)} are mainly concerned with the branching mechanism.
	
	\begin{example} In the following cases, the offspring distribution $p$ satisfies {\rm(C2)} and has finite moments of all order. 
		
		{\rm 1)} $p$ has a bounded support, that is, $\sum_{0\leq j\leq j_0}p_j=1$ for some $j_0\in\mathbb{N}$.
		
		{\rm 2)} $p$ is a Poisson distribution with parameter $\lambda\in[0,1)$.
		
		{\rm 3)} $f$ is a linear fractional offspring generating function given by
		\begin{align*}
			1-f(s)= \frac{a}{1-b} - \frac{as}{1-bs},\quad s\in[0,1],
		\end{align*}
		where $a\geq0,b\in[0,1)$ and $a+b\leq1$.
	\end{example}
	
	\begin{remark}
		The condition {\rm(C3)} holds if the Markov chain with transition probability $P$ is ergodic, while the strongly ergodic of this Markov chain is the sufficient condition for {\rm(C4)}. 
	\end{remark}
	
	The symbol $\stackrel{\text{w}}{\longrightarrow}$ stands for weak convergence throughout this paper. 
	Our first main result is described as follows. 
	\begin{theorem}\label{Yaglom}{\rm (Yaglom's Theorem)}
		Suppose that the conditions {\rm (C1)}-{\rm (C3)} hold. Then there exists a unique probability measure $\mathbf{Q}_{\rm min}$ on $\mathfrak{N}^{\circ}(\mathbb{R})$ such that
		\begin{align*}
			\mathbb{P}_{\mu}(Z_n\in\cdot\mid N_n>0)\xrightarrow[n \rightarrow \infty]{\text{\rm w}} \mathbf{Q}_{\rm min}(\cdot), \quad \mu \in \mathfrak{N}^{\circ}(\mathbb{R}).
		\end{align*}
		In this case, the probability measure $\mathbf{Q}_{\rm min}$ is called the Yaglom limit of $Z$.
	\end{theorem}
	
	More generally, it is meaningful to condition the branching Markov chain on survival up to time $n+l$, and we get the following results.
	\begin{corollary}\label{cor2.7}
		Suppose that the conditions {\rm (C1)}-{\rm (C3)} hold. For fixed $l\in\mathbb{N}$, there exists a unique probability measure $\mathbf{Q}_{(l)}$ on $\mathfrak{N}^{\circ}(\mathbb{R})$ such that
		\begin{align*}
			\mathbb{P}_{\mu}(Z_n\in\cdot\mid N_{n+l}>0)\xrightarrow[n \rightarrow \infty]{\text{\rm w}} \mathbf{Q}_{(l)}(\cdot), \quad \mu \in \mathfrak{N}^{\circ}(\mathbb{R}).
		\end{align*} 
		Moreover, let $L_{\mathbf{Q}_{\rm min}}$ be the Laplace functional of $\mathbf{Q}_{\rm min}$, then the Laplace functional of the weak limit is given by 
		\begin{align}\label{2.7}
			L_{\mathbf{Q}_{(l)}} [h] = \frac{ L_{\mathbf{Q}_{\rm min}}[h]- L_{\mathbf{Q}_{\rm min}}\big[h-\log f_l(0)\big] }{m^{l} },\quad h\in\mathcal{B}_b^+(\mathbb{R}).
		\end{align}
	\end{corollary}
	
	\begin{corollary}\label{cor2.8}
		Suppose that the conditions {\rm (C1)}-{\rm (C3)} hold. 
		\begin{itemize}
			\item [{\rm(1)}] For any $\mu \in \mathfrak{N}^{\circ}(\mathbb{R})$, there exists a unique probability measure $\mathbf{Q}_{\infty}$ on $\mathfrak{N}^{\circ}(\mathbb{R})$ such that $\lim_{n\rightarrow\infty} \lim_{l\rightarrow\infty} \mathbb{P}_{\mu}(Z_n\in\cdot\mid N_{n+l}>0)= \mathbf{Q}_{\infty}(\cdot)$ by weak convergence and its Laplace functional satisfies 
			\begin{align*}
				L_{\mathbf{Q}_{\infty}} [h] = \varphi(0)\int_{\mathfrak{N}^{\circ}(\mathbb{R})} \nu(\mathbb{R}) \mathrm{e}^{-\nu(h)} \mathbf{Q}_{\rm min}(\mathrm{d} \nu),\quad h\in\mathcal{B}_b^+(\mathbb{R}),
			\end{align*}
			where $\varphi(0)=\lim_{n\rightarrow\infty}(1-f_n(0))/m^n$.
			
			\item [{\rm (2)}] For any $\mu \in \mathfrak{N}^{\circ}(\mathbb{R})$, we have  $\lim_{l\rightarrow\infty} \lim_{n\rightarrow\infty} \mathbb{P}_{\mu}(Z_n\in\cdot\mid N_{n+l}>0)= \mathbf{Q}_{\infty}(\cdot)$ by weak convergence, where $\mathbf{Q}_{\infty}(\cdot)$ is defined as in {\rm(1)}.
		\end{itemize}
	\end{corollary}
	
	
	For any probability measure $\mathbf{Q}$ on $\mathfrak{N}^{\circ}(\mathbb{R})$, we write $(\mathbf{Q}\mathbb{P})(\cdot):=\int_{\mathfrak{N}^{\circ}(\mathbb{R})}\mathbb{P}_{\mu}(\cdot)\mathbf{Q}(\mathrm{d}\mu)$. A probability measure $\mathbf{Q}$ on $\mathfrak{N}^{\circ}(\mathbb{R})$ is called the {\it QSD} of $Z$ if
	\begin{align}\label{qsd1}
		(\mathbf{Q}\mathbb{P})(Z_n\in B\mid N_n>0)=\mathbf{Q}(B),\quad n\in\mathbb{N},\ B\in\mathcal{B}(\mathfrak{N}^{\circ}(\mathbb{R})).
	\end{align} 

	The next theorem gives a equation characterization of the QSD in terms of its p.g.fl.
	\begin{proposition}\label{lem1}
		Suppose that $\mathbf{Q}$ is a probability measure on $\mathfrak{N}^{\circ}(\mathbb{R})$. Then $\mathbf{Q}$ is the QSD of $Z$ if and only if its p.g.fl. $G_{\mathbf{Q}}$ satisfies
		\begin{align}\label{qsd3}
			G_{\mathbf{Q}}\big[G[h](\cdot)\big]-G_{\mathbf{Q}}\big[f(0)\big]=m^{\alpha} G_{\mathbf{Q}}\big[h \big],\quad h\in\mathcal{V}_0(\mathbb{R}),
		\end{align}
		for some $\alpha\in(0,\infty)$. In this case, the probability measure $\mathbf{Q}$ is called a QSD of $Z$ with eigenvalue $m^{\alpha}$.
	\end{proposition}
	\begin{theorem}\label{QSD}
		Suppose that the conditions {\rm(C1)}, {\rm(C2)} and {\rm(C4)} hold.
		\begin{itemize}
			\item[{\rm(1)}] When $\alpha>1$, there exist no 
				QSDs of $Z$ with eigenvalue $m^{\alpha}$.
			\item[{\rm(2)}] The only QSD of $Z$ with eigenvalue $m$ is the Yaglom limit $\mathbf{Q}_{\rm min}$.
			\item[{\rm(3)}] For $\alpha\in(0,1)$. A probability measure $\mathbf{Q}_{\alpha}$ on $\mathfrak{N}^{\circ}(\mathbb{R})$ is a QSD of $Z$ with eigenvalue $m^{\alpha}$ if and only if its p.g.fl. $G_{\mathbf{Q}_{\alpha}}[h]$ satisfies
			\begin{align}\label{qsd4}
				G_{\mathbf{Q}_{\alpha}}[h]=\int_{(0,\infty)}\Big(\mathrm{e}^{(G_{\mathbf{Q}_{\rm min}}[h]-1)x}-\mathrm{e}^{-x}\Big)\frac{1}{x^{\alpha}}\Delta(\mathrm{d}x),\quad h\in\mathcal{V}_0(\mathbb{R}),
			\end{align}
			where 
			$\Delta$ is a locally finite measure on $(0,\infty)$ satisfying:
			\begin{itemize}
				\item[{\rm(a)}] $\Delta(B)=\Delta(mB)$ for every $B\in\mathcal{B}((0,\infty))$.
				\item[{\rm(b)}] $\int_{(0,\infty)} (1-\mathrm{e}^{-x}) x^{-\alpha} \Delta (\mathrm{d} x)=1$.
			\end{itemize}
		\end{itemize}
	\end{theorem}
	
	\begin{remark}
		By applying \cite[Lemma 15]{M18} to the function $x\mapsto x^{-\alpha} (1-\mathrm{e}^{-x}),\alpha\in(0,1)$, it is easy to get the existence of the locally finite measure $\Delta$ on $(0,\infty)$ satisfying {\rm (a)} and {\rm (b)}.
	\end{remark}
	
	\begin{remark}
		The measure $\Delta$ is uniquely determined by $\mathbf{Q}_{\alpha}$. In the proof of Theorem \ref{QSD}, the measure $x^{-\alpha}\Delta(\mathrm{d}x)$ is constructed as the vague limit ($n\rightarrow\infty$) of the measure $\kappa_n$ on $(0,\infty)$ defined by $\kappa_n(B) := m^{-\alpha n} \mathbf{Q}_{\alpha}(\{\nu\in\mathfrak{N}^{\circ}(\mathbb{R}): \nu(\mathbb{R})\in q_n^{-1}B\})$, where $B\in\mathcal{B}((0,\infty))$ and $q_n=\mathbb{P}(N_n>0)$.
	\end{remark}
	
	\begin{remark}
		The gamma function $\Gamma(z):=\int_{0}^{\infty} t^{z-1}\mathrm{e}^{-t} \mathrm{d}t, z>0$ can be analytically extended so that it is also defined on $\mathbb{C}\backslash\{0,-1,\cdots\}$, where $\mathbb{C}$ is the field of complex numbers. In particular, for $z\in\mathbb{R}\backslash\{0,-1,\cdots\}$, the gamma function satisfies the recursion 
			$$\Gamma(1+z)=z \Gamma(z) \quad\text{and}\quad \Gamma(1/2)=\sqrt{\pi}.$$ 
			Therefore, $\Gamma(-\alpha)<0$ for $\alpha\in(0,1)$. These and other properties related to the gamma function can be found in \cite{GR14,L72}. It follows from letting $\Delta (\mathrm{d} x)=-\frac{1}{x \Gamma (-\alpha)} \mathrm{d} x$ in Theorem \ref{QSD} that $G_{\mathbf{Q}_{\alpha}}[h] = 1-(1-G_{\mathbf{Q}_{\rm min}}[h])^{\alpha}$, these special QSDs in multitype branching processes with a finite set of types were obtained by Hoppe and Seneta \cite{HS78}.  Later, Asmussen and Hering \cite{AH83} extended these results to multitype branching processes with an infinite set of types (when the type space is $\mathbb{R}$, the multitype branching process is the branching Markov chain).
	\end{remark}
	
	In Section \ref{sec3}, we give some results about conditional reduced branching Markov chains. In Section \ref{sec4.1}, we quickly explain the spinal decomposition theorem and the many-to-few formula, which are fundamental tools in this paper. In Section \ref{sec4.2}, we study the limit of the $r$-th moment of $Z_n(h)$ under the probability $\mathbb{P}_{\delta_x}(\ \cdot\mid N_n>0)$. Based on these preparations, we prove Theorem \ref{Yaglom} using the moment method in Section \ref{sec5.1}, while in Section \ref{sec5.2} we prove Corollaries \ref{cor2.7} and \ref{cor2.8} as byproducts of Theorem \ref{Yaglom}. In Section \ref{sec6.1}, we give the proof of Proposition \ref{lem1}. In Section \ref{sec6.2}, we give the proof of Theorem \ref{QSD} by using Propositions \ref{prop6.3}-\ref{prop6.5}, whose proofs are given in Section \ref{sec6.3}.
	\section{Conditional reduced branching Markov chains} 
	\label{sec3}
	By Lemma \ref{lem2.3}, the way we will approach the proof of  Theorem \ref{Yaglom} is through two fundamental steps. The first step deals with the convergence of the Laplace functional for $Z_n$ under the probability $\mathbb{P}_{\delta_{x}}(\ \cdot\mid N_n>0)$. For the second step, we prove that the above limit functional on $\mathcal{B}_b^+(\mathbb{R})$ is continuous with respect to bound pointwise convergence. 
	
	For any $h\in\mathcal{B}_b^+(\mathbb{R})$, in the first step, it is essential to calculate the moments of $Z_n(h)$ under the probability $\mathbb{P}_{\delta_{x}}(\ \cdot\mid N_n>0)$ and use the result about the moment method, which guarantees that the limit distribution is uniquely determined by its moments. The many-to-few formula plays an important role in calculating these moments. The key ingredient in adapting the many-to-few formula is based an important fact that conditional on $N_n>0$, only the particles in generation $k\ (0\leq k\leq n)$ having descendants in the $n$-th generation make contributions to the integral $Z_n(h)$. The process which is composed by the number of these particles at time $k\ (0\leq k\leq n)$ is called the {\it reduced process} in the literature; see, e.g., \cite{FP74,FKR77,FS78,FV99}. 
	\subsection{Conditional reduced GW-processes}
	\label{sec3.1}
	Recall that $N=(N_n:n\in\mathbb{N})$ is a GW-process with offspring distribution $p$. Let $N_{k,n}$ be the number of particles in the $k$-th generation having non-empty offspring at time $n$, then $(N_{k,n}:0\leq k\leq n)$ is called the {\it reduced GW-process}. Let $\hat{N}_{k,n}$ be a random variable with law
	\begin{align*}
		\mathbb{P}_{\delta_{x}}(\hat{N}_{k,n}\in\cdot\ )=\mathbb{P}_{\delta_{x}}(N_{k,n}\in\cdot\mid N_n>0),
	\end{align*}
	where $\hat{N}_n:=\hat{N}_{n,n}$ by convention. We say the process $(\hat{N}_{k,n}:0\leq k\leq n)$ is the {\it conditional reduced GW-process}. 
	
	Based on the results of reduced GW-processes established in Fleischmann and Siegmund-Schultze \cite{FKR77}, we can get the following properties of the conditional reduced GW-process. 
	\begin{proposition}
		The process $(\hat{N}_{k,n}:0\leq k\leq n)$ is a time-inhomogeneous GW-process. Suppose that $\hat{p}(k,n):=\{\hat{p}_l(k,n)\}_{l\in\mathbb{N}}$ is the offspring distribution of a particle at time $k\ (0\leq k\leq n-1)$, then its generating function is given by
		\begin{align}\label{gs2}
			\hat{f}_{k,n}(s):=\mathbb{E}[s^{\hat{N}_{k+1,n}} \mid \hat{N}_{k,n}=1]
			=\frac{f\big(f_{n-k-1}(0)+s\big(1-f_{n-k-1}(0)\big)\big)-f_{n-k}(0)}{1-f_{n-k}(0)}.
		\end{align}
	\end{proposition}
	
	By differentiating both side of (\ref{gs2}), we get the following proposition.
	\begin{proposition}\label{prop3.2}
		Suppose that the offspring distribution $p$ has finite moments of all orders. For $r\in\mathbb{N}_+$, let $f^{(r)}$ and $\hat{f}_{k,n}^{(r)}$ denote respectively the $r$-th derivative of $f$ and $\hat{f}_{k,n}$. Then we have
		\begin{align*}
			\hat{f}_{k,n}^{(r)}(1)=\frac{\big(1-f_{n-k-1}(0)\big)^{r}}{1-f_{n-k}(0)}f^{(r)}(1).
		\end{align*}
	\end{proposition}
	
	For $j\in\mathbb{N}_+$, we denote by $\hat{m}_j(k,n)$ the $j$-th moment of  $\hat{p}(k,n)$. It is easy to calculate the first and second moments:
	\begin{align}
		&\hat{m}_1(k,n)=m \frac{1-f_{n-k-1}(0)}{1-f_{n-k}(0)},\label{5} \\
		&\hat{m}_2(k,n)=\frac{(1-f_{n-k-1}(0))^2}{1-f_{n-k}(0)} f^{(2)}(1)+m \frac{1-f_{n-k-1}(0)}{1-f_{n-k}(0)}.  \label{6}
	\end{align}
	\subsection{Conditional reduced branching Markov chains}
	\label{sec3.2}
	Recall that under $\mathbb{P}_{\delta_x}$, $\mathbb{T}$ is a GW-tree with offspring distribution $p$, which starts from one particle. By the definition in \cite{FKR77}, we can get a {\it reduced GW-tree} $\mathbb{T}_n$ (rooted at $\varnothing$) by removing all branches of $\mathbb{T}$, which don't extend to generation $n$. Let $\hat{\mathbb{T}}_n$ be a random tree rooted at $\varnothing$ with law
	\begin{align*}
		\mathbb{P}_{\delta_{x}}(\hat{\mathbb{T}}_n\in\cdot\ )=\mathbb{P}_{\delta_{x}}(\mathbb{T}_n\in\cdot\mid N_n>0).
	\end{align*}
	For any $0\leq k\leq n$, let
	$$\hat{\mathbb{G}}_{k,n}:=\{u\in\hat{\mathbb{T}}_n:|u|=k\}$$
	be the set of individuals in generation $k$. The {\it conditional reduced branching Markov chain} $(V(u):u\in\hat{\mathbb{G}}_{k,n},0\leq k\leq n)$ under $\mathbb{P}_{\delta_x}$ is characterized by the following properties:
	\begin{itemize}
		\item[{\rm (1)}] It starts with a single particle $\varnothing$ born at time $0$ and position $x$, which forms the $0$-th generation, i.e., $\hat{\mathbb{G}}_{0,n}=\{\varnothing\}$.
		
		\item[{\rm (2)}] For each $0\leq k\leq n-1$, each particle $u\in\hat{\mathbb{G}}_{k,n}$ (after one unit time) splits independently of the other particles of the $k$-th generation into a random number of children with law $\hat{p}(k,n)$. Conditional on $V(u)=y$, these children are born at positions which are given by independent copies of the transition probability $P(y,\cdot)$.
	\end{itemize}
	
	Let $\hat{Z}_{k,n}(\cdot)$ denote the counting measure of particles of generation $k$ in the above system. Then $(\hat{Z}_{k,n}:0\leq k\leq n)$ is a $\mathfrak{N}(\mathbb{R})$-valued Markov process. 
	Let $\hat{Z}_{k,n}(B)$ denote the number of particles in the generation $k$ located in $B$. In particular, $\hat{N}_{k,n}=\hat{Z}_{k,n}(\mathbb{R})$, which forms a conditional reduced process described in Section \ref{sec3.1}. 
	
	We make the convention that $\hat{Z}_n:=\hat{Z}_{n,n}$. It is easy to see that
	\begin{align*}
		\mathbb{P}_{\delta_{x}}(\hat{Z}_n\in\cdot\ )=\mathbb{P}_{\delta_{x}}(Z_n\in\cdot\mid N_n>0).
	\end{align*}
	Let $\hat{Z_n}(h):=\int_{\mathbb{R}} h(x) \hat{Z}_n(\mathrm{d}x)$, then the p.g.fl. $\hat{G}_{n}$ of the point process $\hat{Z}_{n}$ under $\mathbb{P}_{\delta_{x}}$ is given by
	\begin{align*}
		\hat{G}_{n}[h](x)=\mathbb{E}_{\delta_{x}}\big[\mathrm{e}^{\hat{Z}_{n}(\log h)}\big]=\mathbb{E}_{\delta_{x}}\big.\big[\mathrm{e}^{{Z}_{n}(\log h)}\ \big|\ N_n>0\big],\quad h\in\mathcal{V}_0(\mathbb{R}).
	\end{align*}
	Therefore we only need to consider the moments of $\hat{Z}_n(h)$, and then obtain the convergence of $\hat{G}_{n}$, which completes the first step.
	\section{Moments and Changes of measures}
	\label{sec4}
	The many-to-few formula which is useful to calculate higher moments for branching models has been discussed in many other related contexts, see, for example, \cite{A13,B18,BH15,HS09,ZS15} for the many-to-one formula and \cite{HR17,HL23,HY23} for the many-to-few formula. Here we state the spinal decomposition theorem and the many-to-few formula for conditional reduced branching Markov chains, the above results about time-homogeneous branching models were obtained by Harris and Roberts \cite{HR17}. However the conditional reduced branching Markov chain is time-inhomogeneous, checking out the proof of \cite[Lemma 8]{HR17} and taking the martingale $\zeta\equiv1$, one only need to replace the offspring distribution and the motion distribution of the particle $u$ at
	time $k$ with  $\{\hat{p}_l(k,n)\}_{l\in\mathbb{N}}$ and $P(V(u),\cdot)$. 
	
	\subsection{The spinal decomposition theorem and the many-to-few formula}
	\label{sec4.1} 
	For any $r\in\mathbb{N}_+$, we now attach $r$ additional random distinguished lines of descent $w^1,\cdots,w^r$ to the conditional reduced branching Markov chain, which are called {\it spines}. 
	For any $x\in\mathbb{R}$, let us introduce the following system, which is referred to as the conditional reduced branching Markov chain starting from one particle located at $x$ with $r$ spines $w^1,\cdots,w^r$:
	\begin{itemize}
		\item[(1)] At time $0$, there is one particle at position $x$ carrying $r$ marks $1,\cdots,r$, which forms the $0$-th generation. 
		
		\item[(2)] For each $1\leq i\leq r$, we regard the line of descent carrying mark $i$ as the spine $w^i$. Note that $w^i=(w_1^i,w_2^i,\cdots)$, where $w_k^i$ denotes the particle carrying mark $i$ at time $k$. We write $\mathfrak{X}_k^i$ for its position, i.e., $\mathfrak{X}_k^i=V(w_k^i)$.
		
		\item[(3)] A particle $u$ in generation $k$ carrying $j$ marks (after one unit time) splits independently of all others into a random number of children according to the $j$-th {\it sized-biased} distribution $\big\{\frac{l^j \hat{p}_l(k,n)}{\hat{m}_j(k,n)}\big\}_{l\in\mathbb{N}}$. Conditional on $V(u)=y$, these children are born at positions which are given by independent copies of the law $P(y,\cdot)$. Given that $l$ particles are born at such a branching event, each of $j$ marks chooses a particle to follow independently and uniformly at random among these $l$ variables.
		
		\item[(4)] Particles which carry no marks at time $k$ have children according to the offspring distribution $\{\hat{p}_l(k,n)\}_{l\in\mathbb{N}}$ and the one-step transition probability $P$, just as under $\mathbb{P}_{\delta_x}$.
		
		\item[(5)] The children of particles in generation $k$ and their positions (including the marked particles) form the $(k+1)$-th generation. 
	\end{itemize} 
	
	We denote by $\mathbb{Q}_{\delta_x}^{[r]}$ the law of this new system. In other words, under $\mathbb{Q}_{\delta_x}^{[r]}$, spine particles give birth to children according to the size-biased distribution, but the $n$-step transition probability of the {\it spine Markov process} $(\mathfrak{X}_k^i:k\geq0)$ will be unchanged, which is still $P_n$. The number of children depends on how many marks the spine particle carrying, but the motion does not. 
	
	For $u,v\in\mathcal{U}$, recall that $|v|$ denotes the generation of $v$, $\overleftarrow{v}$ is the parent of $v(\neq\varnothing)$ and $v\prec u$ implies $v$ is an ancestor of $u$. For fixed $u^1,\cdots,u^r\in\mathcal{U}$, let
	\begin{align*}
		&{\rm skel}_{u^1,\cdots,u^r}(k)=\{v\in\mathcal{U}:|v|\leq k,\exists\ 1\leq j\leq r \text{ with } v\prec u^j \},\\
		&D_{u^1,\cdots,u^r}(v)=\#\{1\leq j\leq r:v\prec u^j\},
	\end{align*}
	where $\#B$ denotes the cardinality of the set $B$. Let $(\hat{\mathcal{F}}_{k,n})_{0\leq k\leq n}$ be the natural filtration of the conditional reduced branching Markov chain. We define an $\hat{\mathcal{F}}_{k,n}$-adapted process $(W^{[r]}(k,n):0\leq k\leq n)$ by
	\begin{align}\label{b3}
		W^{[r]}(k,n):=\sum_{u^1,\cdots,u^r\in\hat{\mathbb{G}}_{k,n}}\prod_{v\in{\rm skel}_{u^1,\cdots,u^r}(k)\backslash\{\varnothing\}}\big(\hat{m}_{D_{u^1,\cdots,u^r}(\overleftarrow{v})}(|\overleftarrow{v}|,n)\big)^{-1}.
	\end{align}
	
	The following lemma establishes the link between $\mathbb{P}_{\delta_x}$ and $\mathbb{Q}^{[r]}_{\delta_x}$, and calculates the probability that the $r$ marks are carried by particles $u^1,\cdots,u^{r}$ at time $k$; see, e.g., \cite{HR17}.
	\begin{lemma}\label{lem4.1}{\rm(The spinal decomposition theorem)}
		
		{\rm (1)} The derivative of $\mathbb{Q}^{[r]}_{\delta_x}$ with respect to $\mathbb{P}_{\delta_x}$ on $\hat{\mathcal{F}}_{k,n}$ is $W^{[r]}(k,n)$ defined by {(\ref{b3})}, that is,
		\begin{align*}
			\bigg.\frac{\mathrm{d} \mathbb{Q}_{\delta_x}^{[r]}}{\mathrm{d} \mathbb{P}_{\delta_x}}\bigg|_{\hat{\mathcal{F}}_{k,n}}=W^{[r]}(k,n).
		\end{align*}
		
		{\rm (2)} For any $0\leq k\leq n$ and any $u^1,\cdots,u^{r}\in\hat{\mathbb{G}}_{k,n}$, we have
		\begin{align*}
			\mathbb{Q}_{\delta_x}^{[r]}(w_k^1=u^1,\cdots,w_k^r=u^r \mid \hat{\mathcal{F}}_{k,n}) = \frac{1}{W^{[r]}(k,n)} \prod_{v\in{\rm skel}_{u^1,\cdots,u^r}(k)\backslash\{\varnothing\}}\big(\hat{m}_{D_{u^1,\cdots,u^r}(\overleftarrow{v})}(|\overleftarrow{v}|,n)\big)^{-1}.
		\end{align*} 
	\end{lemma}
	
	Let us introduce the many-to-few formula for the conditional reduced branching Markov chain. Let $\hat{\mathcal{G}}_{k,n}^{[r]}$ be the filtration containing all information about the conditional reduced branching Markov chain (including the $r$ spines) up to time $k$. It follows from \cite{HR17} that if $F$ is measurable with respect to $\hat{\mathcal{G}}_{k,n}^{[r]}$, then it can be expressed as
	$$
	F=\sum_{u^1,\cdots,u^r\in\hat{\mathbb{G}}_{k,n}}F(u^1,\cdots,u^r)\mathbf{1}_{\{w_k^1=u^1,\cdots,w_k^r=u^r\}},
	$$
	where for each $u^1,\cdots,u^r\in\hat{\mathbb{G}}_{k,n}$ the random variable $F(u^1,\cdots,u^r)$ is $\hat{\mathcal{F}}_{k,n}$-measurable. 
	For the sake of simplicity, we still write $\mathbb{Q}_{\delta_x}^{[r]}$ for its associated expectation in this paper.

	\begin{lemma}\label{lem6}{\rm(The many-to-few formula)}
		For $r\in\mathbb{N}_+$ and $0\leq k\leq n$, suppose that $F$ is measurable with respect to $\hat{\mathcal{G}}_{k,n}^{[r]}$. Then we have 
		\begin{align}\label{many1}
			\mathbb{E}_{\delta_x}\bigg[\sum_{u^1,\cdots,u^r\in\hat{\mathbb{G}}_{k,n}}F(u^1,\cdots,u^r)\bigg]=\mathbb{Q}_{\delta_x}^{[r]}\bigg[F(w_k^1,\cdots,w_k^r)\prod_{v\in\text{\rm skel}(k)\backslash\{\varnothing\}}\hat{m}_{D(\overleftarrow{v})}(|\overleftarrow{v}|,n)\bigg],
		\end{align}
		where {\rm skel}$(k)$ denotes the set of particles which carry at least one mark up to time $k$, and  $D(v)$ is the total number of marks carried by the particle $v$.
	\end{lemma}
	\subsection{Moments}
	\label{sec4.2}
	In this section we obtain some results for the moments of $\hat{Z}_n(h)$. Here and in the sequel, we make the convention that $\sum_{j=l}^k=0$ and $\prod_{j=l}^k=1$ for any $k<l$. Let $\hat{\mathbb{G}}_{n}:=\hat{\mathbb{G}}_{n,n}$ for simplicity. 
	It was proved in \cite{J67} that when $0<m<1$, the sequence $(1-f_n(s))/(m^n(1-s))$ is a monotone decreasing in $n$ and
	\begin{align}\label{13}
		\lim_{n\rightarrow\infty}\frac{1-f_n(s)}{m^n(1-s)}=\varphi(s)\in[0,1],\quad \text{ $\varphi(s)$ is a non-decreasing function}.
	\end{align}
	Furthermore, $\varphi(0)>0$ if and only if $\sum_{j}(j\log j)p_j<\infty$; see, e.g., \cite[Theorem 2]{HS67}.
	
	\begin{lemma}\label{lem3.4}{\rm (The first moment)}
		Suppose that the condition {\rm(C1)} holds. For any $h\in\mathcal{B}_b^+(\mathbb{R})$, we have
		\begin{align}\label{14}			\mathbb{E}_{\delta_x}[\hat{Z}_n(h)]=\frac{m^n}{1-f_n(0)}P_nh(x).
		\end{align}
	\end{lemma}
	
	\noindent{\it Proof.}
	Note that $F(w_n^1)=h(V(w_n^1))=h(\mathfrak{X}_n^1)$ implies $F(u^1)=h(V(u^1))$. Then by applying (\ref{many1}) for $r=1,k=n$ and using (\ref{5}) we obtain
	\begin{align*}
		\mathbb{E}_{\delta_x}\big[\hat{Z}_n(h)\big]=\mathbb{E}_{\delta_x}\Big[\sum_{u^1 \in \hat{\mathbb{G}}_{n}} h\big(V(u^1)\big)\Big]=\frac{m^n}{1-f_n(0)}\mathbb{Q}^{[1]}_{\delta_x}\big[h(\mathfrak{X}_n^1)\big],
	\end{align*}
	where $(\mathfrak{X}_n^1:n\geq0)$ is the spine Markov process whose $n$-step transition probability is $P_n$. Hence the  desired result follows from this property.
	$\hfill\square$ 
	
	In view of (\ref{13}), then an application of Lemma \ref{lem3.4} gives the following result.
	\begin{corollary}\label{cor3.5}{\rm (The limit of the first moment)}
		If the conditions {\rm(C1)} and {\rm(C3)} hold, we have
		\begin{align}\label{17}			\lim_{n\rightarrow\infty}\mathbb{E}_{\delta_x}[\hat{Z}_n(h)]=\frac{\pi(h)}{\varphi(0)},\quad h\in\mathcal{B}_b^+(\mathbb{R}).
		\end{align}
	\end{corollary}

	\begin{lemma}\label{lem3.6}{\rm (The second moment)}
		Suppose that the condition {\rm(C1)} holds. For any $h\in\mathcal{B}_b^+(\mathbb{R})$, we have
		\begin{align}\label{15}
			\mathbb{E}_{\delta_x}\big[\big(\hat{Z}_n(h)\big)^2\big]=\frac{m^n}{1-f_n(0)}P_nh^2(x)+\frac{m^nf^{(2)}(1)}{m(1-f_n(0))}\sum_{j=0}^{n-1}m^j \int_{\mathbb{R}}\big(P_{j+1}h(y)\big)^2 P_{n-j-1}(x,\mathrm{d} y).
		\end{align}
	\end{lemma}
	
	\noindent{\it Proof.}
	If  $F(w_n^1,w_n^2)=h(V(w_n^1))h(V(w_n^2))=h(\mathfrak{X}_n^1)h(\mathfrak{X}_n^2)$, then $F(u^1,u^2)=h(V(u^1))h(V(u^2))$. We denote by $T_1^{[2]}$ the first split time at which marks 1 and 2 are carried by different particles. Then applying  (\ref{many1}) for $r=2$ and $k=n$ gets
	\begin{align}\label{2bgs5}
		\mathbb{E}_{\delta_x}\big[\big(\hat{Z}_n(h)\big)^2\big] 
		=&\ \mathbb{E}_{\delta_x}\Big[\sum_{u^1,u^2\in\hat{\mathbb{G}}_{n}} h\big(V(u^1))h(V(u^2)\big) \Big] \nonumber\\
		=&\ \mathbb{Q}_{\delta_x}^{[2]}\bigg[\mathbf{1}_{\{T_1^{[2]}>n\}}\big(h(\mathfrak{X}_{n}^{1})\big)^2\prod_{i=0}^{n-1}\hat{m}_2(i,n)\bigg] \nonumber\\
		&+\sum_{j=1}^{n}\mathbb{Q}_{\delta_x}^{[2]}\bigg[h(\mathfrak{X}_{j-1,n}^{1})h(\mathfrak{X}_{j-1,n}^{2})\mathbf{1}_{\{T_1^{[2]}=j\}}\prod_{i=0}^{j-1}\hat{m}_2(i,n)\prod_{i=j}^{n-1}\big(\hat{m}_1(i,n)\big)^2\bigg],
	\end{align} 
	where $(\mathfrak{X}_{n}^{1}:n\geq0)$ is the spine Markov process carrying two marks, $(\mathfrak{X}_{j-1,n}^{1}:n\geq0)$  and $(\mathfrak{X}_{j-1,n}^{2}:n\geq0)$ are two spine Markov processes having the same trajectory before time $j$ and splitting at time $j$. 
	
	Under $\mathbb{Q}_{\delta_x}^{[2]}$, the initial particle carrying marks $1$ and $2$ splits at time $1$ into $l$ particles with probability $l^2 \hat{p}_l(0,n)/\hat{m}_2(0,n)$. At such a branching event, the two marks follow the same particle with probability $1/l$. Thus
	\begin{align*}
		\mathbb{Q}_{\delta_x}^{[2]}\big(T_1^{[2]}>1\big)=\frac{\hat{m}_1(0,n)}{\hat{m}_2(0,n)}.
	\end{align*}
	Similarly, under $\mathbb{Q}_{\delta_x}^{[2]}$ the particle in generation $i$ carrying two marks splits at time $i+1$ into $l$ particles with law $l^2 \hat{p}_l(i,n)/\hat{m}_2(i,n)$ and then the two marks follow the same particle with probability $1/l$. Hence,
	\begin{align}\label{7}
		\mathbb{Q}_{\delta_x}^{[2]}\big(T_1^{[2]}>j\big)=\prod_{i=0}^{j-1}\frac{\hat{m}_1(i,n)}{\hat{m}_2(i,n)},\quad j=1,2,\cdots.
	\end{align}
	Then, for any integer $j\geq 1$,
	\begin{align}\label{8}
		\mathbb{Q}_{\delta_x}^{[2]}\big(T_1^{[2]}=j\big)=\mathbb{Q}_{\delta_x}^{[2]}\big(T_1^{[2]}>j-1\big)-\mathbb{Q}_{\delta_x}^{[2]}\big(T_1^{[2]}>j\big)=\prod_{i=0}^{j-2}\frac{\hat{m}_1(i,n)}{\hat{m}_2(i,n)}\bigg(1-\frac{\hat{m}_1(j-1,n)}{\hat{m}_2(j-1,n)}\bigg),
	\end{align}
	with $\prod_{i=0}^{-1}=1$ by convention. From (\ref{5}), (\ref{6}), (\ref{7}) and (\ref{8}) it is easy to see that
	\begin{align*}
		\prod_{i=0}^{n-1} \hat{m}_2(i,n) \mathbb{Q}_{\delta_x}^{[2]}\big(T_1^{[2]}>n\big)=\frac{m^n}{1-f_n(0)}
	\end{align*}
	and
	\begin{align*}
		\prod_{i=0}^{j-1}\hat{m}_2(i,n) \prod_{i=j}^{n-1}\big(\hat{m}_1(i,n)\big)^2 \mathbb{Q}_{\delta_x}^{[2]}\big(T_1^{[2]}=j\big)=\frac{m^{2 n-1}}{1-f_n(0)} f^{(2)}(1) m^{-j}.
	\end{align*}
	Since the reproduction and motion of an individual are independent, we obtain
	\begin{align*}
		\begin{small}
			\mathbb{E}_{\delta_x}\big[\big(\hat{Z}_n(h)\big)^2\big]=\frac{m^n}{1-f_n(0)} \mathbb{Q}_{\delta_x}^{[2]}\big[\big(h(\mathfrak{X}_{n}^{1})\big)^2\big]
			+\frac{m^{2n-1}f^{(2)}(1)}{1-f_n(0)}  \sum_{j=1}^n m^{-j}\mathbb{Q}_{\delta_x}^{[2]}\big[h(\mathfrak{X}_{j-1,n}^1)h(\mathfrak{X}_{j-1,n}^2)\big].
		\end{small}
	\end{align*}
	According to the properties of spine Markov processes $\mathfrak{X}_{n}^{1}$, $\mathfrak{X}_{j-1,n}^1$ and $\mathfrak{X}_{j-1,n}^2$,  then
	\begin{align*}
		\mathbb{E}_{\delta_x}\big[\big(\hat{Z}_n(h)\big)^2\big]=&\ \frac{m^n}{1-f_n(0)}P_nh^2(x)+\frac{m^{2 n-1}}{1-f_n(0)} f^{(2)}(1) \sum_{j=1}^n m^{-j} \int_{\mathbb{R}}\big(P_{n-j+1}h(y)\big)^2 P_{j-1}(x, \mathrm{d} y)\\
		=&\ \frac{m^n}{1-f_n(0)}P_nh^2(x)+\frac{m^nf^{(2)}(1)}{m(1-f_n(0))} \sum_{j=0}^{n-1} m^j \int_{\mathbb{R}}\big(P_{j+1}h(y)\big)^2 P_{n-j-1}(x, \mathrm{d} y).
	\end{align*}
	That gives the desired result.
	$\hfill\square$ 
	
	\begin{corollary}\label{cor3.7}{\rm (The limit of the second moment)}
		If the conditions {\rm(C1)} and {\rm(C3)} hold, we have
		\begin{align}\label{18}
			\lim_{n\rightarrow\infty}	\mathbb{E}_{\delta_x}\big[\big(\hat{Z}_n(h)\big)^2\big]=\frac{\pi(h^2)}{\varphi(0)}+\frac{f^{(2)}(1)}{m\varphi(0)} \sum_{j=0}^{\infty}m^j \int_{\mathbb{R}}\big(P_{j+1}h(y)\big)^2 \pi(\mathrm{d} y),\quad h\in\mathcal{B}_b^+(\mathbb{R}).
		\end{align}
	\end{corollary}
	\noindent{\it Proof.}
	By (\ref{13}), it is not hard to see that the first term on the right-hand side of (\ref{15}) converges to $\pi(h^2)/\varphi(0)<\infty$ as $n\rightarrow\infty$ under {\rm (C1)} and {\rm (C3)}.
	Let 
	\begin{align*}
		a_{j,n}:=\mathbf{1}_{\{0\leq j\leq n-1\}}m^j \int_{\mathbb{R}}\big(P_{j+1}h(y)\big)^2 P_{n-j-1}(x, \mathrm{d} y),
	\end{align*}
	then the second term on the right-hand side of (\ref{15}) is equal to 
	\begin{align}\label{19}
		\frac{m^nf^{(2)}(1)}{m(1-f_n(0))} \sum_{j=0}^{\infty}a_{j,n}.
	\end{align}
	Under the conditions {\rm (C1)} and {\rm (C3)}, it is easy to see that $\sum_{j}m^j<\infty$ and
	\begin{align*}
		0\leq a_{j,n}\leq \|h\|^2m^j,\quad \lim_{n\rightarrow\infty}a_{j,n}=m^j \int_{\mathbb{R}}\big(P_{j+1}h(y)\big)^2 \pi(\mathrm{d} y).
	\end{align*}
	Then by dominated convergence we have
	\begin{align*}
		\lim_{n\rightarrow\infty}\sum_{j=0}^{\infty}a_{j,n}=\sum_{j=0}^{\infty}m^j \int_{\mathbb{R}}\big(P_{j+1}h(y)\big)^2 \pi(\mathrm{d} y).
	\end{align*}	
	This combined with (\ref{13}) and (\ref{19}) gives the second term on the right-hand side of (\ref{18}).
	$\hfill\square$ 
	
	Before giving the $r$-th moment of $\hat{Z}_n(h)$, we introduce some notations. For $r_1,l_1,\cdots,r_k,l_k,r\in\mathbb{N}_+$, let
	\begin{align}\label{23}
		C_r(r_1,l_1,\cdots,r_k,l_k):=\frac{r !}{\left(r_{1} !\right)^{l_1} \cdots\left(r_{k} !\right)^{l_k}} \cdot \frac{1}{l_{1} ! \cdots l_{k} !}.
	\end{align} 
	When $l_1 r_1+\cdots+l_k r_k=r$ and $l_1+\cdots+l_k=i$, then $C_r(r_1,l_1,\cdots,r_k,l_k)$ is the number of partitions of dividing $r$ different elements into $i$ indistinguishable blocks, where there are $l_1$ blocks with $r_1$ elements, $\cdots$, $l_k$ blocks with $r_k$ elements. For $i\geq2$, let the sum $\sum_{r_\cdot,l_\cdot,k}$ be over all sequences of positive integers satisfying 
	\begin{equation}\label{16}
		\left\{
		\begin{aligned}
			&l_1 r_1+\cdots+l_k r_k=r,\\
			&l_1+\cdots+l_k=i,\\
			&0< r_1<\cdots<r_k\leq r-i+1,\\
			&1\leq k\leq i.
		\end{aligned}
		\right.
	\end{equation}  
	Let $M_h(0,x,r):=Ph^{r}(x)$ and define $M_h(k_1,x,r)$ for $k_1\in\mathbb{N}$ inductively
	\begin{small}
		\begin{align*}
			M_h(k_1,x,r)\!\!:=\!\!P_{k_1+1}h^{r}(x)\!\!
			+\!\!\sum_{i=2}^{r}\!\!\frac{f^{(i)}(1)}{m}\!\!\sum_{r_\cdot,l_\cdot,k}\!\!\!C_{r}(r_1,l_1,\cdots,r_k,l_k)\!\!\sum_{k_2=0}^{k_1-1}\!\!m^{(i-1)k_2}\!\!\!\int_{\mathbb{R}}\!\bigg[\!\prod_{j=1}^{k}\!\!\big(M_h(k_2,y,r_j)\big)^{l_j}\!\bigg]\!P_{k_1-k_2}(x,\mathrm{d}y),
		\end{align*}
	\end{small}with the convention that $\sum_{i=2}^r=0$ for $r\leq 1$. The following result shows that the $r$-th moment can be represented as a combination of the moments of low order than $r$. 
	
	\begin{lemma}\label{lem4.6}{\rm (The $r$-th moment)}
		Suppose that the condition {\rm(C1)} holds. For any $r\in\mathbb{N}_+$ and any $h\in\mathcal{B}_b^+(\mathbb{R})$, we have 
		\begin{align}\label{12}
			\mathbb{E}_{\delta_x}\big[\big(\hat{Z}_n(h)\big)^r\big]=\Lambda_{r,1}(n,h)+\Lambda_{r,2}(n,h)
		\end{align}
		where $\Lambda_{r,1}(n,h)\!\!=\frac{m^n}{1-f_n(0)}P_nh^r(x) $
		and
		\begin{small}
			\begin{align*}
				\Lambda_{r,2}(n,h)\!\!=\!\!\sum_{i=2}^{r}\!\!\frac{m^nf^{(i)}(1)}{m(1-f_n(0))}\!\!\sum_{r_\cdot,l_\cdot,k}\! C_r(r_1,l_1,\cdots,r_k,l_k)\!\!\sum_{k_1=0}^{n-1}\!\!m^{(i-1)k_1}\!\!\!\!\int_{\mathbb{R}}\!\bigg[\!\prod_{j=1}^{k}\big(M_h(k_1,y,r_j)\big)^{l_j}\!\bigg]\!P_{n-k_1-1}(x,\mathrm{d}y).
			\end{align*}
		\end{small}In particular, we have
		\begin{align}\label{30}
			\int_{\mathbb{R}}\mathbb{E}_{\delta_y}\big[\big(\hat{Z}_n(h)\big)^r\big]P(x,\mathrm{d}y)=\frac{m^n}{1-f_n(0)}M_h(n,x,r).
		\end{align}
	\end{lemma}
	
	\noindent{\it Proof.}
	We shall give the proof by induction in $r\geq1$. For $r=1$ the result follows from (\ref{14}). For $r=i=2$, it is easy to see that $(r_1,l_1,k)=(1,2,1)$ is the unique positive integer solution of (\ref{16}). Note that 
	\begin{align*}
		C_2(1,2)=1\quad\text{and}\quad M_h(k_1,y,1)=P_{k_1+1}h(y).
	\end{align*}
	These combined with (\ref{15}) yield that the result holds for $r=2$.
	Now suppose that (\ref{12}) holds for  $r-1\geq 1$. It is clearly sufficient to consider the $r$-th moment. Let 
	\begin{align*}
		F(w_n^1,\cdots,w_n^r)=h(V(w_n^1))\cdots h(V(w_n^r))=h(V(\mathfrak{X}_n^1))\cdots h(V(\mathfrak{X}_n^r)),
	\end{align*}
	then
	\begin{align*}
		F(u^1,\cdots,u^r)=h(V(u^1))\cdots h(V(u^r)).
	\end{align*}
	Hence by applying (\ref{many1}) for $k=n$ we have
	\begin{small}
		\begin{align}\label{11}
			\mathbb{E}_{\delta_x}\big[\big(\hat{Z}_n(h)\big)^r\big]
			=\mathbb{E}_{\delta_x}\bigg[\sum_{u^1,\cdots,u^r\in\hat{\mathbb{G}}_{n}}\!\!\!\!\!F(u^1,\cdots,u^r)\bigg]
			=\mathbb{Q}_{\delta_x}^{[r]}\bigg[F(w_n^1,\cdots,w_n^r)\!\!\!\!\prod_{v\in\text{\rm{skel}}(n)\backslash\{\varnothing\}}\!\!\!\!\hat{m}_{D(\overleftarrow{v})}(|\overleftarrow{v}|,n)\bigg].
		\end{align}	
	\end{small}We next to express the product in the expectation of (\ref{11}) according to the number of marks carried by particles. Let $T_1^{[r]}$ be the first time at which $r$ marks are carried by different particles. Then the event $\{T_1^{[r]}>n\}$ implies that all $r$ marks are carried by the same particle up to time $n$ and 
	\begin{align*}
		\{T_1^{[r]}\!\leq\! n\}\!=\!\bigcup_{i=2}^{r}\{\text{the $r$ marks are carried by $i$ different particles at the time ($\leq n$) they first split} \}.
	\end{align*}
	More specifically, suppose that $\{k,r_1,l_1,\cdots,r_k,l_k\}$ is a sequence of positive integer satisfying (\ref{16}). Let $T^{[r]}_1(r_1,l_1,\cdots,r_k,l_k)$ be the first time at which $r$ marks are carried by $i$ different particles, of which there are $l_1$ particles carrying $r_1$ marks, $\cdots$, $l_k$ particles carrying $r_k$ marks. Therefore,
	\begin{align*}
		\{T_1^{[r]}\leq n\}=\bigcup_{i=2}^{r} \bigcup_{r_\cdot,l_\cdot,k} \bigcup_{k_1=1}^{n} \{ T_1^{[r]}(r_1,l_1,\cdots,r_k,l_k)=k_1 \}.
	\end{align*}
	Hence, we have
	\begin{small}
		\begin{align}
			\mathbb{E}&_{\delta_x}\big[\big(\hat{Z}_n(h)\big)^r\big]=\mathbb{Q}_{\delta_x}^{[r]}\bigg[\big(h(\mathfrak{X}_{n}^{1})\big)^r\mathbf{1}_{\{T_1^{[r]}>n\}}\prod_{j=0}^{n-1}\hat{m}_r(j,n)\bigg]\nonumber\\
			&+\sum_{i=2}^{r}\sum_{r_\cdot,l_\cdot,k}\sum_{k_1=1}^{n}\mathbb{Q}_{\delta_x}^{[r]}\bigg[h(\mathfrak{X}_{k_1-1,n}^{1})\cdots h(\mathfrak{X}_{k_1-1,n}^{r})\mathbf{1}_{\{T_1^{[r]}(r_1,l_1,\cdots,r_k,l_k)=k_1\}}\!\!\!\!\!\prod_{v\in\text{\rm{skel}}(n)\backslash\{\varnothing\}}\!\!\!\hat{m}_{D(\overleftarrow{v})}(|\overleftarrow{v}|,n)\bigg],\label{21}
		\end{align}
	\end{small}where $(\mathfrak{X}_{n}^{1}:n\geq0)$ is the spine Markov process carrying $r$ marks, and  $(\mathfrak{X}_{k_1-1,n}^{1}:n\geq0),\cdots,(\mathfrak{X}_{k_1-1,n}^{r}:n\geq0)$ are spine Markov processes characterized by the following properties:
	\begin{itemize}
		\item[{\rm (1)}] These $r$ spine Markov processes have the same trajectory up to time $k_1-1$ and split into $i$ independent different trajectories at time $k_1$.
		
		\item[{\rm (2)}] At time $k_1$, these $r$ marks are carried by $i$ different particles, of which there are $l_1$ particles carrying $r_1$ marks, $\cdots$, $l_k$ particles carrying $r_k$ marks.
		
		\item[{\rm (3)}] From time $k_1$ to time $n$, the trajectories of spine Markov processes carrying  $r_j$ marks are equal in law to those in the many-to-$r_j$ formula for the branching process  $(\hat{Z}_{j,n}^{\theta^{k_1}}:0\leq j\leq n-k_1)$, where $\hat{Z}_{j,n}^{\theta^{k_1}}$ is defined as follows. 
	\end{itemize}
	Suppose that $\{\hat{N}^{\theta^{k_1}}(u)\}$ is a family of independent and identically distributed integer-valued random variable with the same law $\hat{p}(k_1+|u|,n)$. Let $\hat{\mathbb{T}}_n^{\theta^{k_1}}$ be the {\it shift tree} of $\hat{\mathbb{T}}_n$ with defining elements $\{\hat{N}^{\theta^{k_1}}(u)\}$: 
	\begin{align*}
		\hat{\mathbb{T}}_n^{\theta^{k_1}}:=\bigcup_{j=0}^{n-k_1}\hat{\mathbb{G}}_{j,n}^{\theta^{k_1}}\quad\text{with}\quad\hat{\mathbb{G}}_{0,n}^{\theta^{k_1}}:=\{\varnothing\}\text{ and\ }\hat{\mathbb{G}}_{j+1,n}^{\theta^{k_1}}:=\{ui:u\in\hat{\mathbb{G}}_{j,n}^{\theta^{k_1}},1\leq i\leq \hat{N}^{\theta^{k_1}}(u)\}.
	\end{align*}
	Then we can get a branching Markov chain $(V(u):u\in\hat{\mathbb{G}}^{\theta^{k_1}}_{j,n},0\leq j\leq n-k_1)$ by endowing each $u$ with the position $V(u)$ in the same way as that in the construction of the branching Markov chain starting from one particle, which is descried in Section \ref{sec2.1}. The corresponding point process is defined by
	\begin{align*}
		\hat{Z}_{j,n}^{\theta^{k_1}}(B):=\sum_{u\in\hat{\mathbb{G}}_{j,n}^{\theta^{k_1}}}\delta_{V(u)}(B),\quad 0\leq j\leq n-k_1,\quad B\in\mathcal{B}(\mathbb{R}).
	\end{align*}
	We write $\hat{Z}_{n-k_1}^{\theta^{k_1}}:=\hat{Z}_{n-k_1,n}^{\theta^{k_1}}$ for simplicity. 
	
	For the first term on the right-hand side of (\ref{21}), by arguments similar to those used for (\ref{7}) we get
	\begin{align}\label{b5}
		\mathbb{Q}_{\delta_x}^{[r]}\big(T_1^{[r]}>n\big)=\prod_{j=0}^{n-1} \frac{\hat{m}_1(j,n)}{\hat{m}_r(j,n)}.
	\end{align}
	It then follows (\ref{5}) that 
	\begin{align*}
		\mathbb{Q}_{\delta_x}\bigg[\big(h(\mathfrak{X}_{n}^{1})\big)^r\mathbf{1}_{\{T_1^{[r]}>n\}}\prod_{j=0}^{n-1}\hat{m}_r(j,n)\bigg]=\frac{m^n}{1-f_n(0)}P_nh^r(x),
	\end{align*}
	which is equal to $\Lambda_{r,1}(n,h)$. For the second term on the right-hand side of (\ref{21}), note that
	\begin{align*}
		\{T_1^{[r]}&(r_1,l_1,\cdots,r_k,l_k)=k_1\}
		=\{\text{$r$ marks follow the same particle until time $k_1-1$; at time $k_1$} \\
		&\text{they are carried by $i$ different particles, of which there are $l_1$ particles carrying $r_1$ marks,}  \\
		&\text{$\cdots$, $l_k$ particles carrying $r_k$ marks} \}.
	\end{align*}
	Under $\mathbb{Q}_{\delta_x}^{[r]}$ the particle in generation $k_1$ carrying $r$ marks splits at time $k_1+1$ into $l$ particles with probability $l^r\hat{p}_l(k_1,n)/\hat{m}_r(k_1,n)$. At such a branching event, the $r$ marks are carried by $i$ different particles, of which there are $l_1$ particles carrying $r_1$ marks, $\cdots$, $l_k$ particles carrying $r_k$ marks with probability 
	\begin{align*}
		C_r(r_1,l_1,\cdots,r_k,l_k)\cdot\frac{l(l-1)\cdots(l-i+1)}{l^r},
	\end{align*}
	where $C_r(r_1,l_1,\cdots,r_k,l_k)$ is defined by (\ref{23}). By the above comments and (\ref{b5}) we have
	\begin{align*}
		&\ \mathbb{Q}_{\delta_x}^{[r]}\big(T_1^{[r]}(r_1,l_1,\cdots,r_k,l_k)=k_1\big)\nonumber\\
		=&\ C_r(r_1,l_1,\cdots,r_k,l_k)\Big[\sum_{l\geq1}\frac{l^r\hat{p}_l(k_1-1,n)}{\hat{m}_{r}(k_1-1,n)}\times\frac{l(l-1)\cdots(l-i+1)}{l^r}\Big]\times\prod_{j=0}^{k_1-2}\frac{\hat{m}_{1}(j,n)}{\hat{m}_{r}(j,n)}\nonumber\\
		=&\ C_r(r_1,l_1,\cdots,r_k,l_k)\frac{\hat{f}_{k_1-1,n}^{(i)}(1)}{\hat{m}_{r}(k_1-1,n)}\times\prod_{j=0}^{k_1-2}\frac{\hat{m}_{1}(j,n)}{\hat{m}_{r}(j,n)}\nonumber\\
		=&\ C_r(r_1,l_1,\cdots,r_k,l_k)\frac{f^{(i)}(1)}{\hat{m}_{r}(k_1-1,n)}\frac{\big(1-f_{n-k_1}(0)\big)^i}{1-f_{n-k_1+1}(0)}\times\prod_{j=0}^{k_1-2}\frac{\hat{m}_{1}(j,n)}{\hat{m}_{r}(j,n)},
	\end{align*}
	where the last equality comes from Proposition \ref{prop3.2}. Recall that the reproduction and motion of a particle are independent. Then by the property (4) of spine Markov processes, it is easy to infer that the second term on the right-hand side of (\ref{21}) is equal to
	\begin{small}
		\begin{align*}
			\sum_{i=2}^{r}\sum_{r_\cdot,l_\cdot,k}\sum_{k_1=1}^{n}\!\!\mathbb{Q}_{\delta_x}^{[r]}\!\big(T_1^{[r]}(r_1,l_1,\cdots,r_k,l_k)\!=\!k_1\!\big)\!\bigg[\prod_{j=0}^{k_1-1}\!\!\!\hat{m}_r(j,n)\bigg]\!\!\int_{\mathbb{R}}\Big[\prod_{j=1}^{k}\!\!\Big(\!\!\int_{\mathbb{R}}\!\!\!\mathbb{E}_{\delta_y}\!\big[\big(\hat{Z}_{n-k_1}^{\theta^{ k_1}}\!(h)\big)^{r_{j}}\big]P(\widetilde{y},\mathrm{d}y)\Big)^{l_{j}}\Big]\!P_{k_1-1}(x,\mathrm{d}\widetilde{y}).
		\end{align*}
	\end{small}By the definition of $\hat{Z}_{n-k_1}^{\theta^{k_1}}$ and $M_h(k_1,x,r)$, then we use the induction hypothesis to get
	\begin{align*}
		\int_{\mathbb{R}}\mathbb{E}_{\delta_y}\big[\big(\hat{Z}_{n-k_1}^{\theta^{k_1}}(h)\big)^{r_{j}}\big]P(\widetilde{y},\mathrm{d}y) =\frac{m^{n-k_1}}{1-f_{n-k_1}(0)}M_h(n-k_1,\widetilde{y},r_j),\quad 1\leq k_1\leq n.
	\end{align*}
	In view of (\ref{5}), we obtain
	\begin{align*}
		\frac{\big(1-f_{n-k_1}(0)\big)^i}{1-f_{n-k_1+1}(0)}\times\bigg[\prod_{j=0}^{k_1-2}\hat{m}_1(j,n)\bigg]\times\bigg[\prod_{j=1}^{k}\Big(\frac{m^{n-k_1}}{1-f_{n-k_1}(0)}\Big)^{l_j}\bigg]=\frac{m^n}{m(1-f_n(0))} m^{(i-1)(n-k_1)}.
	\end{align*}
	Combining the four preceding equalities gives that the second term on the right-hand side of (\ref{21}) is $\Lambda_{r,2}(n,h)$. Then we complete the proof.
	$\hfill\square$ 
	
	We now give some characterizations of the Yaglom limit for GW-processes, which play a very important role in studying the limit of the $r$-th moment. Recall that $\nu_{\text{min}}$ is the Yaglom limit defined in Section \ref{sec1}. Let $c_1^{-1}=\varlimsup_{n\rightarrow\infty}(p_n)^{\frac{1}{n}}$, then the condition {\rm(C2)} implies $c_1>1$ and 
	\begin{align*}
		\sum_{n\in\mathbb{N}}p_n|s|^n<\infty\quad\text{for all}\quad|s|<c_1.
	\end{align*}
	Hence, the generating function $f(s)=\sum_{n\in\mathbb{N}} p_n s^n$ is well-defined for $|s|<c_1$.
	\begin{proposition}\label{prop4.7}
		Suppose that the conditions {\rm(C1)} and {\rm(C2)} hold. 
		\begin{itemize}
			\item[{\rm(1)}] There exits a constant $c_2\in(1,c_1)$ such that $f_n(s)$ is monotone decreasing in $n$ and $f_n(s)\downarrow 1$ as $n \rightarrow \infty$ for all $s\in(1,c_2)$.
			\item[{\rm(2)}] $\lim_{n\rightarrow\infty}\mathbb{E}[\mathrm{e}^{c_3\hat{N}_n}]<\infty$ for some constant $c_3>0$.
			\item[{\rm(3)}] {\rm (Cramér's condition)} Let $\psi(t)=\sum_{n\in\mathbb{N}_+}\nu_{\text{\rm min}}(n)\mathrm{e}^{tn}$, then $\psi(t)<\infty$ for all $t\in(-c_3,c_3)$ where $c_3>0$ is a constant.
		\end{itemize} 
	\end{proposition}
	
	\noindent{\it Proof.}  {\rm (1)} Since $m>0$, the generating function $f(s)$ is strictly increasing on $s\in[0,c_1)$. Take $c_2\in(1,c_1)$ such that $f(c_2)<c_2$, then we have $1<f(s)<s$ for $s\in(1,c_2)$. Iterating this inequality one can see 
	\begin{align*}
		1<f_n(s)<f_{n-1}(s)<\cdots<f(s)<s<c_2,\quad\text{for all}\quad n\geq 1.
	\end{align*}
	Thus $f_n(s) \downarrow c_{\infty}\in[1,c_2)$. Since $f$ is a continuous function, we may take limits through the equation $f_{n+1}(s)=f(f_n(s))$ and conclude that $c_{\infty}=f(c_{\infty})$. Note that $c_{\infty}<c_2$ implies that there are no roots of $s=f(s)$ in $(1,c_2)$, and hence $c_{\infty}=1$.
	
	{\rm (2)} Note that
	\begin{align*}
		\mathbb{E}\big[\mathrm{e}^{c_3\hat{N}_n}\big]=\mathbb{E}\big.\big[\mathrm{e}^{c_3 N_n}\ \big|\ N_n>0\big]=1+\frac{f_n(\mathrm{e}^{c_3})-1}{1-f_n(0)},
	\end{align*}
	then it suffices to show there exists $\widetilde{c}_3>1$ so that the sequence $(f_n(\widetilde{c}_3)-1)/(1-f_n(0))$ converges. 
	For $s\in(1,c_2)$, by the Taylor's formula of $f$ at $s=1$ we have 
	\begin{align}\label{4.1}
		\frac{f(s)-1}{s-1}=m+\epsilon(s),
	\end{align}
	where
	\begin{align}\label{4.4}
		\epsilon(s):=\frac{f^{(2)}(1+\beta(s-1))}{2}(s-1),\quad \beta\in(0,1).
	\end{align}
	One can easily check that $\epsilon(s)$ is monotone increasing for $s\in(1,c_2)$ and $\lim_{s\downarrow1}\epsilon(s)=0$. Replacing $s$ by $f_{k-1}(s)$ in (\ref{4.1}) yields 
	\begin{align*}
		\frac{f_k(s)-1}{f_{k-1}(s)-1}=m\Big[1+\frac{\epsilon(f_{k-1}(s))}{m}\Big],\quad s\in(1,c_2).
	\end{align*}
	Taking the product on the both side of the above equality for $k=1,\cdots,n$ we obtain
	\begin{align}\label{4.2}
		\frac{f_n(s)-1}{s-1}=m^n \times \prod_{k=0}^{n-1}\Big[1+\frac{\epsilon\big(f_k(s)\big)}{m}\Big],\quad s\in(1,c_2).
	\end{align}
	It is clear that the sequence $(f_n(s)-1)/(m^n(s-1))$ is monotone increasing in $n$ and thus converges to a non-decreasing function $\widetilde{\varphi}(s)$. By (\ref{4.4}), we see that $\widetilde{\varphi}(s)<\infty$ if and only if
	\begin{align*}
		\sum_{k=0}^{\infty} f^{(2)}\big(1+\beta(f_k(s)-1)\big)\big(f_k(s)-1\big)<\infty.
	\end{align*}
	Given a constant $0<\varepsilon<c_2-1$ satisfying $f^{(1)}(1+\varepsilon)<1$, then by the assertion {\rm(1)}, there exists $k_0\geq 1$ such that $1<f_k(s)<1+\varepsilon$ for $k\geq k_0$. It then follows that 
	\begin{align*}
		f^{(2)}\big(1+\beta(f_{k+k_0}(s)-1)\big) \leq f^{(2)}(1+\varepsilon),\quad s\in(1,c_2).
	\end{align*}
	For $s\in(1,c_2)$, by the mean-value theorem, there exists a constant $\beta_1\in(0,1)$ such that
	\begin{align*}
		f_{k+k_0}(s)-1=&\ f^{(1)}\big(1+\beta_1(f_{k+k_0-1}(s)-1)\big) \times\big(f_{k+k_0-1}(s)-1\big)\\
		\leq&\ f^{(1)}(1+\varepsilon) \times\big(f_{k+k_0-1}(s)-1\big)\\
		\leq&\ \varepsilon\Big[f^{(1)}(1+\varepsilon)\Big]^k.
	\end{align*}
	Combining the three preceding inequalities gives that $\widetilde{\varphi}(s)<\infty$ for $s\in(1,c_2)$. Thus we get the desired result.
	
	{\rm (3)} It is easy to check that $\psi(t)<\infty$ for $t\leq 0$. If $t>0$, let $c_3>0$ be chosen as in the assertion {\rm(2)}. By Fatou's Lemma and Skorokhod's theorem we get
	\begin{align*}
		\psi(c_3)\leq\lim_{n\rightarrow\infty}\mathbb{E}\big[\mathrm{e}^{c_3\hat{N}_n}\big].
	\end{align*}
	Then the assertion {\rm(3)} follows by {\rm (2)}.
	$\hfill\square$
	
	We continue to study the asymptotic behavior of the $r$-th moment.
	\begin{corollary}\label{cor4.8}{\rm (The limit of the $r$-th moment)}
		Suppose that the conditions {\rm(C1)}-{\rm(C3)} hold. For any $r\in\mathbb{N}_+$ and any $h\in\mathcal{B}_b^+(\mathbb{R})$, then 
		\begin{align}\label{20}
			\lim_{n\rightarrow\infty}\mathbb{E}_{\delta_x}\big[\big(\hat{Z}_n(h)\big)^r\big]=\Lambda_{r,1}(h)+\Lambda_{r,2}(h),
		\end{align}
		where $\Lambda_{r,1}(h)=\frac{\pi(h^r)}{\varphi(0)}$ and 
			\begin{align*}
				\Lambda_{r,2}(h)=\sum_{i=2}^r \frac{f^{(i)}(1)}{m\varphi(0)}\sum_{r_\cdot,l_\cdot,k}\! C_r(r_1,l_1,\cdots,r_k,l_k)\sum_{k_1=0}^{\infty}\!m^{(i-1)k_1}\!\!\int_{\mathbb{R}}\!\bigg[\prod_{j=1}^{k}\big(M_h(k_1,y,r_j)\big)^{l_j}\bigg]\pi(\mathrm{d}y).
			\end{align*}
	Moreover, let $\Lambda_r(h):=\Lambda_{r,1}(h)+\Lambda_{r,2}(h)$, then we have $\Lambda_r(h)<\infty$.
	\end{corollary}
	
	\noindent{\it Proof.} For $r\in\mathbb{N}_+$ and $h\in\mathcal{B}_b^+(\mathbb{R})$, suppose that $\Lambda_{r,1}(n,h)$ and $\Lambda_{r,2}(n,h)$ are defined as in Lemma \ref{lem4.6}. From (\ref{13}) it is easy to get that $\lim_{n\rightarrow\infty}\Lambda_{r,1}(n,h)=\Lambda_{r,1}(h)<\infty$ under {\rm(C3)}.  We next consider the limit of $\Lambda_{r,2}(n,h)$. Let
	\begin{align*}
		a_{i,k_1,n}(h):=\mathbf{1}_{\{0\leq k_1\leq n-1\}}m^{(i-1)k_1} \int_{\mathbb{R}}\bigg[\prod_{j=1}^{k}\big(M_h(k_1,y,r_j)\big)^{l_j}\bigg] P_{n-k_1-1}(x, \mathrm{d} y).
	\end{align*}
	Then we have
	\begin{align}\label{4.5}
		\lim_{n\rightarrow\infty}\Lambda_{r,2}(n,h)=&\ \frac{1}{\varphi(0)}\times\lim_{n\rightarrow\infty}\sum_{i=2}^{r}\frac{f^{(i)}(1)}{m}\sum_{r_\cdot,l_\cdot,k} C_r(r_1,l_1,\cdots,r_k,l_k)\sum_{k_1=0}^{\infty}a_{i,k_1,n}(h)\nonumber\\
		=&\ \sum_{i=2}^{r}\frac{f^{(i)}(1)}{m\varphi(0)}\sum_{r_\cdot,l_\cdot,k} C_r(r_1,l_1,\cdots,r_k,l_k)\times\lim_{n\rightarrow\infty}\sum_{k_1=0}^{\infty}a_{i,k_1,n}(h).
	\end{align}
	From (\ref{30}) we see 
	\begin{align*}
		0\leq M_h(k_1,y,r_j)\leq\|h\|^{r_j} \mathbb{E}\big[\big(\hat{N}_{k_1}\big)^{r_j}\big]\frac{1-f_{k_1}(0)}{m^{k_1}}.
	\end{align*}
	It then follows that
	\begin{align*}
		0\leq a_{i,k_1,n}(h)\leq\|h\|^{r} m^{(i-1)k_1}\bigg(\frac{1-f_{k_1}(0)}{m^{k_1}}\bigg)^{i}\times\bigg[\prod_{j=1}^{k}\big(\mathbb{E}\big[\big(\hat{N}_{k_1}\big)^{r_j}\big]\big)^{l_j}\bigg].
	\end{align*}
	Under the conditions {\rm(C1)} and {\rm(C2)}, from the assertion {\rm(2)} of Proposition \ref{prop4.7} we have $\lim_{n\rightarrow\infty}\mathbb{E}[(\hat{N}_n)^r]<\infty$ for $r\geq1$, and then $\sup_{n}\mathbb{E}[(\hat{N}_n)^r]:=d_r<\infty$. In view of (\ref{13}), we also have that $(1-f_n(0))/m^n\leq 1$. Thus
	\begin{align}\label{4.3}
		0\leq a_{i,k_1,n}(h)\leq\|h\|^{r} m^{(i-1)k_1}  \times \bigg[\prod_{j=1}^{k}d_{r_j}^{l_j}\bigg]
	\end{align}
	Note that $\sum_{k_1}m^{(i-1)k_1}<\infty$ for $2\leq i\leq r$ and
	\begin{align*}
		\lim_{n\rightarrow\infty}a_{i,k_1,n}(h)=m^{(i-1)k_1} \int_{\mathbb{R}}\bigg[\prod_{j=1}^k \big(M_h(k_1,y,r_j)\big)^{l_j}\bigg] \pi(\mathrm{d} y):=a_{i,k_1,\infty}(h),\quad\text{(under {\rm (C3)})}.
	\end{align*}
	Then by dominated convergence we get $\sum_{k_1}a_{i,k_1,n}(h)$ converges to $\sum_{k_1}a_{i,k_1,\infty}(h)$ as $n\rightarrow\infty$. This combined with (\ref{4.5}) implies $\lim_{n\rightarrow\infty}\Lambda_{r,2}(n,h)=\Lambda_{r,2}(h)$. Hence, (\ref{20}) holds. In addition, by (\ref{4.5}) and (\ref{4.3}) we see
	\begin{align*}
		\Lambda_{r,2}(h)\leq\frac{\|h\|^{r}}{m\varphi(0)}\sum_{i=2}^{r}\frac{f^{(i)}(1)}{1-m^{i-1}}\sum_{r_\cdot,l_\cdot,k} C_r(r_1,l_1,\cdots,r_k,l_k) \times \bigg[\prod_{j=1}^{k}d_{r_j}^{l_j}\bigg]<\infty.
	\end{align*}
	Therefore $\Lambda_r(h)=\Lambda_{r,1}(h)+\Lambda_{r,2}(h)<\infty$ for any $h\in\mathcal{B}_b^+(\mathbb{R})$.
	$\hfill\square$ 
	\section{Proofs of Theorem \ref{Yaglom}, Corollaries \ref{cor2.7} and \ref{cor2.8}}
	\label{sec5}
	\subsection{Proof of Theorem \ref{Yaglom}}
	\label{sec5.1}
	Suppose that the conditions {\rm (C1)}-{\rm (C3)} hold, then by Corollary \ref{cor4.8}  we see that
	\begin{align*} \Lambda_r(h)=\lim_{n\rightarrow\infty}\mathbb{E}_{\delta_x}\big[\big(Z_n(h)\big)^r\ \big|\ N_n>0\big]<\infty,\quad h\in\mathcal{B}_b^+(\mathbb{R}). 
	\end{align*} 
	The following result is an application of the moment method, which describes the weak convergence of the conditional distribution of $Z_n(h)$ given $N_n>0$.
	\begin{lemma}\label{lem5.1}
		Suppose that the conditions {\rm(C1)}-{\rm(C3)} hold. For any $x\in\mathbb{R}$ and any $h\in\mathcal{B}_b^+(\mathbb{R})$, then the distribution of $Z_n(h)$ under $\mathbb{P}_{\delta_x}(\ \cdot\mid N_n>0)$ converges weakly as $n\rightarrow\infty$ to that of the random variable $\eta_h$, where the distribution of $\eta_h$ is uniquely determined by its moments $\Lambda_r(h)$.
	\end{lemma}
	\noindent{\it Proof.} Since $\lim_{n\rightarrow\infty}\mathbb{E}_{\delta_x}[(Z_n(h))^r\mid N_n>0]=\Lambda_r(h)<\infty$ under {\rm(C1)}-{\rm(C3)}, by \cite[Theorem 3.3.12]{D10} it suffices to show that
	\begin{align}\label{5.8}
		\varlimsup_{r\rightarrow\infty}\frac{(\Lambda_{2r}(h))^{1/2r}}{2r}<\infty.
	\end{align}
	Clearly, $\mathbb{E}_{\delta_x}[(\hat{Z}_n(h))^r]\leq \|h\|^r\mathbb{E}[(\hat{N}_n)^r]$ implies that $\Lambda_r(h)\leq\|h\|^r\Lambda_r(1)$. Thus we only need to prove
	\begin{align}\label{5.1}
		\varlimsup_{r\rightarrow\infty}\frac{(\Lambda_{2r}(1))^{1/2r}}{2r}<\infty.
	\end{align}
	By the assertion {\rm(2)} of Proposition \ref{prop4.7}, we can easily get that $\sup_n\mathbb{E}[(\hat{N}_n)^r]<\infty$ for all $r\geq1$ under {\rm(C1)} and {\rm(C2)}. Recall that $\hat{N}_n$ converges in distribution to $\eta_1$ as $n\rightarrow\infty$. Then 
	\begin{align*}
		\mathbb{E}[(\eta_1)^r]=\lim_{n \rightarrow \infty}\mathbb{E}\big[\big(\hat{N}_n\big)^r\big]=\Lambda_r(1).
	\end{align*}
	Hence, by using the Taylor's expansion of $\mathrm{e}^x$ and monotone convergence, we see
		\begin{align*}
			\mathbb{E}[\mathrm{e}^{t\,\eta_1}] = \mathbb{E}\Big[ \sum_{r=0}^{\infty} \frac{(t \eta_1)^r}{r!} \Big] =  \sum_{r=0}^{\infty} \frac{t^r}{r!} \mathbb{E}\big[ (\eta_1)^r\big] = \sum_{r=0}^{\infty} \frac{t^r}{r!} \Lambda_r(1).
		\end{align*}
		From the assertion {\rm(3)} of Proposition \ref{prop4.7}, the conditions {\rm(C1)} and {\rm(C2)} imply that $\mathbb{E}[\mathrm{e}^{t\eta_1}]<\infty$ for all $t\in(-c_3,c_3)$, where $c_3>0$. This means the radius of convergence of $\sum_{r=0}^{\infty} (t^r \Lambda_r(1))/r!$ is greater that $0$, i.e.,
		\begin{align*}
			\varlimsup_{r\rightarrow\infty} \Big(\frac{\Lambda_r(1)}{r!}\Big)^{1/r} <\infty.
		\end{align*}
		Stirling's formula (see \cite{F68}, p.52) tells us that as $r\rightarrow\infty$, 
		\begin{align*}
			r! \sim r^r \mathrm{e}^{-r} \sqrt{2 \pi r}.  
		\end{align*}
		Combining the two preceding results gives (\ref{5.1}). Thus (\ref{5.8}) holds.
	$\hfill\square$ 
	
	\begin{remark}\label{zhu1.1}
		Taking $h\equiv1$ in Lemma \ref{lem5.1} we can see that the distribution of $\eta_1$ is the Yaglom limit defined in Section \ref{sec1}. Furthermore, Lemma \ref{lem5.1} shows that the Yaglom limit has finite moments $\Lambda_r(1)$ of all orders under the conditions {\rm(C1)} and {\rm(C2)}. 
	\end{remark}
	
	For any $h\in\mathcal{B}_b^+(\mathbb{R})$, suppose that the conditions {\rm(C1)}-{\rm(C3)} hold and $F_h$ denotes the distribution of the random variable $\eta_h$ given by Lemma \ref{lem5.1}. Since $Z_n(h)$ is a $[0,\infty)$-valued random variable, then by Lemma \ref{lem5.1} and the continuity theorem for Laplace transforms we have
	\begin{align}\label{5.3}
		\mathbb{E}_{\delta_x}\big.\big[\mathrm{e}^{-\lambda Z_n(h)}\ \big|\ N_n>0\big]\xrightarrow[n\rightarrow\infty]{}\int_{[0,\infty)}\mathrm{e}^{-\lambda y}\mathrm{d}F_h(y),\quad x\in\mathbb{R},\quad \lambda\geq0,
	\end{align}
	where the distribution $F_h$ is uniquely determined by its moments $\Lambda_r(h)$. In particular, 
	\begin{align}\label{5.4}
		\mathbb{E}_{\delta_x}\big.\big[\mathrm{e}^{-Z_n(h)}\ \big|\ N_n>0\big]\xrightarrow[n\rightarrow\infty]{}\int_{[0,\infty)}\mathrm{e}^{-y}\mathrm{d}F_h(y),\quad x\in\mathbb{R}.
	\end{align}
	
	\begin{lemma}\label{lem5.2}
		Suppose that the conditions {\rm(C1)}-{\rm(C3)} hold. For any $x\in\mathbb{R}$, there exists a unique probability measure $\mathbf{Q}_{\rm min}$ on $\mathfrak{N}^{\circ}(\mathbb{R})$ such that 
		\begin{align}\label{5.6}
			\mathbb{P}_{\delta_x}(Z_n\in\cdot\mid N_n>0)\xrightarrow[n \rightarrow \infty]{\text{\rm w}} \mathbf{Q}_{\rm min}(\cdot)
		\end{align}
		and its Laplace functional is given by
		\begin{align}\label{b1}
			L_{\mathbf{Q}_{\rm min}} [h] = \int_{[0,\infty)}\mathrm{e}^{-x}\mathrm{d}F_h(x),\quad h\in\mathcal{B}_b^+(\mathbb{R}),
		\end{align}
		where the distribution $F_h$ is uniquely determined by its moments $\Lambda_r(h)$.
	\end{lemma}
	\noindent{\it Proof.} Assume that the conditions {\rm(C1)}-{\rm(C3)} hold. For any $x\in\mathbb{R}$, we shall apply Lemma \ref{lem2.3} for $\mathbf{Q}_n(\cdot)=\mathbb{P}_{\delta_x}(Z_n\in\cdot\mid N_n>0)$ to prove the existence of $\mathbf{Q}_{\rm min}$. By (\ref{5.4}) we obtain that the assertion {\rm(1)} of Lemma \ref{lem2.3} holds and
	\begin{align*}
		L[h] = \int_{[0,\infty)}\mathrm{e}^{-x}\mathrm{d}F_h(x),\quad h\in\mathcal{V}_0(\mathbb{R}),
	\end{align*}
	where $F_h$ is uniquely decided by its moments $\Lambda_r(h)$. We next prove the assertion (2) of Lemma \ref{lem2.3}. Suppose that the sequence  $\{h,h_1,h_2,\cdots\}\subset\mathcal{B}_b^+(\mathbb{R})$ satisfies $h_n\rightarrow h$ and $\sup_n \|h_n\|<\infty$, then it suffices to show $L(h_n)\rightarrow L(h)$. 
	Let $c_4:=\sup_n\|h_n\|<\infty$. For any $r\in\mathbb{N}_+$, by Corollary \ref{cor4.8} we have
	\begin{align*}
		\Lambda_r(h_n)=\Lambda_{r,1}(h_n)+\Lambda_{r,2}(h_n),\quad \Lambda_r( h)=\Lambda_{r,1}(h)+\Lambda_{r,2}(h).
	\end{align*}
	By dominated convergence it is easy to show that $\Lambda_{r,1}(h_n)\rightarrow\Lambda_{r,1}(h)$ as $n\rightarrow\infty$. On the other hand, by similar calculations to (\ref{4.3}) we obtain
	\begin{align*}
		m^{(i-1)k_1} \bigg[\prod_{j=1}^k \big(M_{h_n}(k_1,y,r_j)\big)^{l_j}\bigg]\leq c_4^r m^{(i-1)k_1} \times \bigg[\prod_{j=1}^{k}d_{r_j}^{l_j}\bigg].
	\end{align*}
	We conclude by dominated convergence that
		\begin{align*}
			\lim_{n\rightarrow\infty}\Lambda_{r,2}(h_n)\!=\!\sum_{i=2}^r\!\frac{f^{(i)}(1)}{m\varphi(0)}\!\!\sum_{r_\cdot,l_\cdot,k}\! C_r(r_1,l_1,\cdots,r_k,l_k)\!\!\sum_{k_1=0}^{\infty}m^{(i-1)k_1}\!\!\int_{\mathbb{R}}\lim_{n\rightarrow\infty}\bigg[\prod_{j=1}^k \!\big(M_{h_n}\!(k_1,y,r_j)\big)^{l_j}\!\bigg]\pi(\mathrm{d}y).
		\end{align*} 
	Since $\sup_{n}\|h_n\|=c_4<\infty$, $\mathbb{E}[(\hat{N}_{k_1})^{r_j}]<\infty$ and $h_n \rightarrow h$, by dominated convergence we have
		\begin{align*}
			\lim_{n\rightarrow\infty} \int_{\mathbb{R}} \mathbb{E}_{\delta_{\widetilde{y}}} \big[ \big(\hat{Z}_{k_1} (h_n)\big)^{r_j} \big] P(y, \mathrm{d}\widetilde{y}  ) = \int_{\mathbb{R}} \mathbb{E}_{\delta_{\widetilde{y}}} \big[\big( \hat{Z}_{k_1} (h)\big)^{r_j} \big] P(y, \mathrm{d}\widetilde{y}  ).
		\end{align*}
		This combined with (\ref{30}) implies 
	\begin{align*}
		\lim_{n\rightarrow\infty}M_{h_n}(k_1,y,r_j)= M_{h}(k_1,y,r_j).
	\end{align*}
	Then we get $\lim_{n\rightarrow\infty}\Lambda_{r,2}(h_n)=\Lambda_{r,2}(h)$ from the definition of $\Lambda_{r,2}(h)$. In conclusion, 
	\begin{align}\label{5.7}
		\Lambda_r(h_n)\xrightarrow [n\rightarrow\infty]{}\Lambda_r(h),\quad r\in\mathbb{N}.
	\end{align} 
	Recall that $\Lambda_r(h_n)$ and $\Lambda_r(h)$ are the $r$-th moment of $F_{h_n}$ and $F_h$, respectively. Hence, in view of (\ref{5.8}) and (\ref{5.7}), we obtain that $F_{h_n}$ converges weakly to $F_h$ by using \cite[Theorem 3.3.12]{D10}. Then $L(h_n)\rightarrow L(h)$ follows from the fact that the function $x\mapsto\mathrm{e}^{-x}$ on $x\in[0,\infty)$ is bound continuous. Summing up, by Lemma \ref{lem2.3}, there exists a unique probability measure $\mathbf{Q}_{\text{min}}$ on $\mathfrak{N}(\mathbb{R})$ such that (\ref{5.6}) and (\ref{b1}) hold.
	
	We now prove that $\mathbf{Q}_{\text{min}}$ is the probability measure on $\mathfrak{N}^{\circ}(\mathbb{R})$. Note that
	\begin{align*}
		\mathbf{Q}_{\text{min}}\big(\{\mu \in \mathfrak{N}(\mathbb{R}): \mu=0\}\big) 
		=&\ \lim_{\lambda\rightarrow\infty}\int_{\mathfrak{N}(\mathbb{R})}\mathrm{e}^{-\lambda\mu(\mathbb{R})}\mathbf{Q}_{\text{min}}(\mathrm{d}\mu)\\
		=&\ \lim_{\lambda\rightarrow\infty} \lim_{n\rightarrow\infty} \mathbb{E}_{\delta_x}\big.\big[\mathrm{e}^{-\lambda N_n}\ \big|\  N_n>0\big]\quad\text{(by setting $h\equiv1$ in (\ref{5.3}))} \\
		=&\ \lim_{\lambda\rightarrow\infty} \sum_{j\in\mathbb{N}_+}\nu_{\text{min}}(j)\mathrm{e}^{-\lambda j} \quad\text{(by the definition of the Yaglom limit)} \\
		=&\ 0.
	\end{align*}
	Then we obtain the desired result.
	$\hfill\square$

	We shall give the proof of Theorem \ref{Yaglom}, which is an application of Lemma \ref{lem5.2}.
	
	\noindent\textbf{Proof of Theorem \ref{Yaglom}}. Without loss of generality, we can assume that $\mu=\sum_{i=1}^{l}\delta_{x_i}\in\mathfrak{N}^{\circ}(\mathbb{R})$ with $l\geq1, x_i\in\mathbb{R}$. By the construction of the branching Markov chain starting from $\mu$, then $Z_n$ is equal in law (under $\mathbb{P}_{\mu}$) to $Z_{n}^{(1)}+\cdots+Z_{n}^{(l)}$, where $Z_{n}^{(1)},\cdots,Z_{n}^{(l)}$ are mutually independent random variables satisfying 
	\begin{align*}
		\mathbb{P}_{\mu}(Z_n^{(i)}\in\cdot\ )=\mathbb{P}_{\delta_{x_i}}(Z_n\in\cdot\ ),\quad\text{for}\ i\in\{1,\cdots,l\}.
	\end{align*}
	Let $N_n^{(i)}:=Z_n^{(i)}(\mathbb{R})$, then $N_n$ is equal in law (under $\mathbb{P}_{\mu}$) to $N_{n}^{(1)}+\cdots+N_{n}^{(l)}$ and $N_{n}^{(1)},\cdots,N_{n}^{(l)}$ are independent copies of $N_n$. 
	Hence, for any $h\in\mathcal{B}_b^+(\mathbb{R})$, 
	\begin{align}
		&\ \mathbb{E}_\mu\big[\big.\mathrm{e}^{-Z_n(h)}\ \big|\ N_n > 0\big]\nonumber\\
		=&\ \frac{\mathbb{E}_\mu\big[\mathrm{e}^{- Z_n(h)}\big]-\mathbb{P}_\mu(N_n=0)}{\mathbb{P}_\mu(N_n> 0)}\nonumber\\
		=&\ \frac{\mathbb{E}_\mu\big[\prod_{i=1}^{l} \mathrm{e}^{-Z_n^{(i)}(h)}\big]-\mathbb{P}_\mu\big(\sum_{i=1}^{l} N_n^{(i)}=0\big)}{1-\mathbb{P}_\mu\big(\sum_{i=1}^{l} N_n^{(i)}=0\big)}\nonumber\\
		=&\ \frac{\prod_{i=1}^{l}\mathbb{E}_{\delta_{x_i}}\big[ \mathrm{e}^{-Z_n(h)}\big]-\big(\mathbb{P}(N_n=0)\big)^l}{1-\big(\mathbb{P}(N_n=0)\big)^l}\nonumber\\
		=&\ \frac{\prod_{i=1}^{l}\mathbb{E}_{\delta_{x_i}}\big[ \mathrm{e}^{-Z_n(h)}\big]-\big(\mathbb{E}\big[\mathrm{e}^{-Z_n(h)}\big]\big)^l}{1-\big(\mathbb{P}(N_n=0)\big)^l}+\frac{\big(\mathbb{E}\big[\mathrm{e}^{-Z_n(h)}\big]\big)^l-\big(\mathbb{P}(N_n=0)\big)^l}{1-\big(\mathbb{P}(N_n=0)\big)^l}.\label{5.12}
	\end{align}
	We shall prove the first term on the right-hand side of (\ref{5.12}) tends to $0$ as $n\rightarrow\infty$. Using \cite[Lemma 3.4.3]{D10} and the equality
	\begin{align}\label{b6}
		a^l-b^l=(a-b)(a^{l-1}+a^{l-2} b+\cdots+b^{l-1}),
	\end{align}
	we obtain
	\begin{align*}
		&\ \frac{\prod_{i=1}^{l}\mathbb{E}_{\delta_{x_i}}\big[ \mathrm{e}^{-Z_n(h)}\big]-\big(\mathbb{E}\big[\mathrm{e}^{-Z_n(h)}\big]\big)^l}{1-\big(\mathbb{P}(N_n=0)\big)^l}\\
		\leq&\ \frac{\sum_{i=1}^{l} \Big|\mathbb{E}_{\delta_{x_i}}\big[\mathrm{e}^{-Z_n(h)}\big]-\mathbb{E}\big[\mathrm{e}^{-Z_n(h)}\big]\Big|}{\big(1-\mathbb{P}(N_n=0)\big)\times\big(1+\mathbb{P}(N_n=0)+\cdots+\big(\mathbb{P}(N_n=0)\big)^{l-1}\big)}\\
		\leq&\ \Big[\sum_{i=0}^{l-1}\big(\mathbb{P}(N_n=0)\big)^{i}\Big]^{-1}\times\sum_{i=1}^{l}\bigg|\mathbb{E}_{\delta_{x_i}}\big.\big[\mathrm{e}^{-Z_n(h)}\ \big|\ N_n>0\big]-\mathbb{E}\big.\big[\mathrm{e}^{-Z_n(h)}\ \big|\ N_n>0\big]\bigg|.
	\end{align*}
	From Lemmas \ref{lem2.5} and \ref{lem5.2}, under {\rm(C1)}-{\rm(C3)}, there exists a unique probability measure $\mathbf{Q}_{\text{min}}$ on $\mathfrak{N}^{\circ}(\mathbb{R})$ such that
	\begin{align}\label{5.13}
		\lim_{n\rightarrow\infty}\mathbb{E}_{\delta_ x}\big.\big[\mathrm{e}^{-Z_n(h)}\ \big|\ N_n>0\big]=L_{\mathbf{Q}_{\text{min}}}[h],\quad\text{for all }x\in\mathbb{R}.
	\end{align}
	Then the first term of (\ref{5.12}) tends to $0$ by the fact that $\lim_{n\rightarrow\infty}\mathbb{P}(N_n=0)=1$. 
	
	For the second term on the right-hand side of (\ref{5.12}), by (\ref{b6}) we have
		\begin{align*}
			&\ \frac{\big(\mathbb{E}\big[\mathrm{e}^{-Z_n(h)}\big]\big)^l-\big(\mathbb{P}(N_n=0)\big)^l}{1-\big(\mathbb{P}(N_n=0)\big)^l}\\
			=&\  \frac{\mathbb{E}\big[\mathrm{e}^{-Z_{n}(h)}\big]- \mathbb{P}(N_n=0)}{ 1- \mathbb{P}(N_n=0) } \times \frac{\sum_{i=0}^{l-1}\big(\mathbb{E}\big[\mathrm{e}^{-Z_{n}(h)}\big]\big)^i\big(\mathbb{P}(N_n=0)\big)^{l-i-1}}{\sum_{i=0}^{l-1}\big(\mathbb{P}(N_n=0)\big)^i} \\
			=&\ \mathbb{E}\big.\big[\mathrm{e}^{ -Z_{n} (h)}\ \big|\ N_n>0\big] \times \frac{\sum_{i=0}^{l-1}\big(\mathbb{E}\big[\mathrm{e}^{-Z_{n}(h)}\big]\big)^i\big(\mathbb{P}(N_n=0)\big)^{l-i-1}}{\sum_{i=0}^{l-1}\big(\mathbb{P}(N_n=0)\big)^i}.
		\end{align*}
		Note that both $\mathbb{P}(N_n=0)$ and $\mathbb{E}[\mathrm{e}^{-Z_n(h)}]$ converge to $1$ as $n\rightarrow\infty$. Therefore by (\ref{5.13}), the right-hand side of the above equality tends to $L_{\mathbf{Q}_{\text{min}}}[h]$.
	$\hfill\square$ 
	\subsection{Proofs of Corollaries \ref{cor2.7} and \ref{cor2.8}}
	\label{sec5.2}
	
	\noindent\textbf{Proof of Corollary \ref{cor2.7}.} For fixed $l\in\mathbb{N}$, $\mu\in\mathfrak{N}^{\circ}(\mathbb{R})$ and  $h\in\mathcal{B}_b^+(\mathbb{R})$,  we have
	\begin{align}\label{5.15}
		&\ \mathbb{E}_\mu \big.\big[\mathrm{e}^{-Z_n(h)}\ \big|\  N_{n+l}>0\big] \nonumber\\
		=&\ \frac{\mathbb{P}_\mu(N_n>0)}{\mathbb{P}_\mu(N_{n+l}>0)} \times \frac{\mathbb{E}_\mu \big[\mathrm{e}^{- Z_n(h)} ; N_{n+l}>0\ \big]}{\mathbb{P}_\mu(N_n>0)} \nonumber\\
		=&\ \frac{\mathbb{P}_\mu(N_n>0)}{\mathbb{P}_\mu(N_{n+l}>0)} \times \mathbb{E}_\mu \big.\big[\mathrm{e}^{- Z_n(h)} (1-f_l(0)^{N_n})\ \big|\ N_n>0\big]\quad \text{(by the branching property)} \nonumber\\
		=&\ \frac{\mathbb{P}_\mu(N_n>0)}{\mathbb{P}_\mu(N_{n+l}>0)} \times \bigg( \mathbb{E}_\mu \big.\big[\mathrm{e}^{- Z_n(h)}\ \big|\ N_n>0\big] - \mathbb{E}_\mu \big.\big[\mathrm{e}^{- Z_n(h)} f_l(0)^{N_n}\ \big|\ N_n>0\big]  \bigg).
	\end{align}
	From Lemma \ref{lem2.5} and Theorem \ref{Yaglom} we see that under {\rm(C1)}-{\rm(C3)},
	\begin{align}\label{5.16}
		\lim_{n\rightarrow\infty} \mathbb{E}_{\mu} \big.\big[ \mathrm{e}^{-Z_n(h)}\ \big|\ N_n>0 \big] = L_{\mathbf{Q}_{\text{min}}}[h],\quad h\in\mathcal{B}_b^+(\mathbb{R}),
	\end{align}
	where $L_{\mathbf{Q}_{\text{min}}}$ is the Laplace functional of $\mathbf{Q}_{\text{min}}$. Since $f_l(0)\geq p_0>0$, by replacing $h$ with $h-\log f_l(0)$ in the above equality we obtain
	\begin{align*}
		\lim_{n\rightarrow\infty} \mathbb{E}_{\mu} \big.\big[ \mathrm{e}^{-Z_n(h)} f_l(0)^{N_n} \ \big|\ N_n>0 \big] = L_{\mathbf{Q}_{\text{min}}}\big[h-\log f_l(0)\big].
	\end{align*}
	Note that $f_n(0)\rightarrow1$ as $n\rightarrow\infty$ implies 
	$\frac{ 1-f_n(0)^{\mu(\mathbb{R})} }{ 1-f_n(0) }\rightarrow\mu(\mathbb{R})$, then
	\begin{align*}
		\frac{\mathbb{P}_\mu(N_n>0)}{\mathbb{P}_\mu(N_{n+l}>0)} = \frac{1-f_n(0)^{\mu(\mathbb{R})}}{1-f_n(0)} \times \frac{1-f_n(0)}{1-f_l\big(f_n(0)\big)} \times \frac{1-f_{n+l}(0)}{1-f_{n+l}(0)^{\mu(\mathbb{R})}} \xrightarrow[n\rightarrow\infty]{} m^{-l}.
	\end{align*}
	Combining the four preceding equalities gives
	\begin{align}\label{5.17}
		\lim_{n\rightarrow\infty} \mathbb{E}_\mu \big.\big[\mathrm{e}^{-Z_n(h)}\ \big|\ N_{n+l}>0\big] = m^{-l} \Big(  L_{\mathbf{Q}_{\text{min}}}[h]- L_{\mathbf{Q}_{\text{min}}}\big[h-\log f_l(0)\big] \Big). 
	\end{align} 
	Let $h_n,h\in\mathcal{B}_b^+(\mathbb{R})$ satisfy $h_n\rightarrow h$ and $\sup_n\|h_n\|<\infty$. By dominated convergence it is easy to show that
	\begin{align}\label{5.18}
		\lim_{n\rightarrow\infty} L_{\mathbf{Q}_{\text{min}}}[h_n] = L_{\mathbf{Q}_{\text{min}}}[h],\quad  \lim_{n\rightarrow\infty} L_{\mathbf{Q}_{\text{min}}}\big[h_n-\log f_l(0)\big] =  L_{\mathbf{Q}_{\text{min}}}\big[h-\log f_l(0)\big].
	\end{align}
	In view of (\ref{5.17}) and (\ref{5.18}), by Lemma \ref{lem2.3}, there exists a unique probability measure $\mathbf{Q}_{(l)}$ on $\mathfrak{N}(\mathbb{R})$ such that
	\begin{align*}
		\mathbb{P}_{\mu} (Z_n\in\cdot \mid N_{n+l}>0) \xrightarrow[n\rightarrow\infty]{\text{\rm w}} \mathbf{Q}_{(l)}(\cdot)
	\end{align*}
	and its Laplace functional satisfies $L_{\mathbf{Q}_{(l)}} [h] = m^{-l} ( L_{\mathbf{Q}_{\text{min}}}[h]- L_{\mathbf{Q}_{\text{min}}}[h-\log f_l(0)])$. Note that
	\begin{align}\label{b7}
		\mathbf{Q}_{(l)}(\{ \mu\in\mathfrak{N}(\mathbb{R}):\mu=0 \}) 
		&= \lim_{\lambda\rightarrow\infty} \int_{\mathfrak{N}(\mathbb{R})} \mathrm{e}^{-\lambda \mu(\mathbb{R})} \mathbf{Q}_{(l)}(\mathrm{d}\mu) 
		= \lim_{\lambda\rightarrow\infty} L_{\mathbf{Q}_{(l)}} [\lambda] \nonumber\\
		&= m^{-l}
		\times \lim_{\lambda\rightarrow\infty} \big( L_{\mathbf{Q}_{\text{min}}}[\lambda]- L_{\mathbf{Q}_{\text{min}}}\big[\lambda-\log f_l(0)\big] \big).
	\end{align}
	Since $\mathbf{Q}_{\text{min}}$ is a probability measure on $\mathfrak{N}^{\circ}(\mathbb{R})$, we have
	\begin{align*}
		&\lim_{\lambda\rightarrow\infty} L_{\mathbf{Q}_{\text{min}}}[\lambda]=\mathbf{Q}_{\text{min}}(\{ \mu\in\mathfrak{N}(\mathbb{R}):\mu=0 \})=0
	\end{align*}
	and $\lim_{\lambda\rightarrow\infty} L_{\mathbf{Q}_{\text{min}}}[\lambda-\log f_l(0)]=0$ similarly. These combined with (\ref{b7}) show that $\mathbf{Q}_{(l)}$ is a probability measure on $\mathfrak{N}^{\circ}(\mathbb{R})$.
	$\hfill\square$ 
	
	\
	
	\noindent\textbf{Proof of Corollary \ref{cor2.8}.} 
	{\rm (1)} For any $\mu\in\mathfrak{N}^{\circ}(\mathbb{R})$ and any $h\in\mathcal{B}_b^+(\mathbb{R})$, in view of (\ref{5.15}) we get
	\begin{align}\label{5.19}
		&\ \mathbb{E}_\mu \big.\big[\mathrm{e}^{-Z_n(h)}\ \big|\  N_{n+l}>0\big] \nonumber\\
		=&\ \frac{\mathbb{P}_\mu(N_n>0)}{\mathbb{P}_\mu(N_{n+l}>0)} \times \mathbb{E}_\mu \big.\big[\mathrm{e}^{- Z_n(h)} (1-f_l(0)^{N_n}) \ \big|\  N_n>0\big]  \nonumber\\
		=&\ \frac{1-f_l(0)}{1-f_{n+l}(0)^{\mu(\mathbb{R})}} \times \mathbb{E}_{\mu} \Big[ \frac{\mathrm{e}^{- Z_n(h)} (1-f_l(0)^{N_n})}{1-f_l(0)}\mathbf{1}_{\{ N_{n}>0 \}} \Big].
	\end{align}
	It is easy to check that
	\begin{align}\label{5.20}
		\frac{1-f_l(0)}{1-f_{n+l}(0)^{\mu(\mathbb{R})}} = \frac{1-f_l(0)}{1-f_n\big( f_l(0) \big)} \times \frac{1-f_{n+l}(0)}{1-f_{n+l}(0)^{\mu(\mathbb{R})}} \xrightarrow[l\rightarrow\infty]{} \frac{1}{m^n \mu(\mathbb{R})}.
	\end{align}
	We first consider the expectation on the right-hand side of (\ref{5.19}). By (\ref{b6}) we have
	\begin{align*}
		0\leq \frac{\mathrm{e}^{- Z_n(h)} (1-f_l(0)^{N_n})}{1-f_l(0)} \mathbf{1}_{\{ N_{n}>0 \}} =  \mathrm{e}^{- Z_n(h)} \big( f_l(0)^{N_n-1}+\cdots+1 \big) \leq N_n.
	\end{align*}
	Since $(1-f_l(0)^{N_n})/(1-f_l(0))\rightarrow N_n$ as $l\rightarrow\infty$ and 
	\begin{align}\label{5.21}
		\mathbb{E}_{\mu} [N_n] = \mu(\mathbb{R})m^n<\infty,
	\end{align}
	by dominated convergence we obtain
	\begin{align*}
		\lim_{l\rightarrow\infty} \mathbb{E}_{\mu} \Big[ \frac{\mathrm{e}^{- Z_n(h)} (1-f_l(0)^{N_n})}{1-f_l(0)} \mathbf{1}_{\{ N_{n}>0 \}}  \Big] = \mathbb{E}_{\mu} \big. \big[ N_n \mathrm{e}^{- Z_n(h)}  \mathbf{1}_{\{ N_{n}>0 \}}  \big] = \mathbb{E}_{\mu} \big. \big[ N_n \mathrm{e}^{- Z_n(h)}  \big].
	\end{align*}
	This together with (\ref{5.19}) and (\ref{5.20}) yields
	\begin{align}\label{5.22}
		\lim_{l\rightarrow\infty} \mathbb{E}_\mu \big.\big[\mathrm{e}^{-Z_n(h)} \ \big|\  N_{n+l}>0\big] = \frac{\mathbb{E}_{\mu}  \big[ N_n \mathrm{e}^{- Z_n(h)} \big]}{m^n \mu(\mathbb{R})} .
	\end{align}
	Let $h_r,h\in\mathcal{B}_b^+(\mathbb{R})$ satisfy $\sup_r\|h_r\|<\infty$ and $\lim_{r\rightarrow\infty}h_r= h$. Note that $Z_n$ is a finite point process, then by dominated convergence it is easy to get that $Z_n(h_r)\rightarrow Z_n(h)$ as $r\rightarrow\infty$. Hence $\lim_{r\rightarrow\infty}N_n \mathrm{e}^{-Z_n(h_r)}= N_n \mathrm{e}^{-Z_n(h)}$. We also note that $0\leq N_n \mathrm{e}^{-Z_n(h_r)} \leq N_n$. Therefore, by (\ref{5.21}) and dominated convergence we get
	\begin{align*}
		\lim_{r\rightarrow\infty} \mathbb{E}_{\mu} \big[ N_n \mathrm{e}^{- Z_n(h_r)}  \big] = \mathbb{E}_{\mu}  \big[ N_n \mathrm{e}^{- Z_n(h)} \big].
	\end{align*}
	According to Lemma \ref{lem2.3}, the above equality combined with (\ref{5.22}) implies that there exits a unique probability measure $\mathbf{Q}_{\{n\}}$ on $\mathfrak{N}(\mathbb{R})$ such that 
	\begin{align}\label{5.23}
		\mathbb{P}_{\mu} (Z_n\in\cdot \mid N_{n+l}>0) \xrightarrow[l\rightarrow\infty]{\text{\rm w}} \mathbf{Q}_{\{n\}}(\cdot)
	\end{align}
	and its Laplace functional is given by
	\begin{align}\label{5.24}
		L_{\mathbf{Q}_{\{n\}}} [h] = \big(m^n \mu(\mathbb{R})\big)^{-1} \mathbb{E}_{\mu}  \big[ N_n \mathrm{e}^{- Z_n(h)} \big],\quad h\in\mathcal{B}_b^+(\mathbb{R}).
	\end{align}
	By dominated convergence it is simple to get 
	\begin{align*}
		\mathbf{Q}_{\{n\}}(\{ \mu\in\mathfrak{N}(\mathbb{R}):\mu=0 \}) 
		= \lim_{\lambda\rightarrow\infty} L_{\mathbf{Q}_{\{n\}}} [\lambda] 
		= \lim_{\lambda\rightarrow\infty}  \big(m^n \mu(\mathbb{R})\big)^{-1} \mathbb{E}_{\mu}  \big[ N_n \mathrm{e}^{- \lambda N_n}  \big]=0.
	\end{align*}
	Summing up, there is a unique probability measure $\mathbf{Q}_{\{n\}}$ on $\mathfrak{N}^{\circ}(\mathbb{R})$ satisfying (\ref{5.23}) and (\ref{5.24}). 
	
	We next consider the asymptotic behavior of $\mathbf{Q}_{\{n\}}$ as $n\rightarrow\infty$. When $h$ is replaced by $\lambda_1+\lambda_2 h$ in (\ref{5.16}), we obtain that under {\rm(C1)}-{\rm(C3)},
	\begin{align*}
		\lim_{n \rightarrow \infty} \mathbb{E}_\mu \big.\big[\mathrm{e}^{-\lambda_1 N_n-\lambda_2 Z_n(h)}\ \big|\ N_n>0\big] = L_{\mathbf{Q}_{\text{min}}} [\lambda_1+\lambda_2 h]= \int_{\mathfrak{N}^{\circ}(\mathbb{R})} \mathrm{e}^{-\nu(\lambda_1+\lambda_2 h)} \mathbf{Q}_{\text{min}}(\mathrm{d} \nu),\ \lambda_1,\lambda_2\geq0.
	\end{align*}
	Therefore, it is easy to check that
	\begin{align*}
		\lim_{n \rightarrow \infty} \mathbb{E}_{\mu} \big.\big[ N_n \mathrm{e}^{- Z_n(h)} \ \big|\ N_n>0 \big] = \int_{\mathfrak{N}^{\circ}(\mathbb{R})} \nu(\mathbb{R}) \mathrm{e}^{-\nu(h)} \mathbf{Q}_{\text{min}}(\mathrm{d} \nu).
	\end{align*} 
	Under {\rm(C1)}, (\ref{13}) implies $m^{-n}\mathbb{P}_{\mu}(N_n>0)\rightarrow \mu(\mathbb{R})\varphi(0)$ as $n\rightarrow\infty$. Consequently, letting $n\rightarrow\infty$ in (\ref{5.24}) we see
	\begin{align}\label{5.25}
		\lim_{n\rightarrow\infty} L_{\mathbf{Q}_{\{n\}}} [h] 
		&= \lim_{n\rightarrow\infty} \frac{\mathbb{P}_{\mu}(N_n>0)}{m^n \mu(\mathbb{R})} \mathbb{E}_{\mu} \big.\big[ N_n \mathrm{e}^{- Z_n(h)} \ \big|\ N_n>0 \big] \nonumber \\
		&=\varphi(0)\int_{\mathfrak{N}^{\circ}(\mathbb{R})} \nu(\mathbb{R}) \mathrm{e}^{-\nu(h)} \mathbf{Q}_{\text{min}}(\mathrm{d} \nu).
	\end{align}
	Suppose that the sequence $\{h,h_1,h_2,\cdots\}\subset\mathcal{B}_b^+(\mathbb{R})$ satisfies $h_n\rightarrow h$ and $\sup_n\|h_n\|<\infty$. Then for all $\nu\in\mathfrak{N}(\mathbb{R})$ we have $\lim_{n\rightarrow\infty}\nu(h_n)=\nu(h)$. On the other hand, by Remark \ref{zhu1.1} we see that under {\rm(C1)} and {\rm(C2)}, $\mathbb{E}[\eta_1]={1}/{\varphi(0)}$. From this and Theorem \ref{Yaglom} it follows that
	\begin{align}\label{kz1}
		\int_{\mathfrak{N}^{\circ}(\mathbb{R})} \nu(\mathbb{R}) \mathbf{Q}_{\text{min}}(\mathrm{d} \nu)=\frac{1}{\varphi(0)}<\infty.
	\end{align}
	Therefore, by dominated convergence we obtain
	\begin{align}\label{5.26}
		\lim_{n\rightarrow\infty} \int_{\mathfrak{N}^{\circ}(\mathbb{R})} \nu(\mathbb{R}) \mathrm{e}^{-\nu(h_n)} \mathbf{Q}_{\text{min}}(\mathrm{d} \nu) = \int_{\mathfrak{N}^{\circ}(\mathbb{R})} \nu(\mathbb{R}) \mathrm{e}^{-\nu(h)} \mathbf{Q}_{\text{min}}(\mathrm{d} \nu).
	\end{align}
	Using (\ref{5.25}), (\ref{5.26}) and Lemma \ref{lem2.3} it is simple to deduce that there exits a unique probability measure $\mathbf{Q}_{\infty}$ on $\mathfrak{N}^{\circ}(\mathbb{R})$ such that $\lim_{n\rightarrow\infty}\mathbf{Q}_{\{n\}}=\mathbf{Q}_{\infty}$ by weak convergence and 
	\begin{align*}
		L_{\mathbf{Q}_{\infty}} [h] = \varphi(0)\int_{\mathfrak{N}^{\circ}(\mathbb{R})} \nu(\mathbb{R}) \mathrm{e}^{-\nu(h)} \mathbf{Q}_{\text{min}}(\mathrm{d} \nu),\quad h\in\mathcal{B}_b^+(\mathbb{R}).
	\end{align*}
	Then the assertion {\rm(1)} follows from (\ref{5.23}).
	
	{\rm (2)} We assume that conditions {\rm(C1)}-{\rm(C3)} hold. Suppose that $\mathbf{Q}_{(l)}$ is given in Corollary \ref{cor2.7} and its Laplace functional $L_{\mathbf{Q}_{(l)}}$ is defined by (\ref{2.7}). Then
	\begin{align*}
		\lim_{l\rightarrow\infty} L_{\mathbf{Q}_{(l)}} [h] 
		= \lim_{l\rightarrow\infty} \frac{ L_{\mathbf{Q}_{\text{min}}}[h]- L_{\mathbf{Q}_{\text{min}}}\big[h-\log f_l(0)\big] }{m^{l} } 
		= \lim_{l\rightarrow\infty} \int_{\mathfrak{N}^{\circ}(\mathbb{R})} \mathrm{e}^{-\nu(h)} \frac{1-f_l(0)^{\nu(\mathbb{R})}}{m^l} \mathbf{Q}_{\text{min}}(\mathrm{d} \nu).
	\end{align*}
	Note that 
	\begin{align*}
		& \frac{1-f_l(0)^{\nu(\mathbb{R})}}{m^l} = \frac{1-f_l(0)^{\nu(\mathbb{R})}}{1-f_l(0)} \times  \frac{1-f_l(0)}{m^l} \xrightarrow[l\rightarrow\infty]{} \nu(\mathbb{R}) \varphi(0), \\
		& \frac{1-f_l(0)^{\nu(\mathbb{R})}}{m^l} = \frac{1-f_l(0)^{\nu(\mathbb{R})}}{1-f_l(0)} \times  \frac{1-f_l(0)}{m^l} \leq \nu(\mathbb{R}).
	\end{align*}
	In view of (\ref{kz1}) and the above results, by dominated convergence we get
	\begin{align*}
		\lim_{l\rightarrow\infty} L_{\mathbf{Q}_{(l)}} [h] = \varphi(0)\int_{\mathfrak{N}^{\circ}(\mathbb{R})} \nu(\mathbb{R}) \mathrm{e}^{-\nu(h)} \mathbf{Q}_{\text{min}}(\mathrm{d} \nu) = L_{\mathbf{Q}_{\infty}} [h],\quad h\in\mathcal{B}_b^+(\mathbb{R}),
	\end{align*}
	which implies the assertion {\rm(2)}.
	$\hfill\square$ 
	\section{Proofs of Proposition \ref{lem1} and Theorem \ref{QSD}}
	\label{sec6}
	For a probability measure $\mathbf{Q}$ on $\mathfrak{N}^{\circ}(\mathbb{R})$, recall that \begin{align*}
		(\mathbf{Q}\mathbb{P})(\cdot):=\int_{\mathfrak{N}^{\circ}(\mathbb{R})}\mathbb{P}_{\mu}(\cdot)\mathbf{Q}(\mathrm{d}\mu),\quad (\mathbf{Q}\mathbb{E})[\cdot]:=\int_{\mathfrak{N}^{\circ}(\mathbb{R})}\mathbb{E}_{\mu}[\cdot]\mathbf{Q}(\mathrm{d}\mu).
	\end{align*} 
	Let $G_{\mathbf{Q}}[h]:=\int_{\mathfrak{N}(\mathbb{R})} \mathrm{e}^{\mu(\log h)} \mathbf{Q}(\mathrm{d}\mu)$ be the p.g.fl. of $\mathbf{Q}$. Suppose that $G_n[h](x):=\mathbb{E}_{\delta_x}[\mathrm{e}^{Z_n(\log h)}]$ is the p.g.fl. of $Z_n$ under $\mathbb{P}_{\delta_x}$ and $G:=G_1$ by convention. For any $n\in\mathbb{N}$ and any $h\in\mathcal{V}_0(\mathbb{R})$, it follows from the independence that
	\begin{align}\label{6.1}
		(\mathbf{Q} \mathbb{E})\big[\mathrm{e}^{ Z_n(\log h)}\big]=&\int_{\mathfrak{N}^{\circ}(\mathbb{R})}\mathbb{E}_{\mu}\big[\mathrm{e}^{ Z_n(\log h)}\big] \mathbf{Q}(\mathrm{d}\mu)\nonumber\\
		=&\int_{\mathfrak{N}^{\circ}(\mathbb{R})}\exp\bigg\{\int_{\mathbb{R}} \log G_n[h](x) \mu(\mathrm{d} x)\bigg\} \mathbf{Q}(\mathrm{d}\mu)\nonumber\\
		=&G_{\mathbf{Q}}\big[G_n[h](\cdot)\big].
	\end{align}
	Taking $h\equiv\mathrm{e}^{-\lambda}$ in (\ref{6.1}) yields $(\mathbf{Q} \mathbb{E})[\mathrm{e}^{ -\lambda N_n}] = G_{\mathbf{Q}}[ f_n(\mathrm{e}^{ -\lambda}) ]$. Then by dominated convergence it is easy to check that
	\begin{align}\label{6.2}
		(\mathbf{Q} \mathbb{P}) (N_n=0)=\lim_{\lambda\rightarrow\infty} (\mathbf{Q} \mathbb{E}) \big[\mathrm{e}^{-\lambda N_n}\big] = \lim_{\lambda\rightarrow\infty} G_{\mathbf{Q}}[ f_n(\mathrm{e}^{ -\lambda}) ] = G_{\mathbf{Q}}[ f_n(0) ].
	\end{align}
	Suppose that $\mathbf{Q}_{\text{min}}$ is the Yaglom limit of $Z$. Then by a combination of Lemmas \ref{lem2.4} and  \ref{lem5.2}, we get the following result, which gives a characterization of Yaglom's theorem in terms of p.g.fl.s: under the conditions {\rm(C1)}-{\rm(C3)},
	\begin{align}\label{bc1}
		G_{\mathbf{Q}_{\text{min}}}\big[h\big] = \lim _{n \rightarrow \infty} \mathbb{E}_{\delta_x}\big.\big[ \mathrm{e}^{Z_{n}(\log h)}\ \big|\ {N_{n}>0}\big] = \lim _{n \rightarrow \infty} \frac{ G_n\big[h \big](x) - f_n(0)} {1-f_n(0)},\  x\in\mathbb{R},\ h\in\mathcal{V}_0(\mathbb{R}).
	\end{align}
	\subsection{Proof of Proposition \ref{lem1}}
	\label{sec6.1}
	\noindent\textbf{Proof of Proposition \ref{lem1}}.
	Assume that $\mathbf{Q}$ is a probability measure on $\mathfrak{N}^{\circ}(\mathbb{R})$. Then by Lemma \ref{lem2.2}, $\mathbf{Q}$ is the QSD of the process $Z$ if and only if its p.g.fl. $G_{\mathbf{Q}}$ satisfies
	\begin{align}\label{qsd2}
		(\mathbf{Q}\mathbb{E})\big[\mathrm{e}^{Z_n(\log h)}\ \big|\ N_n>0\big]=G_{\mathbf{Q}}[h],\quad h\in\mathcal{V}_0(\mathbb{R}).
	\end{align}
	Hence, we only need to prove that (\ref{qsd3}) is equivalent to (\ref{qsd2}).
	
	If (\ref{qsd3}) holds, then taking $h\equiv1$ in (\ref{qsd3}) yields that $G_{\mathbf{Q}}[f(0)]=1-m^{\alpha}$. For  $n\in\mathbb{N}$ and  $h\in\mathcal{V}_0(\mathbb{R})$, replacing $h$ by $G_{n-1}[h]$ in (\ref{qsd3}), then we conclude by Proposition \ref{prop1} that 
	\begin{align*}
		G_{\mathbf{Q}} \big[ G_n[h](\cdot) \big]  = m^{\alpha} G_{\mathbf{Q}}\big[G_{n-1}[h](\cdot)\big]+G_{\mathbf{Q}}[f(0)].
	\end{align*}
	By iterating this equality it is easy to see that
	\begin{align}\label{6.3}
		G_{\mathbf{Q}} \big[ G_n[h](\cdot) \big] =m^{\alpha n} G_{\mathbf{Q}}[h]+1-m^{\alpha n} .
	\end{align}
	Therefore, by a combination of (\ref{6.2}) and (\ref{6.3}) we obtain
	\begin{align}\label{6.4}
		(\mathbf{Q} \mathbb{P}) (N_n=0) &= \lim_{\lambda\rightarrow\infty} G_{\mathbf{Q}}[ f_n(\mathrm{e}^{ -\lambda}) ] = m^{\alpha n} \lim_{\lambda\rightarrow\infty} G_{\mathbf{Q}}[\mathrm{e}^{ -\lambda}]+1-m^{\alpha n} \nonumber\\
		&= 1-m^{\alpha n}\in(0,1).
	\end{align}
	Hence,
	\begin{align*}
		(\mathbf{Q}\mathbb{E})\big.\big[\mathrm{e}^{Z_n(\log h)}\ \big|\ {N_n}>0\big] 
		= &\ \frac{(\mathbf{Q} \mathbb{E})\big[\mathrm{e}^{ Z_n(\log h)}\big] - (\mathbf{Q} \mathbb{P})(N_n=0)} {1-(\mathbf{Q} \mathbb{P})(N_n=0)}\\
		= &\ \frac{G_{\mathbf{Q}}\big[G_n[h](\cdot)\big] - (1-m^{\alpha n}) } {m^{\alpha n}},\quad\text{(by (\ref{6.1}) and (\ref{6.4}))}\\
		= &\ G_{\mathbf{Q}}\big[h\big],\quad\text{by (\ref{6.3})},
	\end{align*}
	which implies (\ref{qsd2}) hold.
	
	Suppose that a probability measure $\mathbf{Q}$ on $\mathfrak{N}^{\circ}(\mathbb{R})$ is a QSD of the process $Z$, then 
	\begin{align*}
		(\mathbf{Q} \mathbb{E}) \big[ \mathrm{e}^{Z_n(\log h)} \big] -(\mathbf{Q} \mathbb{P}) (N_n=0) = (\mathbf{Q} \mathbb{P}) (N_n>0)  G_{\mathbf{Q}}[h].
	\end{align*}
	In view of (\ref{6.1}) and (\ref{6.2}) we get
	\begin{align}\label{b8}
		G_{\mathbf{Q}}\big[G_n[h] (\cdot) \big]-G_{\mathbf{Q}}\big[f_n(0)\big]=(\mathbf{Q} \mathbb{P}) (N_n>0)  G_{\mathbf{Q}}\big[h\big].
	\end{align}
	By a similar proof of \cite[Proposition 2]{MV12}, there exists $\alpha\in(0,\infty)$ such that $(\mathbf{Q} \mathbb{P} )(N_n>0)=m^{\alpha n}$. Then (\ref{qsd3}) follows by taking $n=1$ in (\ref{b8}).
	$\hfill\square$
	
	From the proof of Proposition \ref{lem1} we immediately get  the following result.
	\begin{corollary}\label{cor6.1}
		For any $\alpha\in(0,\infty)$, then the probability measure $\mathbf{Q}_{\alpha}$ on $\mathfrak{N}^{\circ}(\mathbb{R})$ is a QSD of $Z$ with eigenvalue $m^{\alpha}$ if and only if its  p.g.fl. $G_{\mathbf{Q}_{\alpha}}$ satisfies the following equation
		\begin{align}\label{qsd5}
			G_{\mathbf{Q}_{\alpha}}\big[G_n[h](\cdot)\big]-G_{\mathbf{Q}_{\alpha}}\big[f_n(0)\big]=m^{\alpha n} G_{\mathbf{Q}_{\alpha}}\big[h \big],\quad n\in\mathbb{N},\quad h\in\mathcal{V}_0(\mathbb{R}).
		\end{align}	
	\end{corollary}
	
	The relationship between the Yaglom limit and the QSD of $Z$ is described as follows.
	\begin{corollary}\label{prop6.2}
		Suppose the $\mathbf{Q}_{\rm min}$ is the Yaglom limit and  $G_{\mathbf{Q}_{\rm min}}$ is its p.g.fl.. Then $G_{\mathbf{Q}_{\rm min}}$ is the solution of the equation {(\ref{qsd3})} for $\alpha=1$, that is, $\mathbf{Q}_{\rm min}$ is a QSD of $Z$ with eigenvalue $m$.
	\end{corollary}
	
	\noindent{\it Proof.} 
	By (\ref{bc1}) we get that for any $x\in\mathbb{R}$,
	\begin{align*}
		G_{\mathbf{Q}_{\text{min}}}\big[h\big] 
		=&\ \lim _{n \rightarrow \infty} \frac{G_{n+1}[h](x)- f_{n+1}(0) }{1-f_{n+1}(0)} \nonumber \\
		=&\ \lim _{n \rightarrow \infty} \bigg(\frac{G_{n+1}[h](x)-f_n(0)}{1-f_n(0)}-\frac{f_{n+1}(0)-f_n(0)}{1-f_n(0)}\bigg) \times \frac{1-f_n(0)}{1-f_{n+1}(0)}.
	\end{align*}
	Note that $f_{n+1}(0)=G_n[f(0)]$ and $(1-f_n(0))/(1-f_{n+1}(0))\rightarrow 1/m$ as $n\rightarrow\infty$. Therefore,
	\begin{align*}
		G_{\mathbf{Q}_{\text{min}}}\big[h\big] 
		=& \frac{1}{m} \times \lim _{n \rightarrow \infty}  \bigg\{ \frac{G_{n+1}[h](x) -f_n(0)}{1-f_n(0)} - \frac{G_n\big[ f(0) \big]-f_n(0)}{1-f_n(0)}  \bigg\} \\	
		=& \frac{1}{m} \times \lim _{n \rightarrow \infty}  \bigg\{ \frac{G_n\big[ G[h](\cdot)\big](x) -f_n(0)}{1-f_n(0)} - \frac{G_n\big[ f(0) \big]-f_n(0)}{1-f_n(0)}  \bigg\},\quad\text{(by Proposition \ref{prop1})} \\
		=& \frac{1}{m} \Big( G_{\mathbf{Q}_{\text{min}}} \big[ G[h](\cdot) \big] - G_{\mathbf{Q}_{\text{min}}} \big[ f(0) \big] \Big),\quad\text{(by (\ref{bc1}))}.
	\end{align*}
	Then we get the desired result.
	$\hfill\square$
	\subsection{Proof of Theorem \ref{QSD}}
	\label{sec6.2}
	\subsubsection{Proof of the necessity part of Theorem \ref{QSD}}
	For any $\alpha\in(0,\infty)$, suppose that a probability measure $\mathbf{Q}_{\alpha}$ on $\mathfrak{N}^{\circ}(\mathbb{R})$ is a QSD of $Z$ with eigenvalue $m^{\alpha}$. We denote by $G_{\mathbf{Q}_{\alpha}}$ its p.g.fl.. Let $q_n:=\mathbb{P}(N_n>0)=1-f_n(0)$. Then by Corollary \ref{cor6.1} we obtain that for any $h\in\mathcal{V}_0(\mathbb{R})$,
	\begin{align*}
		G_{\mathbf{Q}_{\alpha}} [h] &= m^{-\alpha n}  \Big( G_{\mathbf{Q}_{\alpha}}\big[G_n[h](\cdot)\big]-G_{\mathbf{Q}_{\alpha}}\big[f_n(0)\big] \Big)\\
		&= m^{-\alpha n} \Big( (\mathbf{Q}_{\alpha} \mathbb{E}) \big[ \mathrm{e}^{ Z_n (\log h)} \big] - (\mathbf{Q}_{\alpha} \mathbb{P}) (N_n=0) \Big),\quad\text{(by (\ref{6.1}) and (\ref{6.2}))} \\
		&= m^{-\alpha n}\int_{\mathfrak{N}^{\circ}(\mathbb{R})} \Big( \mathbb{E}_\mu\big[\mathrm{e}^{ Z_n ( \log h)}\big] - \mathbb{P}_{\mu}(N_n=0)  \Big) \mathbf{Q}_\alpha(\mathrm{d} \mu) \\
		&= m^{-\alpha n}\int_{\mathfrak{N}^{\circ}(\mathbb{R})} \bigg[ \exp \bigg(\int_{\mathbb{R}} \log G_n\big[h\big] (x) \mu(\mathrm{d} x)\bigg) - (1-q_n)^{\mu(\mathbb{R})}  \bigg] \mathbf{Q}_\alpha(\mathrm{d} \mu),
	\end{align*}
	where the last equality comes from the independence. For any $n\in\mathbb{N}$ we can define the measure $\gamma_n$ by $\gamma_n(B) = m^{-\alpha n} \mathbf{Q}_\alpha(q_n^{-1} B)$ for $B\in\mathcal{B}(\mathfrak{M}^{\circ}(\mathbb{R}))$. Then a simple calculation shows that
	\begin{align}\label{6.7}
		G_{\mathbf{Q}_{\alpha}} \big[h\big] 
		&= \int_{\mathfrak{M}^{\circ}(\mathbb{R})} \bigg[\exp \bigg(q_n^{-1} \int_{\mathbb{R}} \log G_n[h](x) \nu(\mathrm{d} x)\bigg)-(1-q_n)^{q_n^{-1} \nu(\mathbb{R})}\bigg] \gamma_n(\mathrm{d} \nu) \nonumber \\
		&= \int_{\mathfrak{M}^{\circ}(\mathbb{R})} U_{h,n}(\nu) \widetilde{\gamma}_n(\mathrm{d} \nu),
	\end{align}
	where $U_{h,n}$ is a functional on  $\mathfrak{M}^{\circ}(\mathbb{R})$ defined by
	\begin{align}\label{6.8}
		U_{h, n}(\nu):=\frac{ \exp \big(q_n^{-1} \int_{\mathbb{R}} \log G_n[h](x) \nu(\mathrm{d} x)\big)-(1-q_n)^{q_n^{-1} \nu(\mathbb{R})} }{\nu(\mathbb{R}) \mathrm{e}^{-\nu(\mathbb{R})}}
	\end{align}
	and $\widetilde{\gamma}_n$ is a measure on $\mathfrak{M}(\mathbb{R})$ given by 
	\begin{align}\label{6.9}
		\widetilde{\gamma}_n(\mathrm{d} \nu):=\left\{\begin{array}{ll}
			\nu(\mathbb{R}) \mathrm{e}^{-\nu(\mathbb{R})} \gamma_n(\mathrm{d} \nu) & \text{if}\ \nu \in \mathfrak{M}^{\circ}(\mathbb{R}), \\
			0 & \text{if}\ \nu=\mathbf{0}. 
		\end{array}\right.
	\end{align}
	In the proof of the necessity part of Theorem \ref{QSD}, we will use the following three propositions whose proofs are postponed to Section \ref{sec6.3}.
	
	\begin{proposition}\label{prop6.3}
		Assume that {\rm(C1)}-{\rm(C3)} hold. Suppose that $U_{h,n}$ is defined by {(\ref{6.8})}. Then for any $h\in\mathcal{V}_0(\mathbb{R})$ we have
		\begin{align}\label{6.26}
			\lim_{n\rightarrow\infty}U_{h,n}(\nu)= U_{h,\infty} (\nu)
			:=  \frac{\mathrm{e}^{G_{\mathbf{Q}_{\rm min}}[h] \nu(\mathbb{R})}-1}{\nu(\mathbb{R}) } ,\quad \nu\in\mathfrak{M}^{\circ}(\mathbb{R}),
		\end{align}
		where $\mathbf{Q}_{\rm min}$ is the Yaglom limit and its p.g.fl. is $G_{\mathbf{Q}_{\rm min}}$.
	\end{proposition}
	\begin{proposition}\label{prop6.4}
		Suppose that $\sum_{j} (j\log j)p_j<\infty$. Let  $\widetilde{\gamma}_n$ be a measure on $\mathfrak{M}(\mathbb{R})$ defined by {(\ref{6.9})}, then the sequence $\{\widetilde{\gamma}_n:n\geq 0\}$ is tight in $\mathfrak{M}(\mathbb{R})$ and 
		\begin{align}\label{6.23}
			\sup_{n\in\mathbb{N}} \int_{\mathfrak{M}(\mathbb{R})}  \frac{ \mathrm{e}^{ \nu(\mathbb{R})} - 1 }{\nu(\mathbb{R}) } \widetilde{\gamma}_n(\mathrm{d} \nu) <\infty.
		\end{align} 
	\end{proposition}
	\begin{proposition}\label{prop6.5}
		Assume that {\rm(C1)}, {\rm(C2)} and {\rm(C4)} hold. Let $G_n$ be the p.g.fl. of $Z_n$ defined by (\ref{4}). Suppose that  $\mathbf{Q}_{\rm min}$ is the Yaglom limit and $G_{\mathbf{Q}_{\rm min}}$ is its p.g.fl.. Then we have
		\begin{align*}
			\lim_{n\rightarrow\infty} \frac{G_n[h](x) - f_n(0)}{1-f_n(0)} = G_{\mathbf{Q}_{\rm min}} [h]
		\end{align*} 
		uniformly in $x\in\mathbb{R}$ for each $h\in\mathcal{V}_0(\mathbb{R})$.
	\end{proposition}
	
	\begin{remark}
		In Proposition \ref{prop6.5}, under the strengthened condition {\rm(C4)}, it then follows by using the size-biased measure that the pointwise convergence in Yaglom's theorem is raised to uniform convergence in $x$.  We point out that Proposition \ref{prop6.5} is a key step in giving the integral characterization {(\ref{6.19})} of the p.g.fl. of the QSD with eigenvalue $m^{\alpha}$.
	\end{remark}
	
	\noindent\textbf{Proof of the necessity of Theorem \ref{QSD}.}
	Since $\mathbb{R}$ is a Polish space, so is $\mathfrak{M}(\mathbb{R})$; see, e.g., \cite[Lemma 4.3]{K17}. Then by Proposition \ref{prop6.4} and Prokhorov’s theorem, there exists a subsequence $\{\widetilde{\gamma}_{n_k}\}\subset\{\widetilde{\gamma}_{n}\}$ such that $\widetilde{\gamma}_{n_k}$ converges weakly to a finite measure $\widetilde{\gamma}$ on $\mathfrak{M}(\mathbb{R})$. 
	For any $h\in\mathcal{V}_0(\mathbb{R})$, we claim that
	\begin{align}\label{6.19}
		G_{\mathbf{Q}_{\alpha}}[h] 
		= \int_{\mathfrak{M}^{\circ}(\mathbb{R})} U_{h,\infty} (\nu) \widetilde{\gamma} (\mathrm{d} \nu) + G_{\mathbf{Q}_{\rm min}} [h] \widetilde{\gamma} (\{\mathbf{0}\}),
	\end{align}
	where $U_{h,\infty}$ is defined by (\ref{6.26}). Choose $L\geq1$ with $\widetilde{\gamma}(\{\nu\in\mathfrak{M}(\mathbb{R}):\nu(\mathbb{R})=L\})=0$. It then follows from (\ref{6.7}) that
	\begin{align}\label{6.20}
		G_{\mathbf{Q}_{\alpha}} [h] 
		=&\ \int_{\mathfrak{M}^{\circ}(\mathbb{R})} U_{h,n_k}(\nu) \widetilde{\gamma}_{n_k}(\mathrm{d} \nu) \nonumber \\
		=&\  \int_{ \{\nu\in\mathfrak{M}^{\circ}(\mathbb{R}): \nu(\mathbb{R})\leq L\} } U_{h,n_k}(\nu) \widetilde{\gamma}_{n_k}(\mathrm{d} \nu) +  \int_{ \{\nu\in\mathfrak{M}(\mathbb{R}): \nu(\mathbb{R})>L\} } U_{h,n_k}(\nu) \widetilde{\gamma}_{n_k}(\mathrm{d} \nu).  
	\end{align}
	Let $\mathscr{L}:=\{h\in\mathcal{V}_0(\mathbb{R}): \sup_{x\in\mathbb{R}} h(x)<1\}$. We first prove (\ref{6.19}) holds for $h\in\mathscr{L}$ through the following three steps.

	{\it Step 1.} For $h\in\mathscr{L}$, we shall prove that the second term on the right-hand side of (\ref{6.20}) tends to $0$ as $k\rightarrow\infty$ first and then $L\rightarrow\infty$. Note that 
	\begin{align}\label{6.41}
		G_{n_k}[h] (x) = 1-q_{n_k}\Big(1-\frac{G_{n_k}[h] (x)-(1-q_{n_k})}{q_{n_k}}\Big).
	\end{align}
	Since  $\log (1-x)\leq -x$ for $x<1$, we have
	\begin{align*}
		q_{n_k}^{-1} \log G_{n_k}[h](x) \leq \frac{G_{n_k}[h] (x)-(1-q_{n_k})}{q_{n_k}} -1,\quad x\in\mathbb{R}.
	\end{align*}
	Let $\bar{c}_1:=\sup_{x\in\mathbb{R}} h(x)$, then $h\in\mathscr{L}$ implies $0\leq h(x)\leq \bar{c}_1<1$ for all $x\in\mathbb{R}$. Hence,  by (\ref{13}) it is easy to see that when $\sum(j \log j)p_j<\infty$, we have
	\begin{align*}
		\frac{G_{n_k}[h] (x)-(1-q_{n_k})}{q_{n_k}} \leq \frac{ f_{n_k}(\bar{c}_1)-f_{n_k}(0) }{ 1-f_{n_k}(0) } \xrightarrow[k\rightarrow\infty]{} 1-\bar{c}_2< 1,
	\end{align*}
	where $\bar{c}_2:=\frac{(1-\bar{c}_1)\varphi(\bar{c}_1)}{\varphi(0)}>0$. From the preceding two inequalities one can choose sufficiently large $k\geq1$ to obtain $q_{n_k}^{-1} \log G_{n_k}[h](x) \leq -\bar{c}_2/2$. On the other hand, for sufficiently large $k\geq1$, there exists $\underline{c}\in(0,1)$ such that $(1-q_{n_k})^{q_{n_k}^{-1} \nu(\mathbb{R})} \geq \mathrm{e}^{-(1+\underline{c}) \nu(\mathbb{R})}$. Consequently,
	\begin{align*}
		U_{h,n_k}(\nu) &= \frac{ \exp \big(q_{n_k}^{-1} \int_{\mathbb{R}} \log G_{n_k}[h](x) \nu(\mathrm{d} x)\big)-(1-q_{n_k})^{q_{n_k}^{-1} \nu(\mathbb{R})} }{\nu(\mathbb{R}) \mathrm{e}^{-\nu(\mathbb{R})}} \\ &\leq \frac{\mathrm{e}^{-\bar{c}_2\nu(\mathbb{R})/2}-\mathrm{e}^{-(1+\underline{c}) \nu(\mathbb{R})}}{\nu(\mathbb{R}) \mathrm{e}^{-\nu(\mathbb{R})}}.
	\end{align*}
	It then follows that
	\begin{align}\label{6.22}
		\int_{ \{\nu\in\mathfrak{M}(\mathbb{R}): \nu(\mathbb{R})>L\} }\!\! U_{h,n_k}(\nu) \widetilde{\gamma}_{n_k}(\mathrm{d} \nu) 
		\leq&\ \int_{ \{\nu\in\mathfrak{M}(\mathbb{R}): \nu(\mathbb{R})>L\} } \frac{\mathrm{e}^{-\bar{c}_2\nu(\mathbb{R})/2 }-\mathrm{e}^{-(1+\underline{c}) \nu(\mathbb{R})}}{\nu(\mathbb{R}) \mathrm{e}^{-\nu(\mathbb{R})}} \widetilde{\gamma}_{n_k}(\mathrm{d} \nu)\nonumber \\
		=&\ \int_{ \{\nu\in\mathfrak{M}(\mathbb{R}): \nu(\mathbb{R})>L\} } \Big( \mathrm{e}^{-\bar{c}_2\nu(\mathbb{R})/2}-\mathrm{e}^{-(1+\underline{c}) \nu(\mathbb{R})} \Big) \gamma_{n_k}(\mathrm{d} \nu).
	\end{align}
	For any $n\in\mathbb{N}$, we define two measures $\widetilde{\kappa}_n$ and $\widetilde{\kappa}$  on $[0,\infty)$  by $\widetilde{\kappa}_n(\mathrm{d} x) := \widetilde{\gamma}_n(\{\nu\in\mathfrak{M}(\mathbb{R}):\nu(\mathbb{R})\in\mathrm{d} x\})$ and $\widetilde{\kappa}(\mathrm{d} x) := \widetilde{\gamma}(\{\nu\in\mathfrak{M}(\mathbb{R}):\nu(\mathbb{R})\in\mathrm{d} x\})$, respectively. Let $\kappa_n$ and $\kappa$ be two measures on $(0,\infty)$ given by $\kappa_n(\mathrm{d} x) := x^{-1} \mathrm{e}^{x} \widetilde{\kappa}_n |_{(0,\infty)}(\mathrm{d} x)$ and $\kappa(\mathrm{d} x):=x^{-1} \mathrm{e}^x \widetilde{\kappa}|_{(0, \infty)}(\mathrm{d} x)$, where the symbol $\widetilde{\kappa}|_{(0, \infty)}$ is the restriction of $\widetilde{\kappa}$ to $(0,\infty)$. As a consequence, (\ref{6.22}) can be written as
	\begin{align}\label{6.40}
		\int_{ \{\nu\in\mathfrak{M}(\mathbb{R}): \nu(\mathbb{R})>L\} }\!\! U_{h,n_k}(\nu) \widetilde{\gamma}_{n_k}(\mathrm{d} \nu)  \leq \int_{(L,\infty)} \Big( \mathrm{e}^{-\bar{c}_2 x/2}-\mathrm{e}^{-(1+\underline{c}) x} \Big) \kappa_{n_k}(\mathrm{d} x).
	\end{align}
	Note that $\widetilde{\gamma}_{n_k}$ converges weakly to $\widetilde{\gamma}$ as $k\rightarrow\infty$ implies $\widetilde{\kappa}_{n_k}\xrightarrow[]{\text{\rm w}}\widetilde{\kappa}$. Then it is natural to get that $\kappa_{n_k}$ converges vaguely on $(0,\infty)$ to $\kappa$. Therefore we have 
	\begin{align*}
		\lim_{k\rightarrow\infty} \int_{(L,\infty)} \Big( \mathrm{e}^{-\bar{c}_2 x/2}-\mathrm{e}^{-(1+\underline{c}) x} \Big) \kappa_{n_k}(\mathrm{d} x) = \int_{(L,\infty)} \Big( \mathrm{e}^{-\bar{c}_2 x/2}-\mathrm{e}^{-(1+\underline{c}) x} \Big) \kappa(\mathrm{d} x).
	\end{align*}
	From (\ref{6.23}) we see that $\sup_n \int_{(0,\infty)} (1-\mathrm{e}^{-x}) \kappa_n(\mathrm{d}x)<\infty$. Then the right-hand side of the above equality tends to zero as $L\rightarrow\infty$. Hence the desired result follows from (\ref{6.40}).
	
	{\it Step 2}. For $h\in\mathscr{L}$, by a similar reasoning as in step $1$ we get
	\begin{align*}
		\int_{ \{\nu\in\mathfrak{M}(\mathbb{R}): \nu(\mathbb{R})>L\} } U_{h,\infty}(\nu) \widetilde{\gamma}(\mathrm{d} \nu) \leq \int_{ (L,\infty) }  \Big( \mathrm{e}^{(G_{\mathbf{Q}_{\text{min}}}[\bar{c}_1]-1) x} -\mathrm{e}^{-x} \Big) \kappa(\mathrm{d} x) \xrightarrow[L\rightarrow\infty]{}0.
	\end{align*} 
	
	{\it Step 3}. For $h\in\mathscr{L}$, by 
		setting 
		\begin{align*}
			U_{h,\infty}(\mathbf{0}):=\lim_{\nu\rightarrow\mathbf{0}}U_{h,\infty}(\nu)=G_{\mathbf{Q}_{\text{min}}}[h],
		\end{align*}
		where $\nu\rightarrow\mathbf{0}$ means that $\nu$ converges weakly to $\mathbf{0}$, then the functional $U_{h,\infty}$ on $\mathfrak{M}^{\circ}(\mathbb{R})$ can be continuously extended (with respect to the weak topology) to $\mathfrak{M}(\mathbb{R})$. Based on the steps $1$ and $2$, to get (\ref{6.19}) we only need to show that
	\begin{align}\label{6.25}
		\int_{ \{\nu\in\mathfrak{M}^{\circ}(\mathbb{R}): \nu(\mathbb{R})\leq L\} } U_{h,n_k}(\nu) \widetilde{\gamma}_{n_k}(\mathrm{d} \nu) \xrightarrow[k\rightarrow\infty]{} \int_{ \{\nu\in\mathfrak{M}(\mathbb{R}): \nu(\mathbb{R})\leq L\} } U_{h,\infty} (\nu) \widetilde{\gamma} (\mathrm{d} \nu) 
	\end{align}
	holds for every $L\geq1$ with  $\widetilde{\gamma}(\{\nu\in\mathfrak{M}(\mathbb{R}):\nu(\mathbb{R})=L\})=0$. Recall that $\widetilde{\gamma}_{n_k}(\{\mathbf{0}\})=0$ for any $k\geq0$. Then by a simple calculation we have
	\begin{align*}
		\bigg| \int_{ \{\nu\in\mathfrak{M}^{\circ}(\mathbb{R}): \nu(\mathbb{R})\leq L\} } U_{h,n_k}(\nu) \widetilde{\gamma}_{n_k}(\mathrm{d} \nu) - \int_{ \{\nu\in\mathfrak{M}(\mathbb{R}): \nu(\mathbb{R})\leq L\} } U_{h,\infty} (\nu) \widetilde{\gamma} (\mathrm{d} \nu) \bigg|  \leq J_{1,k} + J_{2,k},
	\end{align*}
	where
	\begin{align*}
		J_{1,k} :=	\bigg| \int_{ \{\nu\in\mathfrak{M}^{\circ}(\mathbb{R}): \nu(\mathbb{R})\leq L\} } \Big[ U_{h,n_k}(\nu) \widetilde{\gamma}_{n_k}(\mathrm{d} \nu) -  U_{h,\infty} (\nu) \Big] \widetilde{\gamma}_{n_k} (\mathrm{d} \nu) \bigg| 
	\end{align*}
	and
	\begin{align*}
		J_{2,k} :=	\bigg| \int_{ \{\nu\in\mathfrak{M}(\mathbb{R}): \nu(\mathbb{R})\leq L\} } U_{h,\infty}(\nu) \widetilde{\gamma}_{n_k}(\mathrm{d} \nu) - \int_{ \{\nu\in\mathfrak{M}(\mathbb{R}): \nu(\mathbb{R})\leq L\} } U_{h,\infty} (\nu) \widetilde{\gamma} (\mathrm{d} \nu) \bigg| .
	\end{align*}
	
	{\it (1) The asymptotic behavior of $J_{2,k}$}. By the definition of $U_{h,\infty}$ it is easy to show 
	\begin{align*}
		J_{2,k} \leq \sup _{0 \leq x \leq L}\bigg|\frac{\mathrm{e}^{G_{\mathbf{Q}_{\text{min}}}[h] x}-1}{x}\bigg| \times \bigg|\widetilde{\kappa}_{n_k}\big([0, L]\big)-\widetilde{\kappa}\big([0, L]\big)\bigg|.
	\end{align*}
	Therefore $\lim_{k\rightarrow\infty}J_{2,k}=0$ follows by the weak convergence of $\widetilde{\kappa}_{n_k}$. 
	
	{\it (2) The asymptotic behavior of $J_{1,k}$}. From the definitions of  $U_{h,n}(\nu)$ and $U_{h,\infty}(\nu)$, then for any $\nu\in\mathfrak{M}^{\circ}(\mathbb{R})$ we have
	\begin{align}\label{6.28}
		&\ U_{h,n_k}(\nu)-U_{h, \infty}(\nu)  \\
		=&\ \frac{\exp\! \big\{q_{n_k}^{-1} \!\! \int_{\mathbb{R}} \log G_{n_k}[h] (x) \nu(\mathrm{d} x)\big\} \!\! - \! \exp\! \big\{(G_{\mathbf{Q}_{\text{min}}}[h]-1) \nu(\mathbb{R})\big\} } {\nu(\mathbb{R}) \mathrm{e}^{-\nu(\mathbb{R})}} \! + \! \frac{(1-q_{n_k})^{q_{n_k}^{-1} \nu(\mathbb{R})} \!-\! \mathrm{e}^{-\nu(\mathbb{R})} }{\nu(\mathbb{R}) \mathrm{e}^{-\nu(\mathbb{R})}}\nonumber.
	\end{align}
	For the second term on the right-hand side of (\ref{6.28}), note that
	\begin{small}
		\begin{align}\label{6.29}
			\int_{ \{\nu\in\mathfrak{M}^{\circ}(\mathbb{R}): \nu(\mathbb{R})\leq L\} } \! \frac{(1-q_{n_k})^{q_{n_k}^{-1} \nu(\mathbb{R})} \!-\! \mathrm{e}^{-\nu(\mathbb{R})} }{\nu(\mathbb{R}) \mathrm{e}^{-\nu(\mathbb{R})}} \widetilde{\gamma}_{n_k} (\mathrm{d} \nu)  =  \int_{ (0,L] } \!\! \frac{ (1-q_{n_k})^{q_{n_k}^{-1} x} \!-\! \mathrm{e}^{-x} }{x \mathrm{e}^{-x}} \widetilde{\kappa}_{n_k}\! ( \mathrm{d} x).
		\end{align}
	\end{small}We denote by $I_{n_k}(x)$ the integrand on the right-hand side of the above equality.  Then the function $x \mapsto I_{n_k}(x)$, continuously extended to $[0,L]$ by setting $I_{n_k}(0)= 1+q_{n_k}^{-1} \log (1-q_{n_k}) $,  converges uniformly on every compact subset of $[0,L]$ to $0$. Hence, by Proposition \ref{prop6.4} we have
	\begin{align}\label{6.30}
		\bigg| \int_{(0,L]}\!\! I_{n_k}(x) \widetilde{\kappa}_{n_k} (\mathrm{d} x) \bigg|\! =\! \bigg| \int_{[0,L]}\!\! I_{n_k}(x) \widetilde{\kappa}_{n_k} (\mathrm{d} x) \bigg| 
		\!\leq\! \sup_{x\in[0,L]} \big| I_{n_k}(x) \big|\! \times\! \sup_{k\in\mathbb{N}} \widetilde{\gamma}_{n_k} (\mathfrak{M}(\mathbb{R}))  \xrightarrow[k\rightarrow\infty]{} 0.
	\end{align}
	
	For the first term on the right-hand side of (\ref{6.28}), based on (\ref{6.41}), by using the inequality $|\mathrm{e}^x-\mathrm{e}^y|\leq |x-y|$ for $x,y\leq0$ it is easily seen that
	\begin{align}\label{6.31}
		&\ \bigg| \exp\! \Big\{q_{n_k}^{-1} \!\! \int_{\mathbb{R}} \log G_{n_k}[h] (x) \nu(\mathrm{d} x)\Big\} \!\! - \! \exp\! \Big\{(G_{\mathbf{Q}_{\text{min}}}[h]-1) \nu(\mathbb{R})\Big\}  \bigg| \nonumber \\
		\leq &\ \bigg| q_{n_k}^{-1} \int_{\mathbb{R}} \log \Big[1-q_{n_k}\Big(1-\frac{\!G_{n_k}[h](x)\!-\!(1-q_{n_k})}{q_{n_k}}\Big)\Big] \nu(\mathrm{d} x) - (G_{\mathbf{Q}_{\text{min}}}[h]-1) \nu(\mathbb{R})  \bigg| \nonumber \\
		\leq &\ \bigg| q_{n_k}^{-1} \int_{\mathbb{R}} \Big\{ \log \!\Big[ \! 1\!-\!q_{n_k}\Big(1\!-\!\frac{G_{n_k}[h](x)\!-\!(1-q_{n_k})}{q_{n_k}}\Big)\Big] \! + \!\Big(\!1-\!\frac{G_{n_k}[h] (x)\!-\!(1-q_{n_k})}{q_{n_k}}\Big) \Big\} \nu(\mathrm{d} x) \bigg| \nonumber \\
		&\ + \bigg| \int_{\mathbb{R}} \Big[ \frac{G_{n_k}[h](x)-(1-q_{n_k})}{q_{n_k}} - G_{\mathbf{Q}_{\text{min}}}[h] \Big] \nu(\mathrm{d} x)  \bigg|.
	\end{align}
	Note that $q_{n_k}\rightarrow0$ as $k\rightarrow\infty$ and $q_{n_k}^{-1}(G_{n_k}[h]-(1-q_{n_k}))\in[0,1]$ for $k\in\mathbb{N}$. Then by Taylor’s expansion of $\log (1-x)$ we see that for sufficiently large $k\geq 1$,
	\begin{align*}
		\log \Big[1-q_{n_k}\Big(1-\frac{G_{n_k}[h](x)-(1-q_{n_k})}{q_{n_k}}\Big)\Big] = -\sum_{i=1}^{\infty} \frac{q_{n_k}^i}{i}\Big(1-\frac{G_{n_k}[h] (x)-(1-q_{n_k})}{q_{n_k}}\Big)^i.
	\end{align*}
	It then follows that
	\begin{align*}
		\Big| \log \!\Big[ \! 1\!-\!q_{n_k}\Big(1\!-\!\frac{G_{n_k}[h](x)\!-\!(1-q_{n_k})}{q_{n_k}}\Big)\Big] \! + \!\Big(\!1-\!\frac{G_{n_k}[h] (x)\!-\!(1-q_{n_k})}{q_{n_k}}\Big) \Big| \leq \sum_{i=2}^{\infty} q_{n_k}^i.
	\end{align*}
	Consequently, 
	\begin{align}\label{6.32}
		&\ \bigg| q_{n_k}^{-1} \int_{\mathbb{R}} \Big\{ \log \!\Big[ \! 1\!-\!q_{n_k}\Big(1\!-\!\frac{G_{n_k}[h](x)\!-\!(1-q_{n_k})}{q_n}\Big)\Big] \! + \!\Big(\!1-\!\frac{G_{n_k}[h] (x)\!-\!(1-q_{n_k})}{q_{n_k}}\Big) \Big\} \nu(\mathrm{d} x) \bigg| \nonumber \\
		\leq &\ \int_{\mathbb{R}} \sum_{i=2}^{\infty} q_{n_k}^{i-1} \nu(\mathrm{d} x) =  \frac{q_{n_k}}{1-q_{n_k}} \nu(\mathbb{R}).
	\end{align}
	This combined with (\ref{6.31}) gives
	\begin{align*}
		\bigg| \exp\! \Big\{q_{n_k}^{-1} \!\! \int_{\mathbb{R}} \log G_{n_k}[h] (x) \nu(\mathrm{d} x)\Big\} \!\! - \! \exp\! \Big\{(G_{\mathbf{Q}_{\text{min}}}[h]-1) \nu(\mathbb{R})\Big\}  \bigg|  
		\leq \Big( \bar{a}_{n_k}(h) + \frac{q_{n_k}}{1-q_{n_k}} \Big)\nu(\mathbb{R}),
	\end{align*}
	where 
	\begin{align*}
		\bar{a}_{n_k}(h):=\sup_{x\in\mathbb{R}}\bigg| \frac{G_{n_k}[h](x)-f_{n_k}(0)}{1-f_{n_k}(0)} - G_{\mathbf{Q}_{\text{min}}}[h] \bigg|.
	\end{align*}
	By Proposition \ref{prop6.5} we see that $\lim_{k\rightarrow\infty}\bar{a}_{n_k}(h)=0$ for every $h\in\mathcal{V}_0(\mathbb{R})$. Proposition \ref{prop6.4} implies that $\sup_{n\in\mathbb{N}} \widetilde{\gamma}_{n} (\mathfrak{M}(\mathbb{R})) <\infty$. It then follows that
	\begin{align*}
		&\ \bigg| \int_{ \{\nu\in\mathfrak{M}^{\circ}(\mathbb{R}): \nu(\mathbb{R})\leq L\} } \frac{ \exp\! \Big\{q_{n_k}^{-1} \!\! \int_{\mathbb{R}} \log G_{n_k}[h] (x) \nu(\mathrm{d} x)\Big\} \!\! - \! \exp\! \Big\{(G_{\mathbf{Q}_{\text{min}}}[h]-1) \nu(\mathbb{R})\Big\} }{ \nu(\mathbb{R}) \mathrm{e}^{-\nu(\mathbb{R})} } \widetilde{\gamma}_{n_k} (\mathrm{d} \nu) \bigg| \\
		\leq&\ \mathrm{e}^L  \Big( \bar{a}_{n_k}(h) + \frac{q_{n_k}}{1-q_{n_k}} \Big) \times \sup_{k\in\mathbb{N}} \widetilde{\gamma}_{n_k} (\mathfrak{M}(\mathbb{R})) \xrightarrow[k\rightarrow\infty]{} 0. 
	\end{align*}
	Then (\ref{6.28})-(\ref{6.30}) and the above inequality imply $J_{1,k}\rightarrow0$ as $k\rightarrow\infty$. Summing up, by the assertions {\rm(1)} and {\rm(2)}, (\ref{6.25}) holds for any $h\in\mathscr{L}$. 
	
	Given $h\in\mathcal{V}_0(\mathbb{R})$ we take $h_n(x):=h(x)\mathbf{1}_{\{h(x)<1-\frac{1}{n}\}}+(1-\frac{1}{n})\mathbf{1}_{\{1-\frac{1}{n}\leq h(x)\leq1\}}$ so that $h_n\in\mathscr{L}$ and $h_n(x)\uparrow h(x)$ as $n\rightarrow\infty$. Then (\ref{6.19}) holds for $h\in\mathcal{V}_0(\mathbb{R})$ by monotone convergence. Therefore, by the definition of $\kappa$ and $\widetilde{\kappa}$, (\ref{6.19}) can be written as
	\begin{align}\label{6.27}
		G_{\mathbf{Q}_{\alpha}}[h] = \int_{(0,\infty)} \big( \mathrm{e}^{(G_{\mathbf{Q}_{\text{min}}}[h]-1) x}-\mathrm{e}^{-x} \big) \kappa(\mathrm{d} x)+\widetilde{\kappa}(\{0\}) G_{\mathbf{Q}_{\text{min}}}[h],\quad h\in\mathcal{V}_0(\mathbb{R}).
	\end{align}
	Furthermore, the theory of Laplace transforms gives that $\kappa$ and $\widetilde{\kappa}(\{0\})$, hence $\widetilde{\kappa}$ is uniquely determined by (\ref{6.27}), so that $\widetilde{\kappa}_{n}$ converges weakly to $\widetilde{\kappa}$. Consequently, $\kappa_{n}$ converges vaguely to $\kappa$ on $(0,\infty)$.
	
	The property of the measure $\kappa$ follows from this convergence: for every compact interval $B\subset(0,\infty)$ whose endpoints are not atoms of $\kappa$,
	\begin{align*}
		\kappa(B)=&\ \lim _{n \rightarrow \infty} \gamma_{n+1}\big(\{\nu \in \mathfrak{M}^{\circ}(\mathbb{R}): \nu(\mathbb{R}) \in B \} \big) \\
		=&\ \lim _{n \rightarrow \infty} m^{-\alpha(n+1)} \mathbf{Q}_{\alpha} \big( \{\nu\in\mathfrak{N}^{\circ}(\mathbb{R}): \nu(\mathbb{R}) \in q_{n+1}^{-1}B\} \big) \\
		=&\ m^{-\alpha} \lim _{n \rightarrow \infty} m^{-\alpha n} \mathbf{Q}_{\alpha} \big( \{\nu\in\mathfrak{N}^{\circ}(\mathbb{R}): \nu(\mathbb{R}) \in q_{n}^{-1} \big(q_n/q_{n+1}\big) B\} \big)  \\
		=&\ m^{-\alpha} \lim _{n \rightarrow \infty} m^{-\alpha n} \gamma_{n}\big(\{\nu \in \mathfrak{M}^{\circ}(\mathbb{R}): \nu(\mathbb{R}) \in \big(q_n/q_{n+1}\big)B \} \big) \\
		=&\ m^{-\alpha} \kappa(m^{-1} B),\quad\text{(by $q_n/q_{n+1}\rightarrow 1/m$)}.
	\end{align*}
	This implies that the measure $\Delta$ defined by $\Delta(\mathrm{d} x):=x^\alpha \kappa(\mathrm{d} x)$ satisfies $\Delta(m B) =\Delta(B)$ for every $B\in\mathcal{B}((0,\infty))$.
	
	It remains to investigate which terms in (\ref{6.27}) vanish for particular values of $\alpha$.  A first constraint comes from the fact that $\mathbf{Q}_{\alpha}$ is a probability on $\mathfrak{N}^{\circ}(\mathbb{R})$, then $G_{\mathbf{Q}_{\alpha}}[1]=1$.
	
	A second constraint comes from the fact that $G_{\mathbf{Q}_{\alpha}}[h]\leq 1<\infty$ for $h\in\mathcal{V}_0(\mathbb{R})$, then the integral in (\ref{6.27}) needs to be finite as well, that is
	\begin{align*}
		\int_{(0,\infty)}\big(\mathrm{e}^{(G_{\mathbf{Q}_{\text{min}}}[h]-1) x}-\mathrm{e}^{-x}\big) \kappa(
		\mathrm{d} x)<\infty.
	\end{align*} 
	Checking the proof of \cite[Corollary 16]{M18}, one can show that this corollary still holds when the generating function $H$ is replaced by the p.g.fl. $G_{\mathbf{Q}_{\text{min}}}$. Then the above inequality holds if and only if $\alpha<1$ or $\kappa=0$.
	
	A third constraint comes from the fact that $G_{\mathbf{Q}_{\alpha}}$ satisfies the equation (\ref{qsd3}). We first recall the following equation for $G_{\mathbf{Q}_{\text{min}}}$:
	\begin{align}
		G_{\mathbf{Q}_{\text{min}}}\big[G[h](\cdot)\big]-G_{\mathbf{Q}_{\text{min}}}\big[f(0)\big] & =m G_{\mathbf{Q}_{\text{min}}}\big[h\big], \quad h\in\mathcal{V}_0(\mathbb{R}),\label{6.33} \\
		G_{\mathbf{Q}_{\text{min}}}\big[f(0)\big] & =1-m,\label{6.34} \\
		G_{\mathbf{Q}_{\text{min}}}\big[G[h](\cdot)\big]-1 & =m(G_{\mathbf{Q}_{\text{min}}}\big[h\big]-1), \quad h\in\mathcal{V}_0(\mathbb{R}) .\label{6.35}
	\end{align} 
	Indeed, (\ref{6.33}) comes from Corollary \ref{prop6.2}. (\ref{6.34}) follows from (\ref{6.33}) by setting $h\equiv1$. (\ref{6.35}) comes from (\ref{6.33}) and (\ref{6.34}) by reordering terms. Then by (\ref{6.27}), for any $h\in\mathcal{V}_0(\mathbb{R})$,
	\begin{align*}
		&\ G_{\mathbf{Q}_{\alpha}}\big[G[h](\cdot)\big]-G_{\mathbf{Q}_{\alpha}}\big[f(0)\big] \\
		=&\ \int_{(0,\infty)}\!\!\Big[\mathrm{e}^{(G_{\mathbf{Q}_{\text{min}}}[G[h](\cdot)]-1) x}\!-\!\mathrm{e}^{(G_{\mathbf{Q}_{\text{min}}}[f(0)]-1) x}\!\Big] \kappa(\mathrm{d} x)\!+\!\widetilde{\kappa}(0)\Big(\!G_{\mathbf{Q}_{\text{min}}}\!\big[G[h](\cdot)\big]\!-\!G_{\mathbf{Q}_{\text{min}}}[f(0)]\Big) \\
		=&\ \int_{(0,\infty)}\Big(\mathrm{e}^{m(G_{\mathbf{Q}_{\text{min}}}[h]-1) x}-\mathrm{e}^{-m x}\Big) \kappa(\mathrm{d} x)+m \widetilde{\kappa}(\{0\}) G_{\mathbf{Q}_{\text{min}}}[h]\quad \text{(by (\ref{6.33})-(\ref{6.35})) } \\
		=&\ m^\alpha \int_{(0,\infty)}\Big(\mathrm{e}^{(G_{\mathbf{Q}_{\text{min}}}[h]-1) x}-\mathrm{e}^{-x}\Big) \kappa(\mathrm{d}x)+m \widetilde{\kappa}(\{0\}) G_{\mathbf{Q}_{\text{min}}}[h]\quad \text{ (by $\kappa(B)=m^{-\alpha} \kappa(m^{-1}B)$)}.
	\end{align*}
	Comparing with (\ref{qsd3}) implies that $\widetilde{\kappa}(\{0\})=0$ unless $\alpha=1$. Summing up, we have the following constraints for (\ref{6.27}):
	\begin{itemize}
		\item $\alpha>1$: there exist no QSDs of $Z$ with eigenvalue $m^{\alpha}$,
		\item $\alpha=1: \kappa=0$ and $\widetilde{\kappa}(\{0\})=1$,
		\item $\alpha<1: \widetilde{\kappa}(\{0\})=0$ and $\int_{(0,\infty)} (1-\mathrm{e}^{x}) \kappa(\mathrm{d} x)=1$.
	\end{itemize}
	This proves the necessity part of the Theorem \ref{QSD}.
	$\hfill\square$
	\subsubsection{Proof of the  sufficiency part of Theorem \ref{QSD}}
	We only need to consider the case $\alpha\in(0,1)$. Suppose that $G_{\mathbf{Q}_{\alpha}}$ is a functional given by (\ref{qsd4}) with a measure $\Delta$ on $(0,\infty)$ satisfying {\rm(a)} and {\rm(b)}. By the above calculations, it is easy to see that $G_{\mathbf{Q}_{\alpha}}$ satisfies (\ref{qsd3}). It remains to show that $G_{\mathbf{Q}_{\alpha}}$ is the p.g.fl. of a probability measure on $\mathfrak{N}^{\circ}(\mathbb{R})$. According to  Maillard \cite{M18}, we can construct a finite point process ($\neq\mathbf{0}$)  whose p.g.fl. is $G_{\mathbf{Q}_{\alpha}}$ as follows. For $x\in(0,\infty)$, suppose that $\mathcal{N}$ is a Poisson random variable with a random parameter drawn according to the measure $x^{-\alpha} \Lambda(\mathrm{d}x)$ and conditioned on non-zero, that is,
		\begin{align*}
			\mathbb{P}(\mathcal{N}=0)=0,\quad \mathbb{P}(\mathcal{N}=k)=\int_{(0,\infty)} \mathrm{e}^{-x} \frac{x^k}{k!} \frac{1}{x^{\alpha}} \Lambda (\mathrm{d}x),\ k\geq 1.
		\end{align*}
		Therefore, the generating function of the positive random variable $\mathcal{N}$ is
		\begin{align*}
			\mathbb{E}\big[ s^{\mathcal{N}} \big] = \int_{(0,\infty)} \big( \mathrm{e}^{(s-1)x}-\mathrm{e}^{-x} \big)\frac{1}{x^{\alpha}} \Lambda (\mathrm{d}x),\quad s\in[0,1].
		\end{align*}
		Suppose that $\{X_i;i\geq1\}$ is a sequence independent and identically distributed point processes with the same p.g.fl. $G_{\mathbf{Q}_{\text{min}}}$. We assume that $\mathcal{N}$ is independent of $\{X_i;i\geq1\}$. Let 
		\begin{align*}
			\bar{X}:=\sum_{i=1}^{\mathcal{N}} X_i.
		\end{align*}
		Then it is easy to check that $\bar{X}$ is a finite point process ($\neq \mathbf{0}$) and its p.g.fl. is $G_{\mathbf{Q}_{\alpha}}$. This finishes the proof of the
		sufficiency part.
	$\hfill\square$ 
	\subsection{Proofs of Propositions \ref{prop6.3}-\ref{prop6.5} }
	\label{sec6.3}
	\noindent\textbf{Proof of Proposition \ref{prop6.3}.} In view of (\ref{6.28}), it suffices to prove that 
	\begin{align}\label{mu1}
		\lim_{n\rightarrow\infty} \frac{\exp\! \big\{q_n^{-1} \!\! \int_{\mathbb{R}} \log G_n[h] (x) \nu(\mathrm{d} x)\big\} \!\! - \! \exp\! \big\{(G_{\mathbf{Q}_{\text{min}}}[h]-1) \nu(\mathbb{R})\big\} } {\nu(\mathbb{R}) \mathrm{e}^{-\nu(\mathbb{R})}} =0
	\end{align}
	and
	\begin{align}\label{mu2}
		\lim_{n\rightarrow\infty} \frac{\big(1-q_n\big)^{q_n^{-1} \nu(\mathbb{R})} \!- \mathrm{e}^{-\nu(\mathbb{R})} }{\nu(\mathbb{R}) \mathrm{e}^{-\nu(\mathbb{R})}}=0
	\end{align}
	hold for $h\in\mathcal{V}_0(\mathbb{R})$ and $\nu\in\mathfrak{M}^{\circ}(\mathbb{R})$. It is no hard to see that $(1-q_n)^{q_n^{-1} \nu(\mathbb{R})}\rightarrow\mathrm{e}^{-\nu(\mathbb{R})}$ as $n\rightarrow\infty$. Then (\ref{mu2}) holds. Note that $\nu(\mathbb{R})<\infty$ and $|q_n^{-1}(G_n[h]-(1-q_n))-G_{\mathbf{Q}_{\text{min}}}[h]|\leq2$ for $n\in\mathbb{N}$. Then in view of (\ref{bc1}), we conclude by dominated convergence that under {\rm(C1)}-{\rm(C3)},
	\begin{align*}
		\bigg| \int_{\mathbb{R}} \Big[ \frac{G_{n}[h](x)-(1-q_{n})}{q_{n}} - G_{\mathbf{Q}_{\text{min}}}[h] \Big] \nu(\mathrm{d} x)  \bigg| \xrightarrow[n\rightarrow\infty]{} 0.
	\end{align*}
	Therefore (\ref{mu1}) follows from (\ref{6.31}) and (\ref{6.32}).
	$\hfill\square$
	
	\
	
	\noindent\textbf{Proof of Proposition \ref{prop6.4}.}
	Suppose that $\widetilde{\gamma}_n$ is a  measure on $\mathfrak{M}(\mathbb{R})$ given by (\ref{6.9}).  
	From the inequality $1-x\leq \mathrm{e}^{-x}$ for $x\in\mathbb{R}$, then taking $h\equiv1$ in (\ref{6.8}) yields
	\begin{align*}
		U_{1,n}(\nu)=\frac{1-(1-q_n)^{q_n^{-1}\nu(\mathbb{R})}}{\nu(\mathbb{R})\mathrm{e}^{-\nu(\mathbb{R})}} \geq \frac{ \mathrm{e}^{\nu(\mathbb{R})}-1 }{\nu(\mathbb{R}) },\quad \nu\in\mathfrak{M}^{\circ}(\mathbb{R}).
	\end{align*}
	On the other hand, by (\ref{6.7}) we get
	\begin{align*}
		\lim_{n\rightarrow\infty} \int_{\mathfrak{M}^{\circ}(\mathbb{R})} U_{1,n}(\nu) \widetilde{\gamma}_n(\mathrm{d} \nu) = G_{\mathbf{Q}_{\alpha}} \big[1\big] <\infty.
	\end{align*}
	Since $\nu \mapsto \frac{ \mathrm{e}^{\nu(\mathbb{R})}-1 }{\nu(\mathbb{R}) }$ can be continuously  extended (with respect to the weak topology) to $\mathfrak{M}(\mathbb{R})$ by setting $\frac{ \mathrm{e}^{\nu(\mathbb{R})}-1 }{\nu(\mathbb{R}) } =1$ for $\nu=\mathbf{0}$. Note that $\widetilde{\gamma}_{n}(\{\mathbf{0}\})=0$, so it is easy to check that
	\begin{align*}
		\int_{\mathfrak{M}(\mathbb{R})}  \frac{ \mathrm{e}^{\nu(\mathbb{R})} - 1 }{\nu(\mathbb{R}) } \widetilde{\gamma}_n(\mathrm{d} \nu) = \int_{\mathfrak{M}^{\circ}(\mathbb{R})}  \frac{ \mathrm{e}^{ \nu(\mathbb{R})} - 1 }{\nu(\mathbb{R}) } \widetilde{\gamma}_n(\mathrm{d} \nu)\leq \int_{\mathfrak{M}^{\circ}(\mathbb{R})}  U_{1,n}(\nu) \widetilde{\gamma}_n(\mathrm{d} \nu) <\infty,\quad n\geq1.
	\end{align*}
	Combining the three preceding equalities implies
	\begin{align}\label{b9}
		\sup_{n\in\mathbb{N}} \int_{\mathfrak{M}(\mathbb{R})}  \frac{ \mathrm{e}^{ \nu(\mathbb{R})} - 1 }{\nu(\mathbb{R}) } \widetilde{\gamma}_n(\mathrm{d} \nu) <\infty.
	\end{align}
	In particular, we have
	\begin{align*}
		\sup_{n\in\mathbb{N}} {\widetilde{\gamma}_n} (\mathfrak{M}(\mathbb{R}))<\infty
	\end{align*}
	and the tightness of $\{\widetilde{\gamma}_n:n\geq0\}$ follows by (\ref{b9}) and \cite[Lemma 3.20]{R87}.
	$\hfill\square$
	
	Recall that $(V(u):u\in\mathbb{G}_{n},n\geq0)$ is a time-homogeneous branching Markov chain with offspring distribution $p$ and one-step transition probability $P$. We next consider a special case of the spinal decomposition theorem for $(V(u):u\in\mathbb{G}_{n},n\geq0)$, which is a key step in proving Proposition \ref{prop6.5}. More specifically, take $r=1$ and the offspring distribution in each generation is $p$, by constructions similar to the conditional reduced branching Markov chain starting from one particle at $x\in\mathbb{R}$ with $r$ spines in \ref{sec4.1}, we obtain a branching Markov chain starting from one particle at $x$ with one spine $w=(w_1,w_2,\cdots)$, which is called a {\it size-biased branching Markov chain}. We denote by $\bar{\mathbb{Q}}_{\delta_x}$ the law of this system. One can refer to \cite{L95,V23,ZS15} for ``size-biased'' branching processes. Let $(\mathcal{F}_n)_{n\in\mathbb{N}}$ be the natural filtration of the original branching Markov chain. From Lemma \ref{lem4.1} we obtain
	\begin{align}\label{6.50}
		\bigg.\frac{\mathrm{d} \bar{\mathbb{Q}}_{\delta_x}}{\mathrm{d} \mathbb{P}_{\delta_x}}\bigg|_{\mathcal{F}_{n}}=\frac{N_n}{m^n}
	\end{align}
	and for any particle $u\in\mathbb{G}_n$,
	\begin{align}\label{6.51}
		\bar{\mathbb{Q}}_{\delta_x}(w_n=u \mid \mathcal{F}_{n}) = \frac{1}{N_n}.
	\end{align}
	
	\noindent\textbf{Proof of Proposition \ref{prop6.5}.}
	For any $h\in\mathcal{V}_0(\mathbb{R})$, by (\ref{6.50}) and (\ref{6.51}) we have
	\begin{align}
		\mathbb{E}_{\delta_x}\big[ \mathrm{e}^{Z_n(\log h)} \mathbf{1}_{\{N_n>0\}} \big] 
		&=m^n \bar{\mathbb{Q}}_{\delta_x} \bigg[ \frac{\mathbf{1}_{\{N_n\geq1\}}}{N_n} \mathrm{e}^{ Z_n(\log h)} \bigg] \nonumber\\
		&= m^n \bar{\mathbb{Q}}_{\delta_x} \bigg[ \frac{\mathbf{1}_{\{V(w_n)\in\mathbb{R}\}}}{N_n} \exp\Big( \sum_{u\in\mathbb{G}_n} \log h (V(u)) \Big) \bigg] \label{6.52}.
	\end{align}
	We now introduce some notations. For every $k\in\mathbb{N}$, let $B(w_{k+1})$ denote the set of brothers of $w_{k+1}$ and $|B(w_{k+1})|$ denote the number of individuals in $B(w_{k+1})$. Recall that $\mathbb{T}$ is a GW-tree rooted at $\varnothing$ defined as in Section \ref{sec2.1}. For $u\in\mathbb{T}$, let $\mathbb{T}^{u}$ be the {\it shifted subtree} of $\mathbb{T}$ generated by $u$ is
	\begin{align*}
		\mathbb{T}^u := \{v: uv\in\mathbb{T}\}
	\end{align*}
	and let $\mathbb{G}^u_k:=\{v\in\mathbb{T}^u:|v|=k\}$ denote the set of individuals of generation $k$ in the subtree $\mathbb{T}^u$. In addition, let
	\begin{align*}
		Z_k^u (\cdot) := \sum_{v\in\mathbb{G}^u_k} \delta_{V(uv)-V(u)} (\cdot)
	\end{align*}
	denote the counting measure of particles of generation $k$ in $\mathbb{T}^u$. By convention, we will write $N_k^u:=Z_k^u(\mathbb{R})$. Let $\rho_u:=V(u)-V(\overleftarrow{u})$, where $\overleftarrow{u}$ denotes the parent of particle $u$. By the spinal decomposition it then follows that
	\begin{align*}
		\sum_{u\in\mathbb{G}_n} \log h(V(u)) = \log h\big(V(w_n)\big) + \sum_{k=0}^{n-1} \sum_{u\in B(w_{k+1})} Z_{n-k-1}^u\big(\log h(V(w_k) + \rho_u + \cdot)\big).
	\end{align*}
	In a similar way, we have
	\begin{align*}
		N_n=1 + \sum_{k=0}^{n-1} \sum_{u\in B(w_{k+1})} N_{n-k-1}^u.
	\end{align*}
	To simplify the formulation, for any $z\in\mathbb{R}$ define
	\begin{align*}
		&Y_{h,n-k}(z) := \sum_{u\in B(w_{k+1})} Z_{n-k-1}^u\big(\log h(z + \rho_u + \cdot)\big),\quad k\in\{0,\cdots,n-1\},\\
		&Y_{h,0}(z) := \log h(z).
	\end{align*}
	Based on these notations, (\ref{6.52}) can be written into
	\begin{align}
		&\ \mathbb{E}_{\delta_x}\big[ \mathrm{e}^{Z_n(\log h)} \mathbf{1}_{\{N_n>0\}} \big] \nonumber\\
		=&\ m^n \bar{\mathbb{Q}}_{\delta_x} \bigg[  \frac{1}{1 + \sum_{k=0}^{n-1} \sum_{u\in B(w_{k+1})} N_{n-k-1}^u}  \exp\Big( \sum_{k=0}^{n} Y_{h,n-k}(V(w_k))  \Big) \bigg] \nonumber\\
		=&\ m^n \bar{\mathbb{Q}}_{\delta_x} \bigg[  \frac{1}{1 + \sum_{k=1}^{n} \sum_{u\in B(w_{n-k+1})} N_{k-1}^u}  \exp\Big( \sum_{k=0}^{n} Y_{h,k}(V(w_{n-k}))  \Big)  \bigg]\label{6.55}.
	\end{align}
	For $K>0$, we consider the event
	\begin{align*}
		\mathcal{I}_n^K:=\{\omega:\sum_{k=K+1}^{n} \sum_{u\in B(w_{n-k+1}(\omega))} N_{k-1}^u(\omega) \geq 1\}.
	\end{align*}
	By the Markov inequality we get
	\begin{align}\label{6.53}
		\bar{\mathbb{Q}}_{\delta_x} (	\mathcal{I}_n^K) = \bar{\mathbb{Q}}_{\delta_x} \Big( \sum_{k=K+1}^{n} \sum_{u\in B(w_{n-k+1})} N_{k-1}^u \geq 1 \Big) \leq \sum_{k=K+1}^{n} \bar{\mathbb{Q}}_{\delta_x} \Big[ \sum_{u\in B(w_{n-k+1})} N_{k-1}^u  \Big].
	\end{align}
	By the construction of $\bar{\mathbb{Q}}_{\delta_x}$, for any $k\in\{1,\cdots,n\}$ we have
	\begin{align}\label{6.54}
		\bar{\mathbb{Q}}_{\delta_x} \Big[ \sum_{u\in B(w_{n-k+1})} N_{k-1}^u  \Big] = \bar{\mathbb{Q}}_{\delta_x} \big[| B(w_{n-k+1}) |\big] \times  \mathbb{E}_{\delta_x} [ N_{k-1} ] = (m^{-1} \sigma^2 -1) m^{k-1},
	\end{align} 
	where $\sigma^2:= \sum_j j^2 p_j<\infty$ under {\rm (C1)}. In the last equality we use the fact that the law of $1+ | B(w_{n-k+1}) |$ under $\bar{\mathbb{Q}}_{\delta_x}$ is $\{\frac{jp_j}{m}\}_{j\in\mathbb{N}}$. Therefore, combining (\ref{6.53}) with (\ref{6.54}) yields
	\begin{align}\label{6.56}
		\bar{\mathbb{Q}}_{\delta_x} (\mathcal{I}_n^K) \leq \frac{m^{-1} \sigma^2 -1}{1-m} m^K = o_K(1),\quad K\rightarrow\infty.
	\end{align}
	Since $\log h\leq 0$ for $h\in\mathcal{V}_0(\mathbb{R})$, this together with (\ref{6.55}) and (\ref{6.56}) implies
	\begin{align}
		&\ \mathbb{E}_{\delta_x}\big[ \mathrm{e}^{Z_n(\log h)} \mathbf{1}_{\{N_n>0\}} \big] \nonumber\\
		=&\ m^n \bigg\{  \bar{\mathbb{Q}}_{\delta_x} \bigg[  \frac{\exp\Big( \sum_{k=0}^{n} Y_{h,k}(V(w_{n-k}))\Big) }{1 + \sum_{k=1}^{n} \sum_{u\in B(w_{n-k+1})} N_{k-1}^u} ; \mathcal{I}_n^K  \bigg] + \bar{\mathbb{Q}}_{\delta_x} \bigg[  \frac{\exp\Big( \sum_{k=0}^{n} Y_{h,k}(V(w_{n-k}))\Big) }{1 + \sum_{k=1}^{n} \sum_{u\in B(w_{n-k+1})} N_{k-1}^u} ; (\mathcal{I}_n^K)^c  \bigg] \bigg\} \nonumber\\
		=&\ m^n \bigg\{  \bar{\mathbb{Q}}_{\delta_x} \bigg[  \frac{1}{1 + \sum_{k=1}^{K} \sum_{u\in B(w_{n-k+1})} N_{k-1}^u}  \exp\Big( \sum_{k=0}^{K} Y_{h,k}\big(V(w_{n-k})\big)  \Big) \bigg] + o_K(1) \bigg\}.\label{6.57}
	\end{align}
	For any $z\in\mathbb{R}$, we write
	\begin{align*}
		\widetilde{Y}_{h,k}(z) := \sum_{u\in B(w_{k+1})} Z_{k-1}^u(\log h(z + \rho_u + \cdot)),\quad k\in\{1,\cdots,n-1\}.
	\end{align*}
	and define the function $F_{K,h}$ by
	$$
	\begin{aligned}
		F_{K, h}: \mathbb{R}^{K+1} & \longrightarrow [0,\infty) \\
		(y,y_1, \cdots, y_{K}) & \longmapsto \bar{\mathbb{Q}}\bigg[\frac{1}{1+\sum_{k=1}^{K} \sum_{u \in B(w_{k+1})} N_{k-1}^u} \exp \Big(\log h(y)+ \sum_{k=1}^{K} \widetilde{Y}_{h, k}\big(V(y_k)\big)\Big)\bigg],
	\end{aligned}
	$$
	where $\bar{\mathbb{Q}}:=\bar{\mathbb{Q}}_{\delta_0}$. Let $\mathcal{G}_n$ be the filtration containing all information about the branching Markov chain and the spine $w$ up to time $n$. By the Markov property we have
	\begin{align}\label{6.58}
		&\ \bar{\mathbb{Q}}_{\delta_x} \bigg[  \frac{1}{1 + \sum_{k=1}^{K} \sum_{u\in B(w_{n-k+1})} N_{k-1}^u}  \exp\Big( \sum_{k=0}^{K} Y_{h,k}\big(V(w_{n-k})\big)  \Big) \bigg] \nonumber\\
		=&\ \bigg. \bar{\mathbb{Q}}_{\delta_x} \bigg[ \bar{\mathbb{Q}}_{\delta_x} \bigg[  \frac{1}{1 + \sum_{k=1}^{K} \sum_{u\in B(w_{n-k+1})} N_{k-1}^u}  \exp\Big( \sum_{k=0}^{K} Y_{h,k}\big(V(w_{n-k})\big)  \Big)  \ \bigg|\ \mathcal{G}_{n-k} \bigg]\bigg] \nonumber\\
		=&\ \bar{\mathbb{Q}}_{\delta_{x}}\! \bigg[ \bar{\mathbb{Q}}_{\delta_{V(w_{n-k})}} \bigg[  \frac{1}{1 \!+\! \sum_{k=1}^{K} \sum_{u\in B(w_{k+1})} N_{k-1}^u}  \exp\!\Big( \sum_{k=0}^{K} \sum_{u\in B(w_{k+1})}\!\! Z_{k-1}^u \big(\log h(V(w_k))+\rho_u+\cdot\big)  \Big) \bigg]  \bigg] \nonumber\\
		=& \int_{\mathbb{R}} \bar{\mathbb{Q}}_{\delta_y} \big[ F_{K,h}\big(y,V(w_1),\cdots,V(w_K)\big) \big] P_{n-k}(x,\mathrm{d}y),
	\end{align}
	where the last equality follows from that $\bar{\mathbb{Q}}_{\delta_x}(V(w_{n-k})\in \mathrm{d}y) = P_{n-k}(x,\mathrm{d}y)$ and the branching and motion mechanisms are independent. Hence, we have 
	\begin{align*}
		&\ \ \ \ \frac{G_n[h](x)-f_n(0)}{1-f_n(0)} 
		= \frac{ \mathbb{E}_{\delta_x}\big[ \mathrm{e}^{Z_n(\log h)} \mathbf{1}_{\{N_n>0\}} \big] }{1-f_n(0)}  \\
		&= \frac{m^n}{1-f_n(0)} \!\times\! \bigg\{\! \int_{\mathbb{R}} \bar{\mathbb{Q}}_{\delta_y} \big[ F_{K,h}\big(y,V(w_1),\cdots,V(w_K)\big) \big] P_{n-k}(x,\mathrm{d}y) \!+\!o_K(1)  \bigg\}\ \text{(by (\ref{6.57}) and (\ref{6.58}))} \\
		&\xrightarrow[n\rightarrow\infty]{\text{uniformly in $x$}} \frac{1}{\varphi(0)} \times \bigg\{ \int_{\mathbb{R}} \bar{\mathbb{Q}}_{\delta_y} \big[ F_{K,h}\big(y,V(w_1),\cdots,V(w_K)\big) \big] \pi(\mathrm{d}y) +o_K(1)  \bigg\}\ \text{(by (\ref{13}) and {\rm(C4)})} \\
		&\xrightarrow[K\rightarrow\infty]{} \frac{1}{\varphi(0)} \int_{\mathbb{R}}  D_{\infty,h} (y) \pi(\mathrm{d}y),
	\end{align*}
	where 
	\begin{align*}
		D_{\infty,h} (y) := \bar{\mathbb{Q}}_{\delta_y} \bigg[ \frac{1}{1+\sum_{k=1}^{\infty} \sum_{u \in B(w_{k+1})} N_{k-1}^u} \exp \Big(\log h(y)+ \sum_{k=1}^{\infty} \widetilde{Y}_{h, k}\big(V(w_k)\big)\Big)  \bigg].
	\end{align*}
	$D_{\infty,h} (y)$ is well-defined as we notice that the sums in its expression are finite almost sure. Indeed, for any $k\in\mathbb{N}_+$, by calculations similar to those in (\ref{6.54}) we get
	\begin{align*}
		\bar{\mathbb{Q}}_{\delta_y} \Big( \sum_{u\in B(w_{k+1})} N_{k-1}^u \geq 1 \Big) \leq \bar{\mathbb{Q}}_{\delta_y} \Big[ \sum_{u\in B(w_{k+1})} N_{k-1}^u \Big] = (m^{-1}\sigma^2-1) m^{k-1}.
	\end{align*}
	Thus by the Borel-Cantelli lemma, there almost surely is only a finite number of integers $k\in\mathbb{N}_+$ such that $\sum_{u\in B(w_{k+1})} N_{k-1}^u \geq 1$. 
	
	On the other hand, by (\ref{bc1}) it is easy to see that under the conditions {\rm(C1)}, {\rm(C2)} and {\rm(C4)},
	\begin{align*}
		\frac{G_n[h](x)-f_n(0)}{1-f_n(0)}  \xrightarrow[n\rightarrow\infty]{} G_{\mathbf{Q}_{\text{min}}}[h],\quad x\in\mathbb{R},\quad h\in\mathcal{V}_0(\mathbb{R}).
	\end{align*} 
	Then the desired result follows from the uniqueness of the limit.
	$\hfill\square$


\begin{thebibliography}{99}
		\baselineskip=14pt
		\bibitem{AH83} Asmussen, S and Hering, H. (1983). {\it Branching Processes}, in: Progress in Probability and Statistics, vol. 3, Birkhäuser Boston, Inc., Boston, MA.
		
		
		
		\bibitem{AN72} Athreya, K. B. and Ney, P. E. (1972). {\it Branching Processes}, in: Die Grundlehren der mathematischen Wissenschaften, vol. Band 196, Springer-Verlag, New York-Heidelberg. 
		
		\bibitem{A13} Aïdékon, E. (2013). Convergence in law of the minimum of a branching random walk. {\it Ann. Probab.} \textbf{41}, 1362-1426.
		
		\bibitem{B18} Bansaye, V. (2019). Ancestral lineages and limit theorems for branching Markov chains in varying environment. {\it J. Theoret. Probab}. \textbf{32}, 249-281.
		
		\bibitem{BH15} Bansaye, V. and Huang, C. (2015). 
		Weak law of large numbers for some Markov chains along non homogeneous genealogies. {\it ESAIM Probab. Stat}. \textbf{19}, 307-326.
		
		
		\bibitem{DV03} Daley, D. J. and Vere-Jones, D. (2003). {\it An introduction to the theory of point processes.} Vol. I. Elementary theory and methods. Second edition. Probability and its Applications (New York). Springer-Verlag, New York.
		
		\bibitem{DV08} Daley, D. J. and Vere-Jones, D. (2008). {\it An introduction to the theory of point processes.} Vol. II. General theory and structure. Second edition. Probability and its Applications (New York). Springer, New York.
		
		\bibitem{D10} Durrett, R. (2010). {\it Probability: Theory and Examples}. Fourth edition, volume 31 of {\it Cambridge Series in Statistical and Probabilistic Mathematics}. Cambridge University Press, Cambridge.
		
		\bibitem{EP90}  Evans, S. N. and Perkins, E. (1990). 
		Measure-valued Markov branching processes conditioned on nonextinction. {\it Israel J. Math}. \textbf{71}, 329-337.
		
		\bibitem{F68} Feller, W. (1968). {\it An introduction to probability theory and its applications}. Vol. I. Third edition. John Wiley \& Sons, Inc., New York-London-Sydney.
		
		
		\bibitem{FP74} Fleischmann, K. and Prehn, U. (1974). Ein Grenzwertsatz f$\ddot{u}$r subkritische Verzweigungsprozesse mit endlich vielen Typen von Teilchen. {\it Math. Nachr}. \textbf{64}, 357-362. 
		
		
		
		\bibitem{FKR77} Fleischmann, K. and Siegmund-Schultze, R. (1977). The structure of reduced critical Galton-Watson processes. {\it Math. Nachr}. \textbf{79}, 233-241.
		
		\bibitem{FS78} Fleischmann, K. and Siegmund-Schultze, R. (1978). An invariance principle for reduced family trees of critical spatially homogeneous branching processes. {\it Serdica 4}. \textbf{2-3}, 111-134.
		
		\bibitem{FV99} Fleischmann, K., Vatutin, V. A. (1999). Reduced subcritical Galton-Watson processes in a random environment. {\it Adv. in Appl. Probab.} \textbf{31}, 88-111.
		
		
		\bibitem{G99} Geiger, J. (1999). Elementary new proofs of classical limit theorems for Galton-Watson processes. {\it J. Appl. Probab}. \textbf{36}, 301-309.
		
		\bibitem{GR14} Grafakos, L. (2014). {\it Classical Fourier analysis}. Third edition. Springer, New York.
		
		\bibitem{H77}  Hering, H. (1977). Minimal moment conditions in the limit theory for general Markov branching processes. {\it Ann. Inst. H. Poincaré Sect. B (N.S.)} \textbf{13}, 299-319.
		
		
		\bibitem{H22} Harris, S. C., Horton, E., Kyprianou, A. E. and Wang, M. (2022). 
		Yaglom limit for critical nonlocal branching Markov processes. {\it Ann. Probab}. \textbf{50}, 2373-2408.
		
		\bibitem{HR17} Harris, S. C. and Roberts, M. I. (2017). 
		The many-to-few lemma and multiple spines. {\it Ann. Inst. Henri Poincaré Probab. Stat}. \textbf{53}, 226-242.
		
		\bibitem{HS67} Heathcote, C. R., Seneta, E. and Vere-Jones, D. (1967). A refinement of two theorems in the theory of branching processes. {\it Teor. Verojatnost. i Primenen}. \textbf{12}, 341-346.
		
		\bibitem{HL23} Hong, W. M. and Liang, S. L. (2024).  Conditional central limit theorem for critical branching
		random walk. {\it ALEA, Lat. Am. J. Probab. Math. Stat.}  \textbf{21}, 555-574.
		
		\bibitem{HY23} Hong, W. M. and Yao, D. (2023).  Conditional central limit theorem for subcritical branching random walk. {\it ALEA, Lat. Am. J. Probab. Math. Stat}. \textbf{20},  1411-1432.
		
		\bibitem{H80}  Hoppe, F. M. (1980). On a Schröder equation arising in branching processes. {\it Aequationes Math.} \textbf{20}, 33-37.
		
		\bibitem{HS78} Hoppe, F. M. and Seneta, E. (1978). Analytical methods for discrete branching processes. {\it Adv. Probab. Related Topics}. \textbf{5}, 219-261.
		
		\bibitem{HS09} Hu, Y. and Shi, Z. (2009). Minimal position and critical martingale convergence in branching random walks, and directed polymers on disordered trees. {\it Ann. Probab.} \textbf{37}, 742-789.
		
		
		
		
		
		
		
		
		
		
		
		\bibitem{J67} Joffe, A. (1967). On the Galton-Waston branching process with mean less than one. {\it Ann. Math. Statist}. \textbf{38}, 264-266.
		
		\bibitem{K17} Kallenberg, O. (2017). {\it Random measures, theory and applications}, in: Probability Theory and Stochastic Modelling, vol. 77, Springer, Cham.
		
		\bibitem{L07} Lambert, A. (2007). Quasi-stationary distributions and the continuous-state branching process conditioned to be never extinct. {\it Electron. J. Probab.} \textbf{12}, 420-446.
		
		\bibitem{L72} Lebedev, N. N. (1972). {\it Special functions and their applications}, Dover Publications, Inc., New York.
		
		\bibitem{L11} Li, Z. (2011). {\it Measure-valued branching Markov processes}, in: Probability and its Applications (New York), Springer, Heidelberg.
		
		\bibitem{LR21} Liu, R., Ren, Y.-X., Song, R. and Sun, Z. (2021). Quasi-stationary distributions for subcritical superprocesses. {\it Stochastic Process. Appl.} \textbf{132}, 108-134. 
		
		\bibitem{L95} Lyons, R., Pemantle, R. and Peres, Y. (1995). Conceptual proofs of $L\log L$ criteria for mean behavior of branching processes. {\it Ann. Probab.} \textbf{23}, 1125-1138.
		
		
		
		
		
		\bibitem{M18} Maillard, P. (2018). The  $\lambda$-invariant measures of subcritical Bienaymé-Galton-Watson processes. {\it Bernoulli} \textbf{24}, 297-315. 
		
		\bibitem{MV12} Méléard, S. and Villemonais, D. (2012). Quasi-stationary distributions and population processes. {\it Probab. Surv.} \textbf{9}, 340-410.
		
		\bibitem{R87} Resnick, S. I. (1987). {\it Extreme Values, Regular Variation, and Point Processes.} Applied Probability, vol 4. New York: Spring.
		
		
		
		\bibitem{KS67}  Spitzer, F. (1967). Two explicit Martin boundary constructions. In {\it Symposium on Probability Methods in Analysis (Loutraki, 1966). Lecture Notes in Math.} \textbf{31}, 296-298. Springer, Berlin-New York.
		
		\bibitem{ZS15} Shi, Z. (2015). {\it Branching Random Walks}, volume 2151 of {\it Lecture Notes in Mathematics}. Springer, Cham. 
		
		\bibitem{V23} Valentin, R. (2023). Invariant measures of critical branching random walks in high dimension. {\it Electron. J. Probab.} \textbf{28}, 1-38.
		
		\bibitem{Y47} Yaglom, A. M. (1947). Certain limit theorems of the theory of branching random processes. {\it Dokl. Akad. Nauk SSSR(N.S.)} \textbf{56}, 795-798.
		
		
		
	\end{thebibliography}
\end{document}